\documentclass[11pt]{article}
\usepackage[utf8]{inputenc}
\usepackage{amsmath,amssymb,epsfig,bbm}
\usepackage{stmaryrd,mathabx}
\usepackage{comment}
\usepackage{color}
\usepackage[T1]{fontenc}

\usepackage[textsize=scriptsize,textwidth=2.1cm]{todonotes}
\usepackage{enumitem}
\usepackage{varwidth}
\setlist{nolistsep}
\usepackage{hyperref}

\newtheorem{defi}{Definition}[section]
\newtheorem{prop}[defi]{Proposition}
\newtheorem{theo}[defi]{Theorem}
\newtheorem{conj}[defi]{Conjecture}
\newtheorem{lemm}[defi]{Lemma}
\newtheorem{coro}[defi]{Corollary}
\newtheorem{rema}[defi]{Remark}
\newtheorem{exem}[defi]{Example}
\newtheorem{exems}[defi]{Examples}

\newcommand{\bdefi}{\begin{defi}}
\newcommand{\edefi}{\end{defi}}
\newcommand{\bprop}{\begin{prop}}
\newcommand{\eprop}{\end{prop}}
\newcommand{\btheo}{\begin{theo}}
\newcommand{\etheo}{\end{theo}}
\newcommand{\blemm}{\begin{lemm}}
\newcommand{\brema}{\begin{rema}}
\newcommand{\erema}{\end{rema}}
\newcommand{\bexer}{\begin{exem}}
\newcommand{\eexer}{\end{exem}}
\newcommand{\bexems}{\begin{exems}}
\newcommand{\eexems}{\end{exems}}
\newcommand{\bconj}{\begin{conj}}
\newcommand{\econj}{\end{conj}}
\newcommand{\elemm}{\end{lemm}}
\newcommand{\bcoro}{\begin{coro}}
\newcommand{\ecoro}{\end{coro}}
\newcommand{\dem}{\noindent{\bf Proof. }}
\newcommand{\rem}{\noindent{\bf Remark. }}

\usepackage{mathrsfs}
\renewcommand\mathcal{\mathscr}

\newcommand{\G}{{\cal G}}

\newcommand{\I}{{\cal I}}

\newcommand{\M}{{\cal M}}
\newcommand{\N}{{\cal N}}
\newcommand{\OOO}{{\cal O}}
\renewcommand{\P}{{\cal P}}

\newcommand{\Scal}{{\cal S}}

\newcommand{\X}{{\cal X}}

\newcommand{\maths}[1]{{\mathbb #1}}  

\renewcommand{\AA}{\maths{A}}

\newcommand{\CC}{\maths{C}}

\newcommand{\FF}{\maths{F}}
\newcommand{\HH}{\maths{H}}

\newcommand{\NN}{\maths{N}}

\newcommand{\PP}{\maths{P}}
\newcommand{\QQ}{\maths{Q}}
\newcommand{\RR}{\maths{R}}

\newcommand{\TT}{\maths{T}}

\newcommand{\ZZ}{\maths{Z}}

\newcommand{\aaa}{{\mathfrak a}}

\newcommand{\ddd}{{\mathfrak d}}

\renewcommand{\ggg}{{\mathfrak g}}

\renewcommand{\lll}{{\mathfrak l}}

\newcommand{\ppp}{{\mathfrak p}}

\newcommand{\sss}{{\mathfrak s}}
\newcommand{\uuu}{{\mathfrak u}}

\newcommand{\weakstar}{\overset{*}\rightharpoonup}
\newcommand{\ra}{\rightarrow}
\newcommand{\bs}{\backslash}

\newcommand{\ov}[1]{{\overline #1}} 
\newcommand{\wt}[1]{{\widetilde{#1}}}
\newcommand{\wh}[1]{{\widehat{#1}}}

\newcommand{\ga}{\gamma}
\newcommand{\Ga}{\Gamma}
\newcommand{\ssm}{\!\smallsetminus\!}
\newcommand{\mb}[1]{{\mathbf{#1}}}

\usepackage{mathtools}

\DeclarePairedDelimiter\intbra{\llbracket}{\rrbracket}

\newcommand{\cqfd}{\hfill$\Box$}

\newcommand{\Ad}{\operatorname{Ad}}

\newcommand{\bigO}{\operatorname{O}}
\newcommand{\card}{{\operatorname{Card}}}

\newcommand{\covol}{\operatorname{covol}}

\newcommand{\dbs}{\backslash\!\!\backslash}

\newcommand{\diag}{{\operatorname{diag}}}

\newcommand{\dvol}{\;d\operatorname{vol}}

\newcommand{\id}{\operatorname{id}}

\newcommand{\Leb}{\operatorname{Leb}}
\renewcommand{\ln}{\operatorname{ln}}

\newcommand{\loz}{\lozenge}

\renewcommand{\Re}{{\operatorname{Re}}}

\newcommand{\sys}{\operatorname{sys}}

\newcommand{\vol}{\operatorname{vol}}

\newcommand{\PSL}{\operatorname{PSL}}
\newcommand{\SL}{\operatorname{SL}}
\newcommand{\GL}{\operatorname{GL}}

\newcommand{\PGL}{\operatorname{PGL}}

\newcommand{\n}{\operatorname{\tt n}}

\newcommand{\redn}{\operatorname{\tt N}}

\newcounter{fig}

\title{Equidistribution of divergent diagonal orbits\\
in positive characteristic}
\author{Nguyen-Thi Dang  \and Fr\'ed\'eric Paulin \and Rafael Sayous} 
\date{\today}

\begin{document}
\bibliographystyle{../alphanum}
\maketitle
\begin{abstract}
Given a local field $\wh K$ with positive characteristic, we study the
dynamics of the diagonal subgroup of the linear group $\GL_n(\wh K)$
on homogeneous spaces of discrete lattices in ${\wh K}^{\,n}$. We
first give a function field version of results by Margulis and
Tomanov-Weiss, characterizing the divergent diagonal orbits.
%as the rational ones.
When $n=2$, we relate the divergent diagonal orbits
with the divergent orbits of the geodesic flow in the modular quotient
of the Bruhat-Tits tree of $\PGL_2(\wh K)$. Using the (high) entropy
method by Einsiedler-Lindentraus et al, we then give a function field
version of a result of David-Shapira on the equidistribution of a
natural family of these divergent diagonal orbits, with height given
by a new notion of discriminant of the orbits.
\end{abstract}

\section{Introduction}
\label{sect:intro}

Equidistribution problems of periodic orbits have been widely studied
in many different settings. In hyperbolic dynamical systems, in
particular for closed orbits of geodesic flows in negative curvature,
see for instance \cite{Margulis69, Bowen72b} and many others,
including \cite[\S 9.3]{PauPolSha15} (see references therein). In
homogeneous dynamics for diagonalisable group actions (sometimes in an
arithmetic framework), see for instance \cite{EinLinMicVen11} and many
others, including \cite{DanLi22} (see references therein). See also
\cite{Shapira17,SolYif23,KemPauSha17} (this last one also over
function fields) for possible chaotic behaviors of weak-star limits of
homogeneous measures on periodic orbits, including surprising loss of
mass phenomena.

Much less studied has been the problem of equidistribution of
divergent orbits, as they require noncompact phase spaces and a
specific study of equidistribution of (locally finite) infinite
measures. See \cite{ParPauSay25} for geodesic flows in variable
curvature, as well as \cite{DavSha18,DavSha20} in homogeneous
dynamics.
%These works of David-Shapira will be the main guiding line
%for the second part of our work.
Considering homogeneous dynamics over various local field is important
and fruitful. The first purpose of this paper is to extend to local
fields in positive characteristic works of Margulis, Tomanov-Weiss
\cite {TomWei03} and Tomanov \cite {Tomanov07} on the characterisation
of divergent orbits. The second purpose is to extend David-Shapira
\cite{DavSha18,DavSha20} results on their equidistribution (with the
challenges required for such an extension). The study of divergent
orbits in homogeneous dynamics, through the Dani correspondence, has
strong ties with Diophantine approximation problem, see for instance
\cite{KleWei05,CheChe16,KadKleLinMar17,DasFisSimUrb24,AraKim25,
  BanKimLim25}, these last two references also over function fields.

Throughout this paper, referring to \cite{Goss98,Rosen02} and
Subsection \ref{subsec:functionfield} for definitions and complements,
we fix a function field $K$ of genus $\ggg$ over a finite field
$\FF_q$ of order $q$, a valuation $v$ of $K$ and a uniformiser $\pi_v$
of $v$. We denote by $K_v$ the completion of $K$ for $v$, by $\OOO_v$
its valuation ring, by $q_v$ the order of its residual field, by
$|\cdot|= q_v^{-v(\,\cdot\,)}$ its (normalized) absolute value, by
$R_v$ the affine function ring associated with $v$, by $\zeta_v$ the
Dedekind zeta function of $R_v$, and by $\varphi_v$ the Euler function
of $R_v$.

We fix $n\in \NN\ssm\{0,1\}$. The unimodular group \mbox{$\GL^1_n(K_v)
  = \{g\in\GL_n(K_v):|\det g|=1\}$} is endowed with the Haar measure
giving mass $1$ to the maximal compact subgroup $\GL_n(\OOO_v)$. We
denote by $\X_1$ the $\GL^1_n(K_v)$-homogeneous space of
$R_v$-lattices in $K_v^{\;n}$ with normalized covolume $1$ (identified
with $\GL^1_n(K_v)/\GL_n(R_v)$ when pointed at the standard
$R_v$-lattice $R_v^{\;n}$). We endow $\X_1$ with the induced
$\GL^1_n(K_v)$-invariant measure, that we denote by ${\tt m}_{\X_1}$
and which is finite.

We denote by $A_1$ the diagonal subgroup of $\GL^1_n(K_v)$, and we
normalize its Haar measure to give mass $1$ to its maximal compact
subgroup $A_1\cap\GL_n(\OOO_v)$. The diagonal orbit $A_1x$ of an
element $x\in\X_1$ is said to be {\it divergent} if the orbital map
$a\mapsto ax$ from $A_1$ to $\X_1$ is proper. The {\it homogeneous
  measure} on $A_1x$, that we denote by $\overline{\mu}_x$, is then
the (locally finite) pushforward measure by this orbital map of the
Haar measure of $A_1$.

The first main result of this paper (see Corollary
\ref{coro:classdivunimod} for a more general result and Theorem
\ref{theo:defidiv} for the analog result for the projective linear
group $\PGL_n(K_v)\,$) is an algebraic characterisation of the
divergent orbits, saying that they are the ``rational'' ones, that is,
they come from a rational point (in $\GL_n^1(K)$) of $\GL^1_n(K_v)$,
up to the action of an element of $A_1$.

\btheo\label{theo:mainuintro1} Let $x\in\X_1$. The diagonal orbit
$A_1x$ is divergent if and only if there exists $g\in A_1\GL_n^1(K)$
such that $x=g\,R_v^{\,n}$.
\etheo

This result has a long history. In the real field case with $n=2$, an
orbit of the diagonal subgroup of $\PSL_2(\RR)$ on $\PSL_2(\RR)/
\PSL_2(\ZZ)$ is well known to be divergent if and only if it
corresponds to a modular group orbit of a geodesic line in the upper
halfspace model $\HH^2_\RR$ of the real hyperbolic plane both of whose
endpoints are rational points (that is, are in $\PP_1(\QQ)$) of the
circle at infinity $\PP_1(\RR)$ of $\HH^2_\RR$. In the function field
case with $n=2$, the corresponding result is also well known: The
quotient of the Bruhat-Tits tree of $\PGL_2(K_v)$ by $\PGL_2(R_v)$
replaces the quotient of the real hyperbolic plane by $\PSL_2(\ZZ)$,
and the set of rational points at infinity is $\PP_1(K)$ in
$\PP_1(K_v)$ (see for instance \cite{Serre83} and Section
\ref{subsec:zigzag}). The real field case for any integer $n$ is due
to Margulis (see \cite[Appendix]{TomWei03}). It has been extended to
all reductive algebraic groups over $\QQ$ in
\cite[Theo.~1.1]{TomWei03}, and to the $S$-adic case over number
fields in \cite{Tomanov21}. The case of divergent orbits of proper
subgroups of the full maximal $\QQ$-torus subgroup has been studied by
\cite{Taman22}, with surprising differences.

\medskip
For every $k\in\NN\ssm\{0\}$, we identify each element of $K_v^{\;k}$
with the column matrix of its coordinates in the canonical basis of
$K_v^{\;k}$. For every element $\mb{t}\in K_v^{\;n-1}$, we define
\begin{equation}\label{eq:defuuut}
\uuu_{\mb{t}}=\begin{pmatrix} 1 & 0\\\mb{t} & I_{n-1}
\end{pmatrix}\in\SL_n(K_v)\,.
\end{equation}
Note that if $\mb{t}\in K^{n-1}$, then $\uuu_{\mb{t}} \,R_v^{\;n}$ is
an $R_v$-lattice with normalized covolume $1$ whose diagonal orbit is
divergent by Theorem \ref{theo:mainuintro1}.

The second main result of this paper (see Corollary
\ref{coro:mainsgras} for a more general result) is the following
equidistribution result in $\X_1$ for natural families of divergent
diagonal orbits in $\X_1$. We emphazise the fact that the measures
that equidistribute are infinite measures. But for the weak-star
convergence, sequences of locally finite infinite measures may indeed
converge to a finite measure. Such is not the case for the narrow
convergence.

\btheo\label{theo:mainuintro2}
Let $c_{K,n}=\frac{(n-1)!\,\prod_{i=1}^{n-1} \zeta_v(-i)}
{q_v\,(q-1)\,\prod_{i=2}^{n-1}(q_v^{\;i}-1) }$.
%%
%\todo{\tiny jolie constante arithmétique :-)}
%%
For every nonzero $s\in R_v$, let us define $\Lambda(s)=
\{(\frac{r_2}{s},\ldots,\frac{r_n}{s}): r_2,\ldots, r_n\in
R_v,\;\forall j\in\{2,\ldots,n\},\;r_jR_v+sR_v=R_v\} \mod
R_v^{\;n-1}$.  Assume that $R_v$ is principal. For the weak-star
convergence of Radon measures on the locally compact space $\X_1$, we
have
%%
%\todo{\tiny voir s'il faut se restreindre à $R_v$ principal}
%%
%%
%\todo{\tiny purement spéculatif...}
%%
\[
\lim_{|s|\ra+\infty,\;|s|\in q_v^{\;n\ZZ}}\;\;
\frac{c_{K,n}}{(\varphi_v(s)\,\log_{q_v}|s|\,)^{n-1}}\;
\sum_{\mb{t}\in \Lambda(s)}\; \ov\mu_{\uuu_{\mb{t}} R_v^{\;n}}\;=\; {\tt m}_{\X_1}\,.
\]
\etheo

Let us discuss the scope of this result. We believe that the principal
assumption on $R_v$ may not be necessary, since we are only using it
to prove the non-escape of mass property in Section
\ref{sec:noescapmass}, and an approach along the lines of
\cite{DavKimMorSha25} could allow its removal. Over the real field,
this result is due to \cite{DavSha18} in dimension $n=2$ and to
\cite{DavSha20} in general. Starting from Section
\ref{sec:entropytoequidistrib}, we will follow their scheme of
proof.
%But even though we believe that David-Shapira could also have
%obtained these improvements,
Our result has two new aspects, besides the fact that the algebraic
properties of the ring $R_v$ are much more involved than the ones of
$\ZZ$. Firstly, we obtain an explicitely renormalized weak-star
convergence, and not only a projective convergence of the
measures. Secondly, we cover a larger set of types of divergent
orbits, as we now explain.  It follows from Theorem
\ref{theo:mainuintro1} that an $A_1$-orbit $\Theta$ is divergent if
and only if it contains an element containing a sub-$R_v$-lattice
$\Lambda$ of $R_v^{\;n}$. We define the {\it type} of a
sub-$R_v$-lattice $\Lambda$ of $R_v^{\;n}$ as the isomorphism class of
the torsion $R_v$-module $R_v^{\;n}/\Lambda$, and the {\it type} of
$\Theta$ as the finite set of types of the sub-$R_v$-lattices of
$R_v^{\;n}$ with minimal covolume contained in the elements of
$\Theta$.  For instance, for every $\mb{t}\in \Lambda(s)$, the type of
the $A_1$-orbit $A_1\uuu_{\mb{t}} R_v^{\;n}$ is reduced to
$\{(R_v/sR_v)^{n-1}\}$ (see Proposition \ref{prop:descripxt}
\eqref{item3:descripxt}).  We choose this type in this Introduction
for simplicity, but we refer to Corollary \ref{coro:mainsgras} for a
generalisation.
%%
%\todo{\tiny Voir si nous pouvons faire...}
%%

\medskip
The techniques of the second part of this paper, that we now present,
rely in particular on the (high) entropy method in homogeneous
dynamics (see for instance \cite{EinLin10}). Let
\begin{equation}\label{eq:defaaa}
\aaa=\begin{pmatrix} \pi_v^{\;n-1} &
  0\\0 & \pi_v^{\;-1}I_{n-1} \end{pmatrix}\in\SL_n(K_v)\,.
\end{equation}
Let $U^-=\{\uuu_{\mb{t}}:\mb{t}\in K_v^{\;n-1}\}$, which is the
unipotent radical of the parabolic subgroup of $\SL_n(K_v)$ fixing the
hyperplane $\{0\}\times K_v^{\;n-1}$ of $K_v^{\;n}$. Note that for all
$k\in\ZZ$ and $\mb{t}\in K_v^{\;n-1}$, we have $\aaa^k\uuu_{\mb{t}}
\aaa^{-k} =\uuu_{\pi_v^{\;-nk}\mb{t}}$, so that $U^-$ is contained in
(and actually equal to) the unstable horospherical group of the
one-parameter diagonal group $(\aaa^k)_{k\in\ZZ}$.

Using Mahler's criterion that a sequence $(x_n)_{n\in\NN}$ in $\X_1$
goes out of every compact subset if and only if the systole of $x_n$
goes to $0$, the first step (see Section \ref{sec:descripdivdiagorb})
is to describe a canonical ``compact core'' $C_x$ of a divergent
diagonal $A_1$-orbit $x$ by trimming out the parts were the systole of
the elements of $x$ are small. When $n=2$, this correspond to removing
the first and last intersections with a cuspidal ray of the
corresponding geodesic line in the quotient graph of groups by
$\PGL_2(R_v)$ of the Bruhat-Tits tree of $\PGL_2(K_v)$ (see
\cite{Serre83,BroParPau19} for background and Subsection
\ref{subsec:zigzag}). Let us denote by $\ov\nu_x$ the restriction of
the homogeneous measure $\ov\mu_x$ to $C_x$ normalized to be a
probability measure, and by $\ov\mu_s=\sum_{\mb{t}\in \Lambda_s}\ov
\mu_{\uuu_{\mb{t}} R_v^{\;n}}$ and $\ov\nu_s=\frac{1}{\card \Lambda_s}
\sum_{\mb{t}\in \Lambda_s} \ov\nu_{\uuu_{\mb{t}} R_v^{\;n}}$.  It will
follow from Subsection \ref{subsec:massbehavedivdiaorb} (with the help
of computations done in Subsection \ref{subsec:typdiscdivdiaorb}) that
the measures $\frac{c_{K,n}}{(\varphi_v(s)\,\log_{q_v}|s|\,)^{n-1}}
\;\ov\mu_s$ have as $|s|\ra+\infty$ the same weak-star asymptotic
properties as the probability measures $\ov\nu_s$. Furthermore, we
prove that the measures $\ov\nu_s$ on $\X_1$ are averages over a
compact subgroup $C^1_n$ of $\GL^1_n(K_v)$ of natural measures $\nu_s$
on $\SL_n(K_v)/\SL_n(R_v)$.

The second step (see Section \ref{sec:noescapmass}) is to prove that
the measures $\nu_s$ on $\SL_n(K_v)/\SL_n(R_v)$ as $|s|\ra+\infty$ do
not suffer any loss of mass, that is, any weak-star accumulation point
$\nu$ of $\nu_s$ as $|s|\ra+\infty$ is also a probability measure. We
did not try to write our equidistribution result replacing the set
$\Lambda_s$ by a logarithmic full proportion of it, as it is done in
\cite{DavSha18,DavSha20}.  This extension requires a non escape of
mass assumption, that has been lifted in \cite{DavKimMorSha25} when
$n=2$ in the real case.

The third step (see Section \ref{sec:lowboundentrop}) is to prove that
the entropy $h_\nu(\aaa)$, which is well defined by the second step,
of any weak-star accumulation point $\nu$ for the diagonal
transformation $\aaa$ on $\SL_n(K_v)/\SL_n(R_v)$ is equal to the
maximal entropy of this transformation. This requires, as in
\cite{EinLinMicVen12}, a construction of high entropy partitions for
$\aaa$, that are build using dynamical neighborhoods for the action of
$\aaa$ on its unstable horospherical group $U^-$.

The last step (see Section \ref{sec:entropytoequidistrib}) is to apply
Einsiedler-Lindenstraus \cite{EinLin10} uniqueness of the probability
measure of maximal entropy on $\SL_n(K_v)/\SL_n(R_v)$ for $\aaa$,
which is the measure ${\tt m}_{\X_1}$ renormalized to be a probability
measure, and to average back on the above-mentionned compact subgroup
$C^1_n$ in order to prove Theorem \ref{theo:mainuintro2}.

Obtaining an error term in Theorem \ref{theo:mainuintro2} would
require an effective version of the uniqueness of measures of maximal
entropy for diagonal actions in positive characteristic, and would
constitute another project. We believe that our results could be extended
to the $S$-adic case (working with a nonempty finite set of places $S$
instead of just one $v$) or to the adelic setting (for the nonuniform
lattice $\PGL_n(K)$ of $\PGL_n(\AA_K)$, where $\AA_K$ is the adèle
ring of $K$).
%We plan an extension of this paper to general semi-simple algebraic
%groups over $K$.
%%
%\todo{\tiny à voir}
%%

\medskip
\noindent{\small {\it Acknowledgements: } This research was supported
  by the French-Finnish CNRS IEA PaCap. We thank Taehyeong Kim for
  mentionning the paper \cite{DavKimMorSha25} and for his comments
  about the loss of mass problem.}
%%
%\todo{\tiny autres ?}
%%

\section{Background material}
\label{sec:background}

For all $r,r'\in\ZZ$ with $r\leq r'$, we denote $\llbracket
r,r'\rrbracket=[r,r']\cap\ZZ$.

\subsection{Function fields over finite fields}
\label{subsec:functionfield}

For the following notions and complements, we refer to \cite{Goss98,
  Rosen02}, as well as to \cite[\S 14.2]{BroParPau19} whose notation
we follow. Let $\FF_q$ be a finite field of order a positive
power $q$ of a prime $p$.  Let $K$ be a (global) function field over
$\FF_q$ of genus $\ggg$, that is, the function field of a
geometrically connected smooth projective curve ${\mathfrak C}$ of
genus $\ggg$ defined over $\FF_q$. We denote by $h_K$ the number of
divisor classes of degree $0$ on ${\mathfrak C}$. Let $v$ be a
(normalised discrete) valuation of $K$, let $K_v$ be the associated
completion of $K$, let $\OOO_v=\{x\in K_v:v(x)\geq 0\}$ be its
valuation ring, let $\pi_v\in K$ with $v(\pi_v)=1$ be a uniformiser of
$v$, let $q_v$ be the order of the residual field $\FF_{q_v}= \OOO_v/
\pi_v\OOO_v$ (that we identify with its lift in $\OOO_v$), and let
$|\cdot|= q_v^{-v(\,\cdot\,)}$ be the (normalized) absolute value
associated with $v$. We denote by $\deg v\in\NN \ssm\{0\}$ the degree
of the closed point of ${\mathfrak C}$ corresponding to $v$, so that
$q_v=q^{\deg v}$. Let $R_v$ be the affine function ring associated
with $v$, that is, the affine algebra of the curve ${\mathfrak C}$
minus its closed point corresponding to $v$. Recall that $R_v$ is a
Dedekind ring whose field of fractions is $K$. The class number $h_v$
of the Dedekind ring $R_v$ is $h_v=(\deg v)\,h_K$ by
\cite[Coro.~4.1.3]{Goss98}.  In particular $R_v$ is principal if and
only if $h_K=1$ and $\deg v=1$, which occurs in positive genus for
exactly $4$ isomorphisms classes of function fields $K$ (one for each
$(\ggg,q)=(1,2),(1,3),(1,4),(2,2)$) by \cite[Theo.~1.1]{MerSti15} and
\cite[Theo.~2]{MadQue72}. Note that $\OOO_v^\times = \{x\in K_v:
|x|=1\}$ and (see for instance \cite[Eq.~(14.2) and
  (14.3)]{BroParPau19}
\begin{equation}\label{eq:inversiRv}
  R_v\cap\OOO_v=\FF_q\quad\text{and}\quad
  R_v^\times=\FF_q^\times\subset \OOO_v^\times \,.
\end{equation}

The simplest example, used in Section \ref{subsec:zigzag}, is given by
the field $K=\FF_q(Y)$ of rational fractions over $\FF_q$ with one
indeterminate $Y$, with genus $\ggg=0$, endowed with the valuation at
infinity $v$ with $\deg v=1$ defined for all $P,Q\in \FF _q[Y]$ with
$Q\neq 0$ by $v(\frac{P}{Q}) \deg Q-\deg P$. Then $K_v=
\FF_q((Y^{-1}))$ is the field of formal Laurent series in $Y^{-1}$
over $\FF_q$, $\OOO_v=\FF_q[[Y^{-1}]]$ is the local ring of formal
power series in $Y^{-1}$ over $\FF_q$, $\pi_v=Y^{-1}$, $q_v=q$, and
$R_v=\FF_q[Y]$ is the ring of polynomials in $Y$ over $\FF_q$.

\medskip
Let $\I_v^+$ be the semigroup of nonzero ideals of the ring $R_v$. As
usual, $\ppp$ ranges throughout the text over prime ideals in $\I_v^+$
and $\redn (I)=[R_v:I]\in q^\NN$ is the absolute norm of $I\in
\I_v^+$. Let $\redn (s)=\redn (sR_v)$ for every $s\in R_v\ssm \{0\}$,
and note that $\redn (s)=|s|$. For all $r,s\in R_v$, we write as usual
$(r,s)=1$ if $r$ and $s$ are coprime, that is, satisfy
$r\,R_v+s\,R_v=R_v$.

We denote by $\mu_v:\I^+_v\ra \ZZ$ the Möbius function of $R_v$, so
that $\mu_v(I)=0$ if $I$ has a squared prime factor, and otherwise
$\mu_v(I)=(-1)^k$ where $k$ is the number of prime factors of $I$.

We denote by $\varphi_v: \I_v^+ \ra \NN$ the Euler function of $R_v$,
defined by
\begin{equation}\label{eq:defiEulerfunct}
\varphi_v(I)= \card \;\big(R_v/I)^\times= \redn (I)
\prod_{\ppp\,\mid\, I}\big(1-\frac{1}{\redn (\ppp)}\big)
=\redn (I)\sum_{I'\in \I_v^+,\;I'\,\mid\, I}
\frac{\mu_v(I')}{\redn (I')}\,,
\end{equation}
and $\varphi_v(s)=\varphi_v(sR_v)$ for every $s\in R_v$.

The Dedekind zeta function of the Dedeking ring $R_v$ is (see, for
instance, \cite[\S 7.8]{Goss98}) the map $\zeta_v:\{z\in\CC:\Re
\;z>1\}\ra\CC$ defined by
\[
\zeta_v :z\mapsto \sum_{I\in\I_v^+} \frac{1}{\redn (I)^z}\,.
\]
By for instance \cite[page 219, line 2]{Goss98} or \cite[page 244,
  Eq.~(1)]{Rosen02}, it is related to the zeta function $\zeta_K$ of
the field $K$ (which is an Eulerian product over all closed points of
${\mathfrak C}$, including the one corresponding to $v$) by the
formula
\begin{equation}\label{eq:relatzetaKzetav}
  \zeta_K(z)=\frac{1}{1-q_v^{-z}}\;\zeta_v(z)\;.
\end{equation}
By for instance \cite[Theo.~5.9]{Rosen02}, $\zeta_K$ has an analytic
continuation on $\CC \ssm\{0, 1\}$ with simple poles at $z = 0$ and $z
= 1$ (it is actually a rational function of $q^{-z}$). Hence the value
$\zeta_v(-k)$ for every $k\in\NN\ssm\{0\}$ is well defined. We recall the
following counting result. 

\blemm \label{lem:gauscounting}
As $t\ra+\infty$, we have
\[
\card\,\{I\in \I_v^+:\redn(I)\leq t\}=\frac{h_K\; q^{2-\ggg}
  \,(q_v-1)}{(q-1)^2\,q_v}\;q^{\lfloor\log_qt\rfloor}+\bigO(1)\,.
\]
%We have $\card\,\{s\in R_v:|s|\leq t\}= q^{1-\ggg}
% \;q_v^{\lfloor \log_{q_v}t\rfloor}+\bigO(1)$.
%%%
%\todo{\tiny à vérifier}
%%%
\elemm

\dem We give a proof for completeness. Let $z\in\CC$ with
$\Re\;z>1$. For every $n\in\NN$, let $c_n=\card\{I\in \I_v^+ :\redn(I)
= q^n\}$. Since $\redn(I)\in q^\NN$ for every $I\in \I_v^+$, we have
$\zeta_v(z)=\sum_{n=0}^\infty\, c_n\,q^{-nz}$. By for instance
\cite[end of page 52]{Rosen02}, we have
$\zeta_K(z)=\sum_{n=0}^\infty\, b_n\,q^{-nz}$
%%
%\todo{\tiny incoherence d'indice in \cite[page 52]{Rosen02} between
%  line -6 and line 11}
%%
with $b_n=h_K \,\frac{q^{n-\ggg+1}-1}{q-1}$ if $n>2\,\ggg-2$. Hence by
Equation \eqref{eq:relatzetaKzetav}, we have
\begin{align*}
  \zeta_v(z)&=(1-q_v^{-z})\;\zeta_K(z)=(1-q^{-z\deg v})
  \sum_{n=0}^\infty\; b_n\,q^{-nz}\\&=\sum_{n=0}^{\deg v\,-1} b_n\,q^{-nz}
  +\sum_{n=\deg v}^\infty (b_n-b_{n-\deg v})\,q^{-nz}\,.
\end{align*}
Hence by identification, if $n\geq 2\,\ggg+\deg v$, we have
\[
c_n=b_n-b_{n-\deg v}=h_K \;\frac{q^{n-\ggg+1}-q^{n-\ggg+1-\deg v}}{q-1}
=\frac{h_K \;q^{1-\ggg}\,(1-q_v^{-1})}{q-1}\;q^n\,.
\]
Therefore, by a geometric series argument, for every $n\in\NN$, we have
\begin{align*}
&\card\,\{I\in \I_v^+:\redn(I)\leq q^n\}=\sum_{i=0}^n\;c_i=
\sum_{i=2\ggg+\deg v}^nc_i+\bigO(1)\\=\;&
\frac{h_K \;q^{1-\ggg}\,(1-q_v^{-1})}{(q-1)}\;\frac{q^{n+1}}{q-1}+\bigO(1)
=\frac{h_K \;q^{2-\ggg}\,(1-q_v^{-1})}{(q-1)^2}\;q^n+\bigO(1)\,.
\end{align*}
Since $\redn(I)\in q^\NN$ for every $I\in \I_v^+$, this proves the result.
\cqfd
%the first claim of Lemma \ref{lem:gauscounting} follows.
%
%For every $I\in \I_v^+$, let $[I]$ be the class of $I$ modulo the
%principal ideals in the class group of $R_v$. Let $c'_n=\card\{I\in
%\I_v^+ : [I]=0,\;\redn(I) = q^n\}$. For every $n\in\NN$, there are
%exactly $h_K$ divisor classes of degree $n$ by
%\cite[Lem.~5.8]{Rosen02}.
%%%
%%\todo{mais combien de classes d'idéaux de $R_v$ par degré
%%  et en quels degrés ?}
%%%
%Hence by \cite[Lem.~4.1.2]{Goss98}, we have $c'_n=0$ if $n\notin (\deg
%v)\NN$ and $c'_{n\deg v}=\frac{1}{h_K}\, c_{n\deg v}$.
%%%
%\todo{à vérifier}
%%%
%Therefore
%\begin{align*}
%&\card\{I\in\I_v^+ : [I]=0,\;\redn(I) \leq q_v^n\}
%=\frac{1}{h_K}\sum_{i=0}^n\;c_{i\deg v}
%\\=\;&\frac{q^{1-\ggg}\,(1-q_v^{-1})}{(q-1)}\;\frac{q_v^{n+1}}{q_v-1}+\bigO(1)
%=\frac{q^{1-\ggg}}{q-1}\;q_v^n+\bigO(1)\,.
%\end{align*}
%The map $s\mapsto sR_v$ from $R_v\ssm\{0\}$ to the subset of principal
%ideals in $\I_v^+$ is a $(q-1)$-to-$1$ map since $\card\; R_v^\times
%=q-1$ by Equation \eqref{eq:inversiRv}. The second claim of Lemma
%\ref{lem:gauscounting} hence follows.
%\cqfd

\medskip
The following lemma is an effective version of
\cite[Lem.~3]{Poonen97}. Its proof follows the one of
\cite[Th.~328]{HarWri08} given in Chap.~XXII, \S 22.9 where $\ZZ$ is
replaced by $R_v$. Again, we add a proof for completeness. We denote
by $\ga$ the Euler constant.

\blemm\label{lem:PoonenRosen}
If $c'=\frac{q^{\ggg-1}\,(q-1)\,\ln q}{(1-q_v^{-1})\,e^\ga\;h_K}$, then
${\displaystyle \liminf_{\redn (I)\ra+\infty}}
\frac{\varphi_v(I)\;\ln\ln(\redn (I))}{\redn (I)}= c'$.
\elemm

In particular, since $\redn(s)=|s|=q_v^{-v(s)}$ for every $s\in R_v
\ssm \{0\}$, we have
\begin{equation}\label{eq:minovarphiv}
  \exists\; c_{\varphi_v}\in\;]0,1],\quad\forall s\in R_v\ssm\{0\},\qquad
  \;\varphi_v(s)\geq c_{\varphi_v}\frac{|s|}{\max\{1,\ln(-v(s))\}}\;.
\end{equation}

\dem Since we will only use the minoration \eqref{eq:minovarphiv}, we
only prove that the lower limit in the statement of Lemma
\ref{lem:PoonenRosen} is at least $c'$. Let $F:\;]0,+\infty[\;\ra \RR$
    be the map
\[
t\mapsto F(t)=(\ln t)\Big(1-\frac{1}{t}\Big)^{\frac{t}{\ln t}}\prod_{\redn (\ppp)\leq t}
\Big(1-\frac{1}{\redn (\ppp)}\Big)\,.
\]
By \cite[Th.~3]{Rosen99} (which gives an asymptotic expansion of the
partial product over all closed points of ${\mathfrak C}$, including
the one corresponding to $v$, which explains the factor $1-q_v^{-1}$
in the constant $c'$), as $t\ra+\infty$, we have $\prod_{\redn
  (\ppp)\leq t} \big(1-\frac{1}{\redn (\ppp)}\big)\sim c'\frac{1}{\ln
  t}$. Hence as $t\ra+\infty$ we have $F(t)\sim
c'(1-\frac{1}{t})^{\frac{t}{\ln t}} \sim c'(1-\frac{1}{\ln t})$, which
tends to $c'$.  For every $I\in \I_v^+$, let $A_I=\{\ppp: \ppp\mid I,
\;\redn (\ppp)\leq \ln\redn(I)\}$ and $B_I=\{\ppp: \ppp\mid I, \;\redn
(\ppp)> \ln\redn(I)\}$. Since $\redn$ is completely multiplicative, we
have $(\ln\redn (I))^{\card\; B_I}\leq \prod_{\ppp\in B_I} \redn
(\ppp)\leq \redn(I)$, hence $\card\; B_I\leq
\frac{\ln\redn(I)}{\ln\ln\redn(I)}$.  Thus as $\ln\redn(I)\ra+\infty$,
we have
\begin{align*}
&\frac{\varphi_v(I)\,\ln\ln\redn(I)}{\redn (I)} \\=\;&
\ln\ln\redn(I)\prod_{\ppp\,\mid\, I}\Big(1-\frac{1}{\redn (\ppp)}\Big)
\geq \ln\ln\redn(I)\Big(1-\frac{1}{\ln\redn(I)}\Big)^{\card\; B_I}
\prod_{\ppp\in A_I}\Big(1-\frac{1}{\redn (\ppp)}\Big)
\\\geq \;& \ln\ln\redn(I)\Big(1-\frac{1}{\ln\redn(I)}
\Big)^{\frac{\ln\redn(I)}{\ln\ln\redn(I)}}
\prod_{\redn (\ppp)\leq \ln\redn(I)}\Big(1-\frac{1}{\redn (\ppp)}\Big)
=F(\ln\redn(I))\sim c'\,. \;\Box
\end{align*}

\medskip
Denote by $\varpi_v: \I_v^+ \ra \NN$ the omega function counting the
prime factors of ideals:
\[
\varpi_v:I\mapsto \card\{\ppp:\ppp\mid I\}\,.
\]
We define $\varpi_v(s)=\varpi_v(sR_v)$ for every $s\in R_v \ssm
\{0\}$. For every $I\in\I_v^+$, with the notation of the above proof,
we have $\card\; A_I=\bigO\big(\frac{\ln\redn(I)}{\ln\ln\redn(I)}
\big)$ by the prime number theorem in $K$ (see for instance
\cite[Theo.~5.12]{Rosen02}) and $\card\; B_I\leq
\frac{\ln\redn(I)}{\ln\ln\redn(I)}$. Hence as $N(I)\ra+\infty$, we
have $\varpi_v(I)=\card\, (A_I\cup B_I)=\bigO\big(
\frac{\ln\redn(I)}{\ln\ln\redn(I)} \big)$.  In particular, since
$\redn(s)=|s|=q_v^{-v(s)}$ for every $s\in R_v \ssm \{0\}$, and since
$\varpi_v(s)=0$ when $s\in R_v^\times$, we have
\begin{equation}\label{eq:majoomegav}
  \exists\; c_{\varpi_v}>0,\quad\forall s\in R_v\ssm\{0\},\qquad
  \;\varpi_v(s)\leq c_{\varpi_v}\frac{-v(s)}{\max\{1,\ln(-v(s))\}}\;.
\end{equation}

\medskip
We fix $n\in\NN\ssm\{0\}$ throughout this paper. We denote by
$(e_1,\ldots, e_n)$ the canonical basis of the product $K_v$-vector
space $K_v^{\,n}$. Let $\|\;\|:K_v^{\,n}\ra [0,+\infty[$ be the
standard norm $(x_1,\ldots,x_n)\mapsto \max_{1\leq i\leq n}
|\,x_i\,|$.  We denote by $\vol_v$ the normalized Haar measure on the
locally compact additive group $K_v$ such that $\vol_v(\OOO_v)=1$. Let
$\vol_v^{\,n}$ be the normalized Haar measure on $K_v^{\,n}$ such that
$\vol_v^{\,n}(\OOO_v^{\,n}) =1$. Note that for every $g\in
\GL_n(K_v)$, we have
\begin{equation}\label{eq:changevarHaar}
d\vol_v^{\,n} (gx)=|\,\det(g)\,| \; d\vol_v^{\,n}(x)\,.
\end{equation}
In particular, we have $\vol_v (\pi_v\OOO_v) =q_v^{\;-1}$ and
\begin{equation}\label{eq:volvOvtimes}
  \vol_v (\OOO_v^\times)=\vol_v (\OOO_v\ssm\pi_v\OOO_v)
  =1-q_v^{\;-1}\,.
\end{equation}
If $G$ is a discrete subgroup of the additive group $K_v^{\,n}$ (for
instance any nonzero, not necessarily principal, ideal of $R_v$ when
$n=1$), we also denote by $\vol_v^{\,n}$ the unique Haar measure on
the quotient abelian topological group $K_v^{\,n}/G$ such that the
covering map $K_v^{\,n}\ra K_v^{\,n}/G$ locally preserves the measure.

\subsection{Lattices}
\label{subsec:lattices}

An {\it $R_v$-lattice} $L$ in $K_v^{\,n}$ is a free rank-$n$
$R_v$-submodule in $K_v^{\,n}$ that generates $K_v^{\,n}$ as a
$K_v$-vector space. It is a discrete cocompact additive subgroup of
$K_v^{\,n}$. For instance, a nonzero ideal $I$ of $R_v$ is an
$R_v$-lattice in $K_v$ if and only if it is principal.

The {\it covolume} of $L$, denoted by $\covol(L)$, is defined as the
measure of the (compact) quotient space $K_v^{\,n}/L$~:
\[
\covol(L)=\vol_v^{\,n}(K_v^{\,n}/L)\,.
\]
For every $g\in \GL_n(K_v)$, by Equation \eqref{eq:changevarHaar}, we
have
\begin{equation}\label{eq:changevarcovol}
\covol(gL)=|\,\det(g)\,| \; \covol(L)\,.
\end{equation}
In particular, if $\lambda\in K_v^\times$, then
$\covol(\lambda\,L) = |\lambda|^n \covol(L)$.  Since the
set of values of $|\cdot|$ is $\{0\}\cup q_v^{\,\ZZ}$, every
$R_v$-lattice is hence homothetic under $K_v^\times$ to an
$R_v$-lattice with covolume in $[1,q_v^{\,n}]$. For example,
$R_v^{\,n}$ is an $R_v$-lattice in $K_v^{\,n}$, and by for instance
\cite[Lem.~14.4]{BroParPau19}, we have
\begin{equation}\label{eq:covolRv}
\covol(R_v^{\,n})=q^{(\ggg-1)n}\,.
\end{equation}
The {\it normalized covolume} of an $R_v$-lattice $L$ is
$\frac{\covol(L)}{\covol(R_v^{\,n})}$, which belongs to $q_v^{\,\ZZ}$
since $\GL_n(K_v)$ acts transitively on the set of $K_v$-basis of
$K_v^{\;n}$, hence on the set of $R_v$-lattices of $K_v^{\;n}$, and by
using Equation \eqref{eq:changevarcovol}.

An $R_v$-lattice $L$ in $K_v^{\,n}$ is said to be

%Délire !!!
%$\bullet$~ {\it decorated} if it is endowed with an equivalence class
%$[(b_1,\dots,b_n)]$, called its {\it decoration}, of an $R_v$-basis
%$(b_1,\dots,b_n)$ of $L$ such that the determinant of $(b_1,\dots,b_n)$ in
%the canonical basis $(e_1,\dots, e_n)$ of $K_v^{\,n}$ is $1$, where
%two $K_v$-basis of $K_v^{\,n}$ are equivalent if the determinant of
%their transition matrix is $1$ (or equivalently if their wedge products
%are equal).

$\bullet$~ {\it unimodular} if $\covol(L)=\covol(R_v^{\,n})$ (by
Equation \eqref{eq:covolRv} for instance, as well as for other
integrality purposes, it is not appropriate to define them by
requiring $\covol(L)=1$),

$\bullet$~ {\it special unimodular} if $L$ admits an $R_v$-basis
$(b_1,\ldots, b_n)$ such that $b_1\wedge\ldots\wedge b_n$ is equal to
the canonical generator $e_1\wedge\ldots\wedge e_n$ of the $n$-th
exterior power $\text{\Large $\wedge$}\!^n(K_v^{\,n})$ (where, as
already said, $(e_1,\ldots, e_n)$ is the canonical $K_v$-basis of
$K_v^{\,n}$).

$\bullet$~ {\it integral} if $L$ is contained in $R_v^{\,n}$,

$\bullet$~ {\it rational} if $L$ is contained in $K^n$,

$\bullet$~ {\it axial} if for every $i\in\intbra{1,n}$, we have
$(K_v e_i)\cap L\neq \{0\}$.

\noindent
Any element of $\GL_n(K_v)$ mapping the canonical $K_v$-basis
$(e_1,\ldots, e_n)$ of $K_v^{\,n}$ to a $K_v$-basis $(b_1,\ldots,
b_n)$ such that $b_1\wedge\ldots\wedge b_n=e_1\wedge\ldots\wedge e_n$
has determinant $1$. Hence by Equation \eqref{eq:changevarcovol},
special unimodular $R_v$-lattices are unimodular.

For instance, if $I_1,\ldots,I_n$ are nonzero principal ideals of
$R_v$, then $\prod_{i=1}^nI_i$ is an integral $R_v$-lattice in
$K_v^{\,n}$. Note that an integral $R_v$-lattice, being a finite index
subgroup of $R_v^{\,n}$, is axial.  If $x$ is an axial $R_v$-lattice,
since $x$ has an $R_v$-basis $(b_1,\ldots,b_n)$ which is a $K_v$-basis
of $K_v^{\;n}$, by Kramer's formula to solve a system of $n-1$
linearly independent linear equations in $n$ variables in terms of one
of these variables, for every $i\in\llbracket 1,n\rrbracket$, the
intersection $R_ve_i\cap x$ is a rank-$1$ $R_v$-submodule of $K_ve_i$.
Hence there exists $\lambda_i\in K_v^\times$ such that $R_ve_i\cap x=
R_v\lambda_i e_i$.

For every integral $R_v$-lattice $L$, by the structure theorem of
finitely generated torsion modules over a Dedekind ring (see for
instance \cite[Theo.~1.41]{Narkiewicz04} without the uniqueness
statement), there exist unique nonzero ideals $I_1,\ldots,
I_n\in\I_v^+$ such that $I_1\,|\,I_2 \,|\ldots|\,I_n$ and
$R_v^{\,n}/L$ is isomorphic to $\prod_{i=1}^n R_v/I_i$ as an
$R_v$-module. The $n$-tuple $(I_1,\ldots, I_n)$, or the isomorphism
class of the $R_v$-module $R_v^{\;n}/L$, is called the {\it type} of
the integral lattice $L$. If $I_1=s_1R_v,\ldots, I_n=s_nR_v$ are
principal ideals, we will also say that the type of $L$ is
$(s_1,\ldots,s_n)$ (which is well defined modulo
$(R_v^\times)^n$). For instance, the type of $R_v^{\;n}$ is
$(1,\ldots, 1)$. The group $\GL_n(R_v)$ acts on the set of integral
$R_v$-lattices of $K_v^{\;n}$ and two integral $R_v$-lattices are in
the same $\GL_n(R_v)$-orbit if and only if they have the same type.
%%
%\todo{\tiny Référence? Si correct...}
%%

\subsection{Homogeneous spaces of lattices}
\label{subsec:homogspalat}

We denote by $I_n$ the $n\times n$ identity matrix. Let
$PG=\PGL_n(K_v)=\GL_n(K_v)/(K_v^\times I_n)$, which is a totally
disconnected metrisable locally compact topological group, and
$P\Ga=\PGL_n(R_v)=\GL_n(R_v)/(R_v^\times I_n)$, which is a nonuniform
lattice in $PG$. Throughout the paper, for every element
$g=(g_{ij})_{1\leq i,j\leq n} \in \GL_n(K_v)$, we denote by
$[g]=[g_{ij}]_{1\leq i,j\leq n}$ its image in $PG$, when
necessary. Otherwise, abusing notation, we omit the brackets.

Let $P\!\X$ be the homogeneous space $PG/P\Ga$, that identifies
$PG$-equivariantly with the space of the homothety classes $[L]
=K_v^\times L$ under $K_v^\times$ of the $R_v$-lattices $L$ in
$K_v^{\,n}$ by the orbital map $[g]P\Ga\mapsto [g R_v^{\,n}]$.
Contrarily to the case of the real field $\RR$ and the ring of
integers $\ZZ$ of the number field $\QQ$, in positive characteristic,
there is a difference between $\SL_n(K_v)/\SL_n(R_v)$ and $\PGL_n(K_v)
/\PGL_n(R_v)$ and it is sometines preferable to work with the latter
one, or with the following avatar.

\medskip
Let
\[
G_1=\GL_n^1(K_v)=\{g\in\GL_n(K_v):|\det g|=1\}\,,
\]
which is a unimodular totally disconnected metrisable locally compact
topological group with center $ZG_1=\OOO_v^{\,\times} I_n$.  We identify
the image of $G_1$ in $PG$ with $G_1/ZG_1$.  Let us denote $\Ga_1=
\GL_n(R_v)$. By Equation \eqref{eq:inversiRv}, we have $\Ga_1\subset
G_1$ and $P\Ga\subset G_1/ZG_1$. Besides, $\Ga_1$ is a nonuniform
lattice in $G_1$.

Finally, let $G=\SL_n(K_v)$, which is a unimodular closed normal
subgroup of $G_1$ with a split exact sequence of topological groups
\begin{equation}\label{eq:shortexcatseqG}
  1\longrightarrow G\longrightarrow G_1 \longrightarrow \OOO_v^{\,\times}
  \longrightarrow 1
\end{equation}
with section $\xi:\OOO_v^\times\ra G_1$ defined by
\[
\xi:\lambda\mapsto
\big(\begin{smallmatrix}\lambda & 0\\0&I_{n-1}\end{smallmatrix}\big)\,.
\]
Let $\Ga=\SL_n(R_v)$, which is a nonuniform lattice in $G$, with an
induced split exact sequence
\begin{equation}\label{eq:shortexcatseqGam}
  1\longrightarrow \Ga\longrightarrow \Ga_1 \longrightarrow
  R_v^{\,\times} \longrightarrow 1
\end{equation}
with section $\xi_{\mid R_v^\times}$.

We endow $PG$ (respectively $G_1$ and $G$) with its right-invariant
Haar measure ${\tt m}_{PG}$ (respectively ${\tt m}_{G_1}$ and ${\tt
  m}_{G}$) such that its maximal compact-open subgroup $PG(\OOO_v)=
\PGL_n(\OOO_v)$ (respectively $G_1(\OOO_v)=\GL_n(\OOO_v)$ and
$G(\OOO_v) =\SL_n(\OOO_v)$ ) has Haar measure $1$. Equation
\eqref{eq:shortexcatseqG} induces a split exact sequence
$1\longrightarrow G(\OOO_v) \longrightarrow G_1(\OOO_v)
\longrightarrow \OOO_v^{\,\times} \longrightarrow 1$ of compact
groups. Since $\vol_v(\OOO_v^{\,\times})=1-q_v^{\,-1}$ by Equation
\eqref{eq:volvOvtimes}, for all $\lambda\in \OOO_v^{\,\times}$ and
$g\in G$, we hence have
\begin{equation}\label{eq:desintegG1surG}
  d\,{\tt m}_{G_1}(\xi(\lambda)g)=
  \frac{q_v}{q_v-1}\;\dvol_v(\lambda)\;d\,{\tt m}_{G}(g)\,.
\end{equation}

\medskip
Let $\X_1$ be the space of unimodular $R_v$-lattices in $K_v^{\;n}$,
endowed with the Chabauty topology. As justified by Equation
\eqref{eq:changevarcovol}, we identify
homeomorphically and $G_1$-equi\-var\-iant\-ly the homogeneous space
$G_1/\Ga_1$ with $\X_1$ by the orbital map $g\Ga_1\mapsto
g\,R_v^{\,n}$.
%Indeed, for every $g\in G_1$, by Equation \eqref{eq:changevarcovol},
%if $L=g R_v^{\,n}$, we have
%\[
%\covol(L)=|\det g|\covol(R_v^{\,n})= \covol(R_v^{\,n})\;;
%\]
%conversely, if $L$ is a unimodular $R_v$-lattice, if $(e'_1,
%\ldots, e'_n)$ is a free $R_v$-basis of $L$ which is a $K_v$-basis of
%$K_v^{\,n}$, then the unique element $g\in \GL_n(K_v)$ sending
%$(e_1,\ldots, e_n)$ to $(e'_1,\ldots, e'_n)$ maps $R_v^{\,n}$ to $L$,
%and $|\det g|=\frac{\covol(L)}{\covol(R_v^{\,n})}=1$, hence $g\in G_1$.

Since $\Ga_1$ is a discrete subgroup of the unimodular group $G_1$, we
endow the homogeneous space $G_1/\Ga_1$ with the unique
$G_1$-invariant measure such that the orbital map $G_1\ra G_1/\Ga_1$
defined by $g\mapsto g\,\Ga_1$ locally preserves the measure, and we
endow $\X_1$ with the corresponding measure ${\tt m}_{\X_1}$.  By for
instance \cite[Eq.~(41)]{HorPau24} (building on \cite[\S 3]{Serre71}),
we have
\begin{equation}\label{eq:totmassmX_1}
\|{\tt m}_{\X_1}\|=\frac{q_v-1}{q_v(q-1)}
\prod_{i=1}^{n-1}\frac{\zeta_v(-i)}{q_v^{\;i}-1}\,.
\end{equation}

Let $\X$ be the closed subspace of $\X_1$ consisting in the special
unimodular $R_v$-lattices in $K_v^{\;n}$, which is equal to the orbit
in $\X_1$ of the standard $R_v$-lattice $R_v^{\;n}$ under the action
of the subgroup $G$ of $G_1$. The stabiliser of $R_v^{\;n}$ in $G$ is
exactly $\Ga$. The homogeneous space $G/\Ga$ identifies
homeomorphically and $G$-equivariantly with the space $\X$ by the map
$g\Ga\mapsto g\,R_v^{\;n}$.  Since $\Ga=\Ga_1\cap G$, the inclusion
map $G\ra G_1$ induces an injection $G/\Ga\ra G_1/\Ga_1$ which is a
homeomorphism onto its image, and the following diagram is commutative:
\[
\begin{array}{ccc}
G/\Ga& \ra & G_1/\Ga_1\smallskip\\
^\simeq \!\downarrow\;\; & & \;\;\downarrow^\simeq\smallskip\\
\X & \hookrightarrow & \,\X_1\,.
\end{array}
\]
The compact group $\OOO_v^{\,\times}/R_v^{\,\times}$ acts
continuously and freely on the topological space $\X_1$ by $(\lambda
R_v^{\,\times},\Lambda)\mapsto \xi(\lambda) \Lambda$. By Equations
\eqref{eq:shortexcatseqG} and \eqref{eq:shortexcatseqGam}, the
inclusion map $\X\ra\X_1$ induces a homeomorphism
$\X\ra(\OOO_v^{\,\times}/R_v^{\,\times}) \backslash \X_1$.

Since $\Ga$ is a discrete subgroup of the unimodular group $G$, we
endow the homogeneous space $G/\Ga$ with the unique $G$-invariant
measure such that the orbital map $G\ra G/\Ga$ defined by $g\mapsto
g\,\Ga$ locally preserves the measure, and we endow $\X$ with the
corresponding measure ${\tt m}_{\X}$. We denote by $\vol'_v$ the
measure on $\OOO_v^{\,\times}/R_v^{\,\times}$ such that the map
$\OOO_v^{\,\times}\ra \OOO_v^{\,\times}/R_v^{\,\times}$ locally
preserves the measure. Its total mass is $\vol'_v (\OOO_v^{\,\times}/
R_v^{\,\times})=\frac{\vol_v(\OOO_v^{\,\times})}{\card\;R_v^{\,\times}}
=\frac{1-q_v^{-1}}{q-1}$. By Equation \eqref{eq:desintegG1surG}, for
all $\lambda \,R_v^{\,\times}\in \OOO_v^{\,\times}/R_v^{\,\times}$ and
$x\in\X$, we have
\begin{equation}\label{eq:mXmX1}
  d\,{\tt m}_{\X_1}(\xi(\lambda)x)=
  \frac{q_v}{q_v-1}\;\dvol'_v(\lambda \,R_v^{\,\times})\;d\,{\tt m}_\X(x)
\quad\text{and}\quad \|\,{\tt m}_{\X}\|=(q-1)\,\|\,{\tt m}_{\X_1}\|\,.
\end{equation}

\subsection{Systoles of lattices}
\label{subsec:systole}

The {\it (normalized) systole} of an $R_v$-lattice $L$ in $K_v^{\;n}$,
which depends only on the homothety class $[L]$ of $L$ modulo
$K_v^{\;\times}$, is defined by
\begin{equation}\label{eq:defisys}
\sys([L])=\sys(L)=
\Big(\frac{\covol (R_v^{\;n})}{\covol (L)}\Big)^{\frac{1}{n}}
\min_{w\in L\ssm\{0\}}\|w\|\,.
\end{equation}
If $L$ is unimodular, we simply have $\sys(L)=\min_{w\in L\ssm\{0\}}
\|w\|$.  Mahler's compactness criterion (see for instance
\cite[Theo.~1.1]{KleShiTom17}) says that for every $\epsilon>0$, the
{\it $\epsilon$-thick part} of $P\!\X$, defined by
\[
P\!\X^{\geq \epsilon}=
\big\{\;x\in P\!\X:\sys(x)\geq \epsilon\;\big\}\;,
\]
is compact in $P\!\X$. For every compact subset $K$ of $P\!
\X$, there exists $\epsilon>0$ such that $K\subset P\!\X^{\geq\epsilon}$.
Similarly, the $\epsilon$-thick parts 
\[
\X_1^{\geq \epsilon}=\big\{\;L\in\X_1:\sys(L)\geq \epsilon\;\big\}
\quad\text{and}\quad
\X^{\geq \epsilon}=\big\{\;L\in\X:\sys(L)\geq \epsilon\;\big\}
\]
of $\X_1$ and $\X$ respectively are compact. For every compact subset
$K$ of $\X_1$ or $\X$, there exists $\epsilon>0$ such that $K$ is
contained in $\X_1^{\geq\epsilon}$ or $\X^{\geq \epsilon}$. We denote
by $P\!\X^{<\epsilon}=P\!\X\ssm\, P\!\X^{\geq\epsilon}$,
$\X_1^{<\epsilon} =\X_1\ssm\,\X_1^{\geq \epsilon}$ and $\X^{<
  \epsilon} = \X\ssm\,\X^{\geq\epsilon}$ the {\it $\epsilon$-thin
  parts} of $P\!\X$, $\X_1$ and $\X$ respectively.

Since $\GL_n(K_v)$ acts transitively on the set of $R_v$-lattices in
$K_v^{\;n}$ and by Equation \eqref{eq:changevarcovol}, the set of
values of the (continuous) {\it systole function} $\sys:P\!\X\ra\RR$
is contained in $q_v^{\frac{1}{n}\ZZ}$. More precisely, let us prove
that we have
%%
%\todo{\tiny what is optimal upperbound in general?}
%%
\[
q_v^{-\frac{1}{n}\NN}\subset \sys(P\!\X)\subset
q_v^{\frac{1}{n}\ZZ}\cap[0,q_v q^{\ggg-1}]\,.
\]
The left inclusion follows by considering, for every $k\in \NN$, the
lattice $L=g_kR_v^{\;n}$ with $g_k=\begin{psmallmatrix} \pi_v^{-k} & 0
\\ 0 & I_{n-1}\end{psmallmatrix}$ (giving $\sys(L)=q_v^{-\frac{k}{n}}$
since $|\pi_v^{-k}|=q_v^k>1$). In order to prove the right inclusion,
let $[L] \in P\X$. Up to rescaling, we may assume that $\min_{w\in
  L\,\ssm\,\{0\}} \|w\|=1$. Then the closed ball $B(0,q_v^{-1})$ in
$K_v^{\;n}$ injects in $K_v^{\;n}/L$ by the ultrametric triangle
inequality. Thus by Equation \eqref{eq:covolRv} and the line before
Equation \eqref{eq:volvOvtimes}, we have as wanted
\[
\sys(L)= \Big(\frac{\covol (R_v^{\;n})}{\covol
  (L)}\Big)^{\frac{1}{n}} \leq\Big(\frac{q^{(\ggg-1)n}}{\vol_v^n
  (B(0,q_v^{-1}))}\Big)^{\frac{1}{n}}=\frac{q^{\ggg-1}}{\vol_v
  (\pi_v\OOO_v)}=q_v \,q^{\ggg-1}\,.
\]
When $\ggg=0$ and $\deg v=1$, we have $\sys(L)\leq 1$,
which is optimal since $\sys(R_v^{\;n})= 1$.

Since the image $\sys(P\!\X)$ is contained in $q_v^{\frac{1}{n}\ZZ}$,
there is a partition $P\!\X=\bigcup_{k=1}^{n} P\!\X_{k}$ into nonempty
closed and open subsets of $P\!\X$, defined by
\[
\forall\;k\in\intbra{1,n},\quad P\!\X_{k}=
\{x\in P\!\X:n\log_{q_v}(\sys(x))\equiv k-1\!\!\mod n\}\,.
\]

The norm $\|\cdot\|$ having values on $K_v^{\;n}\ssm\{0\}$ in
$q_v^\ZZ$, for every $R_v$-lattice $L$ and every $g\in \GL_n(K_v)$, by
Equation \eqref{eq:changevarcovol}, we have
\begin{align*}
n\log_{q_v}(\sys(gL))&\equiv
\log_{q_v}\Big(\frac{\covol (R_v^{\;n})}{\covol (gL)}\Big)
\equiv \log_{q_v}\Big(\frac{\covol (R_v^{\;n})}{\covol (L)}\Big)
-\log_{q_v}|\det g|
\\&\equiv n\log_{q_v}(\sys L)-\log_{q_v}|\det g|\mod n\,.
\end{align*}
Hence the image $G_1/ZG_1$ of $G_1$ in $PG$ acts transitively on each
one of the strata $P\!\X_{k}$ for $k\in\intbra{1,n}$.  Since
$P\Ga\subset G_1/ZG_1$, the stratum $P\!\X_{1}$ thus identifies
$(G_1/ZG_1)$-equi\-variantly with the homogeneous space
$(G_1/ZG_1)/P\Ga$.  Furthermore, for every $k'\in\NN$, the element
$g_{k'}=\begin{psmallmatrix}\pi_v^{k'}&0\\0&I_{n-1}\end{psmallmatrix}$
maps $P\!\X_{k}$ to $P\!\X_{k''}$ where
$k''\in\llbracket1,n\rrbracket$ satisfies $k''=k+k'\!\!\mod n$.

\subsection{Diagonal subgroups}
\label{subsec:diagosub}

We denote by $\wt A$ the diagonal subgroup of $\GL_n(K_v)$, and by $P\!A$
its image in $PG$, that we also call the {\it diagonal subgroup} of
$PG$. Let
%%
%\todo{appropriate signs of the exponents to be checked}
%%
\[
\wt D=\left\{\begin{pmatrix}\pi_v^{-k_1}
&  & 0 \\  & \ddots & \\ 0 &  & \pi_v^{-k_n}\end{pmatrix}:
k_1,\ldots,k_n\in\ZZ\right\}\subset \wt A\,.
\]
We denote by $P\!D$ the image of $\wt D$ in $PG$. Note that the
diagonal subgroup $P\!A$ is a closed noncompact subgroup of $PG$ which
is the direct product $P\!A=P\!A(\OOO_v)\,P\!D$ of its maximal compact
subgroup $P\!A(\OOO_v)=P\!A\cap \PGL_n(\OOO_v)$ and its discrete
subgroup $P\!D$. As seen at the end of the previous subsection
\ref{subsec:systole}, the group $P\!D$ permutes transitively the
strata $P\!\X_{k}$ for $k\in\intbra{1,n}$.

Also note that if $L$ is an axial $R_v$-lattice in $K_v^{\;n}$, then
$a\,L$ is an axial $R_v$-lattice for every $a\in \wt A$, and in
particular, every $R_v$-lattice homothetic to $L$ is axial. Hence we
may define an {\it axial $P\!A$-orbit} in $P\!\X$ to be a $P\!A$-orbit
which contains the homothety class of an axial $R_v$-lattice, or
equivalently a $P\!A$-orbit all of whose elements are homothety
classes of axial $R_v$-lattices.

\medskip
Let $\exp:\ZZ^n\ra \GL_n(K_v)$ be the map ${\bf k}= (k_1,\ldots,k_n)
\mapsto \begin{pmatrix}\pi_v^{-k_1} & & 0 \\ & \ddots & \\ 0 & &
  \pi_v^{-k_n}\end{pmatrix}$, which is an injective group morphism
with image $\wt D$. We have $\aaa=\exp(1-n,1,\dots, 1)$ by Equation
\eqref{eq:defaaa}.  We will also denote by $\exp$ its restriction to
\[
\ZZ^n_0=\{\mb{k}=(k_1,\ldots,k_n)\in\ZZ^n:k_1+\ldots+k_n=0\}\,.
\]
Note that there exists no global exponential map of matrices in
positive characteristic. The above map $\exp$ is a (very weak) ersatz
for it. 

We define $A_1=\wt A\cap G_1$, $A=\wt A\cap G$ and $D=\wt D\cap G_1=
\wt D\cap G= \exp(\ZZ^n_0)$. Thus $A_1$ and $A$ are the direct
products $A_1= A_1(\OOO_v) D$ and $A= A(\OOO_v) D$ of their maximal
compact subgroups $A_1(\OOO_v)=A_1\cap \GL_n(\OOO_v) =\wt A\cap
\GL_n(\OOO_v)$ and $A(\OOO_v)=A\cap \SL_n(\OOO_v)$ respectively with
their discrete subgroup $D$. The split exact sequence in Equation
\eqref {eq:shortexcatseqG} gives a split exact sequence
\begin{equation}\label{eq:shortexcatseqA}
  1\longrightarrow A\longrightarrow A_1 \longrightarrow \OOO_v^{\,\times}
  \longrightarrow 1
\end{equation}
with section $\xi$ (which has values in $A_1$).

\subsection{Homogeneous measures on diagonal orbits}
\label{subsec:homogdiag}

Recalling that $A_1(\OOO_v)$ and $A(\OOO_v)$ are the maximal
compact-open subgroups of the diagonal groups $A_1$ and $A$, we endow the
abelian locally compact groups $A_1$  and $A$ with their unique Haar measure
${\tt m}_{A_1}$ and ${\tt m}_{A}$ normalized so that
\begin{equation}\label{eq:volA1OOOv}
{\tt m}_{A_1}(A_1(\OOO_v))={\tt m}_{A}(A(\OOO_v))=1\,.
\end{equation}
%We denote by $\Haar_{A_1(\OOO_v)}=\Haar_{A_1}\!\!\mid_{A_1(\OOO_v)}$
%the restriction of the measure $\Haar_{A_1}$ to $A_1(\OOO_v)$.
%
%For instance, if $n=1$, then we have $A_1=A_1(\OOO_v)=\OOO_v^{\times}
%=\FF_{q_v}^{\;\times} +\pi_v\OOO_v$ and therefore $\Haar_{A_1}
%=\Haar_{A_1(\OOO_v)} =\frac{q_v}{q_v-1} \vol_v\!\mid_{\OOO_v^\times}$.
By Equation \eqref{eq:shortexcatseqA}, for all $\lambda\in
\OOO_v^{\,\times}$ and $a\in A$, we have
\begin{equation}\label{eq:decompdmA1}
d\,{\tt m}_{A_1}(\xi(\lambda)a)=
\frac{q_v}{q_v-1} \;d\vol_v(\lambda)\;d\,{\tt m}_A(a)\,.
\end{equation}
We denote by ${\tt m}_{\ZZ^n_0}$ the counting measure on $\ZZ^n_0$, and,
in order to simplify notation,
\[
da=d{\,{\tt m}_{A_1}}_{\mid A_1(\OOO_v)}(a), \quad da=d{\,{\tt m}_{A}}_{\mid
  A(\OOO_v)} (a)\quad \text{and}\quad d\,\mb{k}=d\,{\tt m}_{\ZZ^n_0}
(\mb{k})\,,
\]
which are measures on $A_1(\OOO_v)$, $A(\OOO_v)$ and $\ZZ^n_0$
respectively.

The maps $(a,\mb{k})\mapsto a\exp(\mb{k})$ from $A_1(\OOO_v)\times
\ZZ^n_0$ to $A_1$ and from $A(\OOO_v)\times \ZZ^n_0$ to $A$
are isomorphisms of topological groups, and we have
\[
d\,{\tt m}_{A_1}(a\exp(\mb{k}))=da\;d\,\mb{k}\quad\text{and}\quad
d\,{\tt m}_{A}(a\exp(\mb{k}))=da\;d\,\mb{k}\,.
\]

For every $x\in\X_1$ (resp.~$x\in\X$), with $\theta_x:a\mapsto ax$ the
orbital map from $A_1$ to $A_1x$ (resp.~$A$ to $Ax$), we define the
{\it orbital measure} $\ov\mu_x=\ov\mu_{A_1x}$
(resp.~$\mu_x=\mu_{Ax}$)
%
%\todo{\tiny voir si on garde seulement $\ov\mu_x$ et $\mu_x$}
%
on the orbit $A_1x$ (resp.~$Ax$) by
\[
\ov\mu_x=\ov\mu_{A_1x}=(\theta_x)_*{\tt m}_{A_1}\quad\text{(resp.~}\;
\mu_x=\mu_{Ax}=(\theta_x)_*{\tt m}_{A}\,\text{)}\;,
\]
so that for $a\in A_1(\OOO_v)$ (resp.~$a\in A(\OOO_v)$) and
$\mb{k}\in\ZZ^n_0$, we have
\[
d\,\ov\mu_x(a\exp(\mb{k})\,x)=da\;d\,\mb{k}\quad\text{(resp.~}\;
d\mu_x(a\exp(\mb{k})\,x)=da\;d\,\mb{k}\,\text{)}\,.
\]
If $\theta_x$ is a proper map (see Corollary \ref{coro:classdivunimod}
for characterisations), then $\ov\mu_{x}$ is an $A_1$-invariant
infinite locally finite measure on $\X_1$ with support equal to the
orbit $A_1x$ (resp.~$\mu_{x}$ is an $A$-invariant locally finite
infinite measure on $\X$ with support equal to the orbit $Ax$).

By Equation \eqref{eq:decompdmA1} and by the definition of $\vol'_v$
at the end of Section \ref{subsec:homogspalat}, for every $x\in \X$,
the measure $\ov\mu_x$ on $\X_1$ is an average of orbital measures on
$\X$~:
\begin{equation}\label{eq:desintegovmumu}
  \ov\mu_x=\frac{q_v}{q_v-1}\;
  \int_{\lambda\in\OOO_v^{\,\times}}\xi(\lambda)_*\mu_x\;d\vol_v(\lambda)=
  \frac{q_v(q+1)}{q_v-1}\;\int_{\OOO_v^{\,\times}/R_v^{\,\times}}
  \xi(\lambda)_*\mu_x\;d\vol'_v(\lambda R_v^{\,\times})\,.
\end{equation}

\section{A classification of the divergent diagonal orbits}
\label{sec:classdivdiagorb}
  
The following characterisation of the divergent diagonal orbits is due
to Margulis in the case of the real field $\RR$ and the ring of
integers $\ZZ$, see \cite[Theo.~1.2]{TomWei03}. See for instance
\cite{Weiss04,Weiss06,SolTam23} for complementary studies in the real
case. Our proof follows the same scheme of proof as in the real case.

\btheo \label{theo:defidiv} For every $x=[L]\in P\!\X$, the following
assertions are equivalent.
\begin{enumerate}
\item\label{item1:defidiv}
  The map $a\mapsto a\,x$ from $P\!A$ to $P\!\X$ is a proper map.
\item\label{item2:defidiv}
  The map $d\mapsto d\,x$ from $P\!D$ to $P\!\X$ is a proper map.
\item\label{item3:defidiv}
  There exists $g\in (P\!A)(\PGL_n(K))$ such that $[L]=g[R_v^{\;n}]$.
\item\label{item4:defidiv}
  The orbit $P\!A\,x$ contains the homothety class of an integral $R_v$-lattice.
\item\label{item5:defidiv}
  The orbit $P\!A\,x$ contains the homothety class of an axial $R_v$-lattice.
\item\label{item6:defidiv}
  Every element of the orbit $P\!A\,x$ is the homothety class of an axial
  $R_v$-lattice.
\end{enumerate}
\etheo

If one of the above assertions is satisfied, we say that the orbit $P\!A
\,x$ of $x$ by the diagonal subgroup $P\!A$ is {\it divergent}. Hence the
divergent $P\!A$-orbits are the axial ones (as defined in Subsection
\ref{subsec:diagosub}).

\medskip
\dem Assertion \eqref{item3:defidiv} implies Assertion
\eqref{item4:defidiv}, since $K$ is the field of fractions of $R_v$,
hence for every $g'\in \GL_n(K)$, if $r\in R_v\ssm \{0\}$ is the
product of the denominators of the nonzero entries of $g'$ written as
fractions in $R_v$, then $rg'R_v^{\;n}$ is an integral $R_v$-lattice.
It is immediate that Assertion \eqref{item4:defidiv} implies Assertion
\eqref{item5:defidiv}, since an integral $R_v$-lattice is an axial
$R_v$-lattice. We have already seen in Subsection
\ref{subsec:diagosub} that the Assertions \eqref{item5:defidiv} and
\eqref{item6:defidiv} are equivalent.

Let us prove that Assertion \eqref{item6:defidiv} implies Assertion
\eqref{item1:defidiv}. Indeed, assume that $L$ is axial and normalized
in its homothety class to have covolume between $1$ and
$q_v^{\;n}$. Every element $a\in \wt A$ may be multiplied by a central
element of $\GL_n(K_v)$ in order to have absolute value of its
determinant between $1$ and $q_v^{\;n}$.  Then if $a$ goes to
infinity in $\wt A$, it has a diagonal entry that goes to $0$. Hence
the $R_v$-lattice $aL$ has its covolume remaining between $1$ and
$q_v^{\,2n}$, and has a nonzero vector on the coordinate axis
corresponding to that diagonal entry that goes to $0$.  Thus its image
in $P\!\X$ leaves every compact subset of $P\!\X$ by Mahler's
compactness criterion. Note that Assertion \eqref{item1:defidiv} and
Assertion \eqref{item2:defidiv} are equivalent, since $P\!A(\OOO_v)$
is compact and $P\!A=P\!A(\OOO_v)\,PD$.

It remains to prove that Assertion \eqref{item2:defidiv} implies
Assertion \eqref{item3:defidiv}. We first give two lemmas.

\blemm\label{lem:defidiv1} There exist $c>1$, a bounded open
neighborhood $W$ of $0$ in $K_v^{\;n}$ and a finite subset $F$ of
$D=\wt D\cap \SL_n(K_v)$, such that for every $g''\in \GL_n(K_v)$ with
$\det(g'')\in[1,q_v^{\;n}]\,$, there exists $f\in F$ such that for
every $w\in (g''R_v^{\;n})\cap W$, we have \[\|fw\|\geq c\|w\|\,.\]
\elemm

\dem Each element in $\GL_n(K_v)$ with absolute value of its
determinant in $[1,q_v^{\;n}]$ only multiplies the volume
$\vol_v^{\;n}$ by a constant in $[1,q_v^{\;n}]$ by Equation
\eqref{eq:changevarHaar}. Hence there exists an open ball $W$ centered
at $0$ in $K_v^{\;n}$ with small enough radius such that for every
element $g''\in \GL_n(K_v)$ with $\det(g'')\in[1,q_v^{\;n}]$, the
$K_v$-linear subspace generated by $(g''R_v^{\;n})\cap W$ is a proper
$K_v$-linear subspace of $K_v^{\;n}$.

For every $d\in\intbra{1, n-1}$, let $\operatorname{Gr}_d(K_v^{\;n})$
be the Grassmannian space of $d$-dimensional $K_v$-linear subspaces of
$K_v^{\;n}$, endowed with the Chabauty topology. Let us prove that
there exist $c_d>1$ and a finite subset $F_d$ of elements in $D$ such
that for every $V\in\operatorname{Gr}_d (K_v^{\;n})$, there exists
$f\in F_d$ such that for every $w\in V$, we have $\|fw\|\geq c_d
\|w\|$. This proves the result by taking $F=\bigcup_{1\leq d\leq n-1}
F_d$ and $c=\min_{1\leq d\leq n-1} c_d>1$.  By the compactness of
$\operatorname{Gr}_d(K_v^{\;n})$ and of the unit sphere of
$K_v^{\;n}$, and by the homogeneity of the norm $\|\cdot\|$ (so that
$\|\pi_v^{\,\log_{q_v}\|w\|} w\|=1$ for every $w\in
K_v^{\;n}\ssm\{0\}$), we only have to prove that for every $V\in
\operatorname{Gr}_d (K_v^{\;n})$, there exists $a\in D$ such that for
every $w\in V$ with norm $1$, we have $\|a\,w\|>1$.

Since $d<n$, there exists $i_0\in\intbra{1,n}$ such that $K_ve_{i_0}$
is not contained in $V$. We claim that there exists $\epsilon_0>0$
such that for every $w=(w_1,\ldots,w_n)\in V$ with norm $1$, there
exists $j_w\in\llbracket1, n\rrbracket$ different from $i_0$, such
that $|w_{j_w}|\geq \epsilon_0$. Otherwise, for every $k\in\NN$, there
exists $w^{(k)}=(w_{1,k},\ldots,w_{n,k})\in V$ with $\|w^{(k)}\|=1$
and $|w_{i,k}|\leq \frac{1}{k+1}$ for every $i\in\llbracket1,
n\rrbracket\ssm\{i_0\}$. Up to extracting a subsequence by the
compactness of the unit sphere of $K_v^{\;n}$ and since $V$ is closed
in $K_v^{\;n}$, the sequence $(w^{(k)})_{k\in\NN}$ converges to a unit
norm vector $w_\infty$ in $V\cap (K_ve_{i_0})$, contradicting the fact
that $K_ve_{i_0}$ is not contained in $V$.

Now, for every $k\in\NN$, let $a^k$ be the diagonal matrix with
diagonal coefficients $a^k_{ii}=\pi_v^{-k}$ if $i\neq i_0$ and
$a^k_{i_0i_0} =\pi_v^{(n-1)k}$, which belongs to $D=\exp(\ZZ_0^n)$.
Then, if $k$ is large enough (for instance $k=\lceil -\log_{q_v}
\epsilon_0\rceil+1$), for every element $w=(w_1,\ldots,w_n)\in V$ with
norm $1$, we have
\[
\|a^k\,w\|\geq |\pi_v^{-k}w_{j_w}|=q_v^{\,k}\,|w_{j_w}|
\geq q_v^{\,k}\,\epsilon_0>1\;,
\]
which proves the result.
\cqfd

\blemm\label{lem:defidiv2} For every element $g'\in\GL_n(K_v)\ssm\wt
A\,\GL_n(K)$, for every bounded open neighborhood $W$ of $0$ in
$K_v^{\;n}$, for every finite subset $J$ of $g'R_v^{\;n}\ssm \{0\}$,
for every finite subset $C$ of $D$, there exists $a\in D\ssm C$ such
that
\[
(aJ)\cap W=\emptyset\,.
\]
\elemm 

\dem Let $g',W,J,C$ be as in the statement. We first claim that
there exists $i_0\in\intbra{1,n}$ such that
\[
(K_v e_{i_0})\cap (g'R_v^{\;n})=\{0\}\,.
\]
Assume for a contradiction that for every $i\in\intbra{1,n}$ there
exist $\underline{w}_i\in R_v^{\;n}$ and $a_i\in K_v\ssm\{0\}$ such
that $g'\underline{w}_i=a_i e_i$. Then $(\underline{w}_1,\ldots,
\underline{w}_n)$ is a $K$-basis of the $K$-linear space
$\oplus_{1\leq i\leq n} Ke_i$, and the transition matrix $P$ from the
canonical basis $(e_1,\ldots, e_n)$ to this basis (so that
$Pe_i=\underline{w}_i$ for every $i\in\intbra{1,n}$) belongs to
$\GL_n(K)$. Let $a'\in\wt A$ be the diagonal matrix with diagonal
coefficients $a_1,\ldots, a_n$ in this order. Then the linear map
$(a')^{-1}g'P$ fixes the canonical basis, hence is the identity. Thus
$g'=a'P^{-1}\in \wt A\,\GL_n(K)$, a contradiction to the assumption
that $g'\notin \wt A\,\GL_n(K)$.

Now, for every $k\in\NN$, let $a^k$ be the diagonal matrix with
diagonal coefficients $a^k_{ii}=\pi_v^{-k}$ if $i\neq i_0$ and
$a^k_{i_0i_0} =\pi_v^{(n-1)k}$, which belongs to $D$, and does not
belong to $C$ if $k$ is large enough.  For every $w=(w_1,\ldots, w_n)
\in g'R_v^{\;n}\ssm\{0\}$, by the preliminary claim, there exists $i
\neq i_0$ such that $w_i\neq 0$. Hence $\|a^k\,w\|\geq |\pi_v^{-k}w_i|
= q_v^{\,k}\,|w_{i}|$ tends to $+\infty$ as $k\ra+\infty$.  Since $J$
is finite and $W$ is bounded, this implies that for $k$ large enough,
we have $(a^kJ)\cap W=\emptyset$ and $a^k \in D\ssm C$.
\cqfd

\medskip
\noindent{\bf Proof that Assertion \eqref{item2:defidiv} implies
  Assertion \eqref{item3:defidiv}. } Let us fix $g\in PG \ssm (P\!A)
(\PGL_n(K))$, and let us prove that the map $d\mapsto d g[R_v^{\;n}]$
from $P\!D$ to $P\!\X$ is not proper, which concludes the proof of
Theorem \ref{theo:defidiv}. We fix a representative $g'\in \GL_n(K_v)
\ssm\wt A\,\GL_n(K)$ of $g$ with $\det(g')\in[1,q_v^{\;n}]$. Let us
prove that the map $d\mapsto [dg'R_v^{\;n}]$ from the discrete space
$D$ to $P\!\X$ is not proper, which implies the result.

Let $c,W,F$ be as in Lemma \ref{lem:defidiv1}. Without loss of
generality, we may assume that $\id\in F$. Let
\[
W_0\subset\bigcap_{f\in F\cup F^{-1}} fW
\]
be an open ball centered at $0$ in $K_v^{\;n}$, contained in
$\bigcap_{f\in F\cup F^{-1}} fW$ hence in $W$. Let
\begin{equation}\label{eq:defidCprime}
C'=\Big\{[L]\in P\X: \;\frac{\covol (L)}{\covol (R_v^{\;n})}
\in[1,q_v^{\;n}],\;\; L\cap W_0=\{0\}\Big\}\;,
\end{equation}
which is a compact subset of $P\X$ by Mahler's compactness criterion
(and by the definition \eqref{eq:defisys} of the systole).  For every
finite subset $C$ of $D$, let us prove that there
exists an element $d_C\in D\ssm C$ such that
$[d_Cg'R_v^{\;n}]\in C'$, which concludes the proof.

Let $J=(g'R_v^{\;n}\cap C^{-1}W)\ssm\{0\}$, which is a finite subset
of $g'R_v^{\;n}\ssm\{0\}$ since $C^{-1}W$ is bounded. By Lemma
\ref{lem:defidiv2}, there exists $d_0\in D\ssm C$ such that
\begin{equation}\label{eq:defid0}
(d_0 J)\cap W=\emptyset\,.
\end{equation}
Let us define by induction on $k\in\NN$ an element $d_k\in D$ such
that with the notation $\wt d_k=d_0\ldots d_k$, if $k\geq 1$, then
$d_k\in F$ and
\begin{equation}\label{eq:dilatlarat}
\forall\;w\in (\wt d_{k-1}
\,g'R_v^{\;n})\cap W, \quad \|d_k w\| \geq c\,\|w\|\,.
\end{equation}

Let $k\in\NN$, assume that $d_0,\ldots, d_k$ have been constructed,
and note that $\wt d_k=d_0\ldots d_k$ belongs to $D$.  By Lemma
\ref{lem:defidiv1}, applied with $g''=\wt d_k \,g'$ which has absolute
value of its determinant in $[1,q_v^{\;n}]$, there exists $d_{k+1}\in
F$ such that for every $w\in (\wt d_k \,g'R_v^{\;n})\cap W$, we have
$\|d_{k+1}w\| \geq c\,\|w\|$. This concludes the induction.

For every $k\geq 1$, by the definition of $W_0$ which is contained in
$fW$ for every $f\in F$, we have
\begin{equation}\label{eq:inddecreasnonzerovect1}
(\wt d_k \,g'R_v^{\;n})\cap W_0 \subset (\wt d_k \,g'R_v^{\;n})\cap (d_k W)
=d_k\big((\wt d_{k-1} \,g'R_v^{\;n})\cap W\big)\,.
\end{equation}
Hence by Equation \eqref{eq:dilatlarat}, the minimal norm $\n_k$ of a
nonzero vector in $(\wt d_k \,g'R_v^{\;n})\cap W_0$ is at least $c$
times the minimal norm of a nonzero vector in $(\wt d_{k-1}
\,g'R_v^{\;n})\cap W$. Since $W_0$ is a ball centered at $0$ contained
in $W$, this implies that either $(\wt d_{k-1} \,g'R_v^{\;n})\cap
W_0=\{0\}$ or that $\n_k\geq c\,\n_{k-1}$.  Since $W_0$ is bounded and
$c>1$, this implies by a decreasing induction that there exists
$k\in\NN$ such that $(\wt d_k \,g'R_v^{\;n})\cap W_0=\{0\}$. Let
$k_*\in\NN$ be the smallest element $k\in\NN$ for which this equality
is true, so that
\begin{equation}\label{eq:inddecreasnonzerovect2}
  (\wt d_{k_*} \,g'R_v^{\;n})\cap W_0=\{0\}\quad{\rm and}\quad
  \forall\;k\in\llbracket0, k_*-1\rrbracket,\quad
  (\wt d_k\,g'R_v^{\;n})\cap W_0\neq\{0\}\,.
\end{equation}
By Equation \eqref{eq:changevarcovol}, we have $\frac{\covol(\wt d_k
  \,g'R_v^{\;n})} {\covol(R_v^{\;n})} =|\det(\wt d_k \,g')|\in
[1,q_v^{\;n}]$.  In particular, by the definition of $C'$ in Equation
\eqref{eq:defidCprime}, we have $[\wt d_{k_*} \,g'R_v^{\;n}]\in
C'$. Let us prove that $\wt d_{k_*}\notin C$, which gives the wanted
result using $d_C=\wt d_{k_*}$.

If $k_*=0$, this follows by the construction of $d_0=\wt d_0$. Assume
that $k_*\geq 1$ and for a contradiction that $\wt d_{k_*}\in C$. By
the minimality of $k_*$, let $w$ be a nonzero element of $(\wt
d_{k_*-1} \,g'R_v^{\;n})\cap W_0$. If $k_*\geq 2$, by Equation
\eqref{eq:inddecreasnonzerovect1}, there exists $w'\in (\wt d_{k_*-2}
\,g'R_v^{\;n})\cap W$ such that $w=d_{k_*-1}\,w'$. By Equation
\eqref{eq:inddecreasnonzerovect2}, we have $\|w\|\geq c\, \|w'\|$. Hence
$\|w'\|\leq\|w\|$ since $c\geq 1$, and $w'$ belongs as $w$ to the ball
$W_0$ centered at $0$. Thus $w'\in (\wt d_{k_*-2} \,g'R_v^{\;n})\cap
W_0$ and $w=d_{k_*-1}\,w'$. By induction, we may therefore write $w=\wt
d_{k_*-1}w_0$ for a nonzero element $w_0$ in $(g'R_v^{\;n})\cap W_0$
such that $\wt d_{j}w_0\in W_0$ for all $j\in\llbracket0,k_*-1
\rrbracket$. By the definition of $W_0$ which contains $w$ and is
contained in $f^{-1}W$ for every $f\in F$ and since $k_*\geq 1$, we
have $d_{k_*}w\in W$.  Hence $\wt d_{k_*} w_0=d_{k_*} \wt
d_{k_*-1}w_0= d_{k_*} w$ belongs to $W$. Since $\wt d_{k_*}\in C$,
this implies that $w_0\in ((C^{-1}W)\cap (g'R_v^{\;n}))\ssm\{0\}=J$.
Since $d_{0}w_0=\wt d_{0}w_0\in W_0\subset W$, this contradicts
Equation \eqref{eq:defid0}.  \cqfd

\bcoro \label{coro:classdivunimod}
%%
%\todo{\tiny faire version pour $\X$ ?}
%%
For every $L\in\X_1$, the following assertions are equivalent.
\begin{enumerate}
\item\label{item1:classdivunimod}
  The orbit map $a\mapsto aL$ from $A_1$ to $\X_1$ is a proper map.
\item\label{item2:classdivunimod}
  The orbit  map $d\mapsto dL$ from $D$ to $\X_1$ is a proper map.
\item\label{item3:classdivunimod}
  There exists $g\in A_1\GL_n^1(K)$ such that $L=g\,R_v^{\;n}$.
\item\label{item4:classdivunimod} The orbit $\wt A \,L$ of $L$ by the
  diagonal subgroup $\wt A$ of $\GL_n(K_v)$ contains an integral
  (possibly nonunimodular) $R_v$-lattice.
\item\label{item5:classdivunimod}
  The orbit $A_1L$ contains an axial (unimodular) $R_v$-lattice.
\item\label{item6:classdivunimod} Every element of the orbit $A_1L$ is
  an axial (unimodular) $R_v$-lattice.
\end{enumerate}
\ecoro

\dem Using the notation of Subsection \ref{subsec:systole}, we have a
canonical onto map $\X_1\ra P\!\X_{1}$ which associates to
a unimodular $R_v$-lattice its homothety class. This map is
equivariant with respect to the canonical morphism $G_1\ra G_1/ZG_1$,
hence with respect to the canonical morphisms $A_1\ra P\!A$ and $D\ra
P\!D$.  The above map $\X_1\ra P\!\X_{1}$ is a proper map,
since its fibers are the compact subsets $\OOO_v^\times L$ for $L\in
\X_1$.

The image of $A_1$ in $P\!A$ is a finite index subgroup, with index
$n$ (and representatives of the classes the elements
$\begin{bmatrix}\pi_v^{-k} & 0 \\ \\ 0 & I_{n-1}\end{bmatrix}$ for
$k\in\llbracket0, n-1\rrbracket$). As seen in Subsection
\ref{subsec:diagosub}, the space $P\!\X$ is the finite union of the
strata $P\!\X_{k}$ for $k\in\intbra{1,n}$ that are transitively
permuted by $P\!A$.

Therefore for every $L\in \X_1$, the orbit map $a\mapsto aL$ from
$A_1$ to $\X_1$ (respectively $d\mapsto dL$ from $D$ to $\X_1$) is a
proper map if and only if the orbit map $[a]\mapsto [a]\,[L]$ from
$P\!A$ to $P\!\X$ (respectively $[d]\mapsto [d]\,[L]$ from $P\!D$ to
$P\!\X$) is a proper map.

The result then follows from Theorem \ref{theo:defidiv}.
\cqfd

\section{A description of the divergent diagonal orbits}
\label{sec:descripdivdiagorb}
  
\subsection{Compact core and quasicenters of divergent diagonal orbits}
\label{subsec:compdivdiagorb}

Let $x\in \X_1$ be a unimodular $R_v$-lattice in $K_v^{\;n}$, which is axial,
or equivalently by Corollary \ref{coro:classdivunimod} such that its
orbit in $\X_1$ under the diagonal subgroup $A_1$ is divergent. In
this subsection, we define and study several invariants associated with
$x$ or with its $A_1$-orbit $A_1x$.

%By Corollary \ref{coro:classdivunimod}, there exists $a\in A_1$ such
%that $ax$ contains an integral sublattice. Since the sum of the
%integral sublattices of $ax$ is still an integral sublattice of $ax$,
%there exists a unique maximal integral sublattice $L_a$ contained in
%$ax$. If $a'\in A_1$ is another element such that $ax$ contains an
%integral sublattice, we have ???

For every $i\in\intbra{1,n}$, we define
\[
\sys_i(x)=\log_{q_v}\min\big\{\|w\|:
w\in (x\cap K_ve_i)\ssm\{0\}\big\}\in\ZZ\;,
\]
that we call the (logarithmic) {\it $i$th-directional systole}. As
seen in Subsection \ref{subsec:lattices}, since $x$ is axial, there
exists $\lambda_i\in K_v^\times$ such that $x\cap K_ve_i=R_v\lambda_i
e_i$, hence we have $\sys_i(x)=\log_{q_v}|\lambda_i|$. We define
\begin{equation}\label{eq:defitruncovol}
\tau_x=\tau_{A_1x}=\sum_{i=1}^n\sys_i(x)\in\ZZ\;,
\end{equation}
that we call the (truncated) {\it covolume} of the divergent orbit
$A_1x$, and that we will use as a complexity for divergent orbits.
Assertion \eqref{item2:propriinvar} of Proposition
\ref{prop:propriinvar} below says that the covolume $\tau_x$ is indeed
an invariant of the $A_1$-orbit of $x$. We will illustrate in
Proposition \ref{prop:treecase} when $n=2$ why we think of $\tau_x$ as
the volume of a canonically truncated divergent orbit $A_1x$.  Let
\begin{equation}\label{eq:defiDeltax}
\Delta^x=\{{\bf k}=(k_1,\ldots, k_n)\in\ZZ^n_0:
\forall\;i\in\intbra{1,n},\quad k_i\geq -\sys_i(x)\}\;,
\end{equation}
which is a finite subset of $\ZZ^n_0$, and let $A_1^x=A_1(\OOO_v)
\exp(\Delta^x)$.  The subset
\begin{equation}\label{eq:deficompactcore}
C_x=C_{A_1x}=A_1^xx=A_1(\OOO_v)\exp(\Delta^x)x
\end{equation}
of the $A_1$-orbit of $x$ is compact and open in $A_1x$ since
$A_1(\OOO_v)$ is a compact-open subgroup of $A_1$, and is called the
{\it compact core} of the divergent orbit $A_1x$. Assertion
\eqref{item2:propriinvar} of Proposition \ref{prop:propriinvar} below
says that the compact core $C_x$ is indeed an invariant of the
$A_1$-orbit of $x$.

The {\it coordinate sublattice} of $x$ is
%%
%\todo{\tiny notation to be seen}
%%
\begin{equation}\label{eq:deficoosublat}
x^{\rm coo}= (x\cap K_v e_1)+\ldots + (x\cap K_v e_n)\,.
\end{equation}
It is indeed an $R_v$-lattice contained in $x$, and $a(x^{\rm
  coo})=(ax)^{\rm coo}$ for every $a\in A_1$. In particular, the
covolume of $x^{\rm coo}$ is constant on the $A_1$-orbit of $x$.
%, by the first two assertions \eqref{item1:propriinvar} and
%\eqref{item2:propriinvar} of Proposition \ref{prop:propriinvar} below.

The {\it quasicenter} of the $A_1$-orbit of $x$ is the unique point
$\wh x\in A_1x$ modulo the left action of $A_1(\OOO_v)$ (see Assertion
\eqref{item4:propriinvar} of Proposition \ref{prop:propriinvar} below
for its existence and uniqueness) such that if $\big(P=
\lfloor\frac{\tau_x}{n}\rfloor, Q=\tau_x-n \lfloor\frac{\tau_x}{n}
\rfloor\big)\in\NN^2$ is the Euclidean division (with $0\leq Q<n$) of
$\tau_x=Pn+Q$ by $n$, then
\begin{equation}\label{eq:defiquasicenter}
\forall\;i\in\intbra{1, Q},\quad \sys_i(\wh x)=P+1
\quad\text{and}\quad
\forall\;i\in\intbra{Q+1,n},\quad \sys_i(\wh x)=P\,.
\end{equation}

For instance (see Proposition \ref{prop:descripxt} for other examples), if
$x=R_v^{\;n}$, we have $\sys_i(x)=0$ for every $i\in\intbra{1,n}$, hence
$x^{\rm coo}=\wh x=x$ and
\[
\tau_x=0,\quad\Delta^x=\{0\}\quad\text{and}\quad C_x=A_1(\OOO_v)x\,.
\]

\bprop\label{prop:propriinvar} Let $x\in\X_1$ be an axial unimodular
$R_v$-lattive, and let $i\in\intbra{1,n}$.
\begin{enumerate}
\item\label{item1:propriinvar} We have $q_v^{\sys_i(x)}=q^{1-\ggg}
  \vol_v\big((K_v e_i)/(x\cap K_v e_i)\big)$ and $\covol(x^{\rm coo})=
  q^{n(\ggg -1)}q_v^{\,\tau_x}$.
\item\label{item2:propriinvar}
  For every $a=\diag(a_1,\ldots,a_n)\in A_1$, we have
  \[
  \sys_i(ax)= \sys_i(x) -v(a_i)\,, \quad \tau_{ax}=\tau_x \quad
  \text{and}\quad C_{ax}=C_{x}\,.
  \]
\item\label{item3:propriinvar} We have $\tau_x\in\NN$. Furthermore
  $\tau_x=0$ if and only if $x=x^{\rm coo}$.
\item\label{item4:propriinvar} There exists a quasicenter $\wh x$ of
  the $A_1$-orbit of $x$, unique modulo the action of $A_1(\OOO_v)$.
\item\label{item5:propriinvar} With $c_n=\frac{1}{(n-1)!}$, as
  $\tau_x\ra+\infty$, we have
  \[
  \card\;\Delta^x=c_n\tau_x^{n-1} +\bigO(\tau_x^{n-2})\,.
  \]
  %, hence $\mu_{A_1x}(C_x)=????c_n\tau_x^{n-1}+\bigO(\tau_x^{n-2})$.
\end{enumerate}
\eprop

\dem Let $x$ and $i$ be as in the statement.

\eqref{item1:propriinvar} As said above, there exists $\lambda_i\in
K_v^\times$ such that $x\cap K_ve_i=R_v\lambda_i e_i$ and we have
$\sys_i(x)=\log_{q_v}|\lambda_i|$. By Equations
\eqref{eq:changevarcovol} and \eqref{eq:covolRv}, we have
\[
\vol_v\big((K_ve_i)/(x\cap K_v e_i)\big)=|\lambda_i|\covol(R_v)
=q_v^{\sys_i(x)}q^{\ggg-1}\,.
\]
The first claim of Assertion \eqref{item1:propriinvar} follows. The
second one follows from the first one and the definition of $\tau_x$.

\medskip
\eqref{item2:propriinvar} We have $\min_{w\,\in\, ax\cap K_ve_i} \|w\|
=|a_i| \min_{w\,\in\, x\cap K_ve_i}\|w\|$. Hence the first claim of
Assertion \eqref{item2:propriinvar} follows since $|a_i|=
q_v^{-v(a_i)}$.  The second claim follows by summation since $|\det a|
=1$.  By the definition \eqref{eq:defiDeltax} of $\Delta^x$ and the
first claim, we have
\[
\Delta^{ax}=\{{\bf k}=(k_1,\ldots, k_n)\in\ZZ^n_0:
\forall\;i\in\intbra{1,n},\quad k_i\geq -\sys_i(x)+v(a_i)\}\,.
\]
For every $i\in\intbra{1,n}$, there exists $a'_i\in\OOO_v^\times$
such that $a_i=a'_i\,\pi_v^{-(-v(a_i))}$. We may hence write $a=
a'\exp{\bf k}'$ with $a'\in A_1(\OOO_v)$ and ${\bf k}' =(-v(a_1),
\ldots, -v(a_n))\in\ZZ^n_0$. Therefore
\begin{equation}\label{eq:varaDeltax}
  \Delta^{ax}=\Delta^{x}-{\bf k}'\;,
\end{equation}
so that since $A_1$ is Abelian, we have
\[
C_{ax}=A_1(\OOO_v)\exp(\Delta^{ax})ax=
A_1(\OOO_v)a'\exp(\Delta^{ax}+{\bf k}')x=A_1(\OOO_v)\exp(\Delta^{x})x
=C_x\,.
\]

\medskip
\eqref{item3:propriinvar} Since the unimodular $R_v$-lattice $x$
contains its coordinate sublattice $x^{\rm coo}$, we have $\covol
(x^{\rm coo})\geq \covol (x)=\covol (R_v^{\;n})$, with equality if and
only if $x=x^{\rm coo}$. Hence by Equation \eqref{eq:covolRv} and by
Assertion \eqref{item1:propriinvar}, we have $q_v^{\,\tau_x}=
\frac{\covol (x^{\rm coo})}{\covol (R_v^{\;n})}\geq 1$.  Therefore
$\tau_x\geq 0$, with equality if and only if $x=x^{\rm coo}$.

\medskip
\eqref{item4:propriinvar} Let $P= \lfloor\frac{\tau_x}{n}\rfloor$ and
$Q=\tau_x-n \lfloor\frac{\tau_x}{n} \rfloor$. Let $a=\exp{\mathbf k}$
where
\[
{\mathbf k}=(-\sys_1(x)+P+1,\ldots,
-\sys_Q(x)+P+1,-\sys_{Q+1}(x)+P,\ldots, -\sys_n(x)+P)\,.
\]
It is easy to check that ${\bf k}\in\ZZ^n_0$ by the definitions of
$\tau_x=\sum_{i=1}^n\sys_i(x)$, and of $P$ and $Q$ so that $\tau_x=
nP+Q$. By the first claim of Assertion \eqref{item2:propriinvar} and
by Equation \eqref{eq:defiquasicenter}, the element $\wh x=ax$ is a
quasicenter of $A_1x$.  If $\wh {\wh x}$ is another quasicenter of
$A_1x$, if $a\in A_1$ is such that $\wh {\wh x}=a\wh x$, then for
every $i\in\intbra{1,n}$, we have $|a_i|=\frac{q_v^{\sys_i(\wh {\wh
      x})}} {q_v^{\sys_i(\wh x)}}$ by the first claim of Assertion
\eqref{item2:propriinvar}. Hence $|a_i|=1$ by the definition
\eqref{eq:defiquasicenter} of the quasicenter and since $\tau_x$ is
constant on the $A_1$-orbit of $x$ by Assertion
\eqref{item2:propriinvar}.  Therefore $a\in A_1(\OOO_v)$, thus proving
the uniqueness claim.

\medskip
\eqref{item5:propriinvar} For every $m\in\NN$, let $\wh\Delta(m)=
\{{\bf k} =(k_1,\ldots, k_n)\in\ZZ^n_0: \forall\;i\in\intbra{1,n},
\quad k_i\geq -m\}$. By Equation \eqref{eq:varaDeltax} and the
definition of the quasicenter, we have
\[
\card\big(\wh\Delta\big(\big\lfloor\frac{\tau_x}{n}\big\rfloor
\big)\big)\leq \card(\Delta^x)=\card(\Delta^{\wh x})\leq
\card\big(\wh\Delta\big(\big\lfloor\frac{\tau_x}{n}\big\rfloor+1
\big)\big)\,.
\]

Let us prove that if $c'_n=\frac{n^{n-1}}{(n-1)!}$, then, as
$m\ra+\infty$, we have
\begin{equation}\label{eq:calcsomint}
\card(\wh\Delta(m))= c'_n m^{n-1} +\bigO(m^{n-2})\,.
\end{equation}
This implies  Assertion \eqref{item5:propriinvar} with $c_n=
\frac{c'_n}{n^{n-1}}=\frac{1}{(n-1)!}$. 

We start the proof by the following elementary integral computation.
We consider the Euclidean subspace $\RR_0^{\,n}=\{ {\bf t}=
(t_1,\ldots, t_n)\in\RR^n: \sum_{i=1}^n t_i = 0\}$ of the standard
Euclidean space $\RR^n$, endowed with its Lebesgue measure
$\Leb_{\RR_0^{\,n}}$. By invariance under translation, this measure is
proportional to the measure $dt_1\ldots dt_{n-1}$ on $\RR_0^{\,n}$,
and the proportionality constant is classically computed as
follows. Let $u=\frac{1}{\sqrt{n}}(1,1,\ldots,1)$ which is a unit
normal vector to the hyperplane $\RR_0^{\,n}$.  Let $P$ be the
fundamental polytope of the $\ZZ$-lattice $\ZZ_0^{\,n}$ in
$\RR_0^{\,n}$ generated by the vectors
\[
u_1=(1,-1,0,\ldots, 0), \quad u_2=(0,1,-1,0,\ldots, 0),
\;\ldots\;,\quad u_{n-1}=(0,\ldots, 0,1,-1)\,.
\]
Note that the first $n-1$ coordinates of a point $s_1u_1+\ldots+
s_{n-1} u_{n-1}$ of the polytope $P$, with $(s_1,\ldots, s_{n-1})
\in[0,1]^{n-1}$, are $t_1=s_1$, $t_2=s_2-s_1$, $\ldots$, $t_{n-1}=
s_{n-1}-s_{n-2}$, so that $dt_1\ldots dt_{n-1}= ds_1\ldots ds_{n-1}$
and $dt_1\ldots dt_{n-1}(P)=1$.  Therefore
\[
\frac{d\Leb_{\RR_0^{\,n}}}{dt_1 \ldots dt_{n-1}}
=\frac{\Leb_{\RR_0^{\,n}}(P)}{dt_1 \ldots dt_{n-1}(P)} =
\Leb_{\RR_0^{\,n}}(P)=|\det(u_1,\ldots,u_{n-1},u)|=\sqrt{n}\,.
\]
Note for future use that
\begin{equation}\label{eq:covolZ0n}
  \covol(\RR_0^{\,n}/\ZZ_0^{\,n})=\Leb_{\RR_0^{\,n}}(P)=\sqrt{n}\,.
\end{equation}
For every $\alpha >0$, let $\wh\Delta'(\alpha) =\{ {\bf t}=
(t_1,\ldots,t_{n})\in \RR_0^{\,n} : \forall\; i \in \intbra{1,n},
\;t_i \geq - \alpha\}$.

\blemm\label{lem:volreel}
For every $\alpha >0$, we have $\Leb_{\RR_0^{\,n}}
(\wh\Delta'(\alpha))= \sqrt{n}\;
\frac{(n\alpha)^{n-1}}{(n-1)!}$.
\elemm

\dem Up to using an homothety of ratio $\alpha$, we may assume by
homogeneity that $\alpha=1$.  For every $k \in \intbra{1,n}$, using
the standard conventions that $*^0=1$ and $\sum_{\emptyset}*=0$, we
define a map $g_k:\RR^{n-k}\ra \RR$ by
\[
(t_1,\ldots,t_{n-k}) \mapsto
\frac{1}{(k-1)!}\Big(k-\sum_{i=1}^{n-k}t_i\Big)^{k-1}\,.
\]
Note that $g_1=1$ and $g_n=\frac{n^{n-1}}{(n-1)!}$ are constant. For
all $k\in\intbra{1,n-1}$ and $t_1,\ldots,t_{n-k-1}\in\RR$, by a
straightforward integration, we have
\begin{align*}
\int_{-1}^{k-\sum_{i=1}^{n-k-1}t_i} g_k(t_1,\ldots,t_{n-k-1},s)\;ds & =
\frac{1}{(k-1)!} \Big[ -\frac{1}{k} \Big(k -\sum_{i=1}^{n-k-1} t_i
  -s\Big)^k \;\Big]_{s=-1}^{s=k - \sum_{i=1}^{n-k-1}t_i} \\ & =
g_{k+1}(t_1,\ldots,t_{n-k-1})\,.
\end{align*}
We have $t_i \geq -1$ and $\sum_{i=1}^n t_i = 0$ for every $i\in
\intbra{1,n}$ if and only if  $-1\leq t_i \leq
n-i-\sum_{j=1}^{i-1}t_j$ for every $i\in \intbra{1,n-1}$. Hence by
an easy induction, we have, for every $k\in \intbra{1,n-1}$,
\begin{align*}
\Leb_{\RR_0^{\,n}}(\wh\Delta'(1)) &= \sqrt{n} \int_{-1}^{(n-1)}
\int_{-1}^{(n-2) - t_1} \ldots \int_{-1}^{1-\sum_{i=1}^{n-2}t_i}
dt_{n-1} \ldots dt_2 \, dt_1\\ & =\sqrt{n} \int_{-1}^{(n-1)}
\int_{-1}^{(n-2) - t_1} \ldots \int_{-1}^{k-\sum_{i=1}^{n-k-1}t_i}
g_k(t_1,\ldots,t_{n-k}) \,dt_{n-k} \ldots dt_2 \, dt_1\,.
\end{align*}
When $k=n-1$, we get $ \Leb_{\RR_0^{\,n}}(\wh\Delta'(1)) = \sqrt{n}\,
\int_{-1}^{(n-1)} g_{n-1}(t_1) \, dt_1 =\sqrt{n}\; g_n$, as
wanted.  \cqfd

\medskip
By the standard Gauss counting argument, by Lemma \ref{lem:volreel}
and by Equation \eqref{eq:covolZ0n}, Equation \eqref{eq:calcsomint}
follows since
\[
\card(\wh\Delta(m))\sim\frac{\Leb_{\RR_0^{\,n}}
  (\wh\Delta'(m))}{\covol(\RR_0^{\,n}/\ZZ_0^{\,n})}=
\frac{(nm)^{n-1}}{(n-1)!}\,.\;\;\;\Box
\]

\subsection{The mass behavior of the compact cores of divergent
  diagonal orbits}
\label{subsec:massbehavedivdiaorb}
  
In this subsection, we prove that for continuous functions with
support in a fixed compact subset of $\X_1$, most of their mass for
the orbital measure $\ov\mu_{x}=\ov\mu_{A_1x}$ (defined in Subsection
\ref{subsec:homogdiag}) on a divergent orbit $A_1x$ is carried by the
compact core of $A_1x$ as the truncated covolume goes to infinity.

We keep denoting by $x$ an axial unimodular $R_v$-lattice. We denote by
%
%\todo{\tiny voir si on note seulement $\ov\nu_x$  et $\nu_x$}
%
\begin{equation}\label{eq:defovnux}
\ov\nu_x=\ov\nu_{A_1x}
=\frac{1}{\ov\mu_{x}(C_{x})}\;\ov\mu_{x}\!\mid_{C_{x}}
\end{equation}
the restriction of the orbital measure $\ov\mu_{x}$ to the compact
core $C_{x}=C_{A_1x}$ of the divergent orbit $A_1x$, normalized to be
a probability measure on $\X_1$. It is well defined since $C_{x}$ is a
nonempty compact open subset of $A_1x$, hence $0< \ov\mu_{x}(C_{x}) <
+\infty$.  It is independent of the choice $x$ of an element in the
orbit $A_1x$, and its support is equal to $C_{x}$.  By Equation
\eqref{eq:deficompactcore}, $C_{x}=A_1(\OOO_v)\exp(\Delta^x)\,x$ is
the disjoint union of the clopen subsets $A_1(\OOO_v)(\exp{\bf k})
\,x$ for ${\bf k}\in \Delta^x$. By the normalisation of the Haar
measure of $A_1$ in Equation \eqref{eq:volA1OOOv}, we have
\begin{equation}\label{eq:volconcor}
  \ov\mu_{x}(C_{x})={\tt m}_{A_1}(A_1(\OOO_v)\exp(\Delta^x))
  =\card(\Delta^x)\,.
\end{equation}
This formula, paired with Assertion \eqref{item5:propriinvar} of
Proposition \ref{prop:propriinvar}, says that up to an error term, up
to a multiplicative constant and up to a power constant depending only
on $n$, the truncated covolume $\tau_x$ is the orbital measure of the
compact core of the divergent orbit $A_1x$, again justifying its
name. With the simplified notation of Subsection
\ref{subsec:homogdiag}, for $a\in A_1(\OOO_v)$ and
$\mb{k}\in\Delta^x$, we have
\begin{equation}\label{eq:expressimpnu_x}
d\ov\nu_x(a\exp(\mb{k})x)=
\frac{1}{\card(\Delta^x)}\;da\;d\mb{k}\!\mid_{\Delta^x}\,.
\end{equation}

If $y\in\X$, then $Ay$ is divergent in $\X$ if and only if $A_1y$ is
divergent in $\X_1$, and we then similarly denote by
\[
\nu_y=\nu_{Ay}
=\frac{1}{\mu_{y}(C_{y}\cap Ay)}\;\mu_{y}\!\mid_{C_{y}\cap Ay}
\]
the restriction of the orbital measure $\mu_{y}$ to the compact core
$C_{y}\cap Ay=A(\OOO_v)\exp(\Delta^y)y$ of the divergent orbit $Ay$,
normalized to be a probability measure on $\X$.  For $a\in A(\OOO_v)$ and
$\mb{k}\in \Delta^y$, again with the simplified notation of Subsection
\eqref{subsec:homogdiag} now for $A(\OOO_v)$, we have
\begin{equation}\label{eq:calcpratnuy}
  d\nu_y(a\exp(\mb{k})y)=
  \frac{1}{\card(\Delta^y)}\;da\;d\,\mb{k}\!\mid_{\Delta^y}\,.
\end{equation}
By Equation \eqref{eq:decompdmA1}, the measure $\ov\nu_y$ on $\X_1$ is
an average of normalized restrictions of orbital measures on $\X$~:
\begin{equation}\label{eq:desintegovnunu}
  \ov\nu_y=\frac{q_v}{q_v-1}\;
  \int_{\lambda\in\OOO_v^{\,\times}}\xi(\lambda)_*\nu_y\;d\vol_v(\lambda)\,.
\end{equation}

In the next lemma, we denote by $\|\;\|_\infty$ the uniform norm of
continuous functions with compact support.

\blemm\label{lem:massmostoncompcore} For every $m\in\NN$, for every
continuous function $f\in {\rm C}_c(\X_1)$ with compact support
contained in $\X_1^{\geq q_v^{-m}}$, for every axial element $x\in\X_1$, as 
$\tau_x\ra+\infty$, we have 
\[
\frac{1}{c_n\tau_x^{n-1}}\;\ov\mu_{x}(f)=\ov\nu_{x}(f)+
\bigO\Big(\frac{(m+1)\,\|f\|_\infty}{\tau_x}\Big)\,.
\]
Similarly, for all $y\in\X$ axial and $f\in {\rm C}_c(\X)$ with
support in $\X^{\geq q_v^{-m}}$, as $\tau_y\ra+\infty$, we have
\[
\frac{1}{c_n\tau_y^{n-1}}\;\mu_{y}(f)=\nu_{y}(f)+
\bigO\Big(\frac{(m+1)\,\|f\|_\infty}{\tau_y}\Big)\,.
\]
\elemm

\dem Let $m,f,x$ be as in the first statement (the second one is
similar). We define
\[
\Delta^{x,m}=\big\{{\bf k}=(k_1,\ldots, k_n)\in\ZZ^n_0:
\forall\;i\in\intbra{1,n},\quad k_i\geq -\sys_i(x)-m\big\}\;,
\]
$A_1^{x,m}=A_1(\OOO_v)\exp(\Delta^{x,m})$ and $C_{x,m}=A_1^{x,m}\,x$,
which respectively contain $\Delta^x$, $A_1^x$ and $C_x$. By Equation
\eqref{eq:varaDeltax} and the definition of the quasicenter, as in the
proof of Proposition \ref{prop:propriinvar} \eqref{item5:propriinvar},
we have
\[
\card\big(\wh\Delta\big(\big\lfloor\frac{\tau_x}{n}\big\rfloor+m
\big)\big)\leq \card(\Delta^{x,m})=\card(\Delta^{\wh x,m})\leq
\card\big(\wh\Delta\big(\big\lfloor\frac{\tau_x}{n}\big\rfloor+m+1
\big)\big)\;.
\]
Hence by Equation \eqref{eq:calcsomint}, as $\tau_x\ra+\infty$, we have
\begin{equation}\label{eq:controlDeltaxm}
  \card(\Delta^{x,m})=c_n(\tau_x+nm)^{n-1}
  +\bigO((\tau_x+nm)^{n-2})\,.
\end{equation}

Note that if $a=\diag(a_1,\ldots,a_n)\in A_1\ssm A_1^{x,m}$, then
there exist $i\in\intbra{1,n}$, $k_i\in\ZZ$ and
$a'_i\in\OOO_v^{\,\times}$ such that $a_i=a'_i\pi_v^{-k_i}$ and
$k_i=-v(a_i)<-\sys_i(x)-m$. Hence since $x$ is unimodular, by the
definition \eqref{eq:defisys} of the systole in Subsection
\ref{subsec:systole} and of the logarithmic directional systoles in
this section, and by Proposition \ref{prop:propriinvar}
\eqref{item2:propriinvar}, we have
\[
\log_{q_v}\sys(ax)\leq \min_{1\leq j\leq n}\sys_j(ax)=
\min_{1\leq j\leq n}(\sys_j(x) -v(a_j))\leq \sys_i(x) -v(a_i)<-m\,.
\]
Thus $ax\notin \X_1^{\geq q_v^{-m}}$ and $f(ax)=0$. Therefore, using

$\bullet$~ Equations \eqref{eq:defovnux} and \eqref{eq:volconcor} for
the first and second equalities (and similarly for $C_{x,m}$),

$\bullet$~ Equation \eqref{eq:controlDeltaxm} for the third equality,

$\bullet$~ with a $\bigO(\;)$ which depends only on $n$ for
$\tau_x\geq nm$ for the last equality,

\noindent we have
\begin{align*}
  &\Big|\,\frac{1}{c_n\tau_x^{n-1}}\;\ov\mu_{x}(f)-\ov\nu_{x}(f)\,\Big|
  =\Big|\,\Big(\frac{1}{c_n\tau_x^{n-1}}\ov\mu_{x}\!\mid_{C_{x,m}}-
  \frac{1}{\card(\Delta^{x})}\ov\mu_{x}\!\mid_{C_{x}}\Big)(f)\,\Big|
  \\\leq \;&  \Big(\frac{\ov\mu_{x}(C_{x,m})}{c_n\tau_x^{n-1}}-
  \frac{\ov\mu_{x}(C_{x})}{\card(\Delta^{x})}\Big)\|f\|_\infty
  =\Big(\frac{\card(\Delta^{x,m})}{c_n\tau_x^{n-1}}- 1\Big)\|f\|_\infty
  \\= \;&  \Big(\frac{c_n(\tau_x+nm)^{n-1}
    +\bigO((\tau_x+nm)^{n-2})}{c_n\tau_x^{n-1}}- 1\Big)\|f\|_\infty
  \\= \;&  \Big(\big(1+\frac{nm}{\tau_x}\big)^{n-1}
     +\bigO\big((1+\frac{nm}{\tau_x})^{n-2}\frac{1}{\tau_x}\big)
     - 1\Big)\|f\|_\infty=\bigO\Big(\frac{(m+1)\,\|f\|_\infty}{\tau_x}\Big),
\end{align*}
as wanted.
\cqfd

\subsection{Zigzag length and continued fractions}
\label{subsec:zigzag}
  
We assume in this whole subsection that $n=2$ (and we will then use
$n$ as a variable element of $\NN$), that $K=\FF_q(Y)$ and that $\deg
v=1$ so that $R_v=\FF_q[Y]$ (see Subsection
\ref{subsec:functionfield}). We give in this particular case a
geometric interpretation (using the geodesic flow on a Bruhat-Tits
tree) and an arithmetic interpretation (using continued fraction
expansions) of the quantities defined in Subsection
\ref{subsec:compdivdiagorb}. We refer for instance to \cite[\S 15.1,
  15.2]{BroParPau19} for the background information.
%%
%\todo{case $n=3$ ?}
%%

Let $\TT_v$ be the {\it Bruhat-Tits tree} of $(\PGL_2,K_v)$, see for
instance \cite{Serre83}. It is a regular tree of degree $\card
\;\PP_1(\OOO_{v}/\pi_v\OOO_{v})= q+1$ (since $q_v=q$ here) and its set
of vertices $V\TT_v$ is the set of homothety classes (under
$K_v^{\,\times}$) $[\Lambda]$ of $\OOO_{v}$-lattices $\Lambda$ in $K_v
\times K_v$. We denote by $*_v$ the homothety class of the
$\OOO_{v}$-lattice $\OOO_{v} \times\OOO_{v}$ generated by the
canonical basis of $K_v\times K_v$. The left linear action of
$G_1=\GL_2^1(K_v)$ on $K_v\times K_v$ induces a left action of $G_1$
on $\TT_v$, which preserves and is transitive on the set $V_{\rm even}
\TT_v$ of vertices at even distance from $*_v$.
%The stabilizer in $G_1$ of $*_v$
%is $G_1(\OOO_v)= \GL_2^1(K_v) \cap \GL_2(\OOO_v)$.
Let $\aaa=\begin{pmatrix} \pi_v & 0 \\ 0 &
\pi_v^{\;-1} \end{pmatrix}\in D=\wt D\cap\SL_2(R_v)$, which generates
$D\simeq \ZZ$. The lattice $\Ga_1=\GL_2(R_v)$ in $G_1$ (or its
projective version $P\Ga$) is called the {\it Nagao lattice}, see
\cite{Nagao59,Weil70}.

We identify the projective line $\PP_1(K_v)$ with $K_v \cup
\{\infty\}$ using the map $[x:y] \mapsto xy^{-1}$ as usual, and we
endow $\PP_1(K_v)$ with the projective action of $G_1$. The boundary
at infinity $\partial_\infty\TT_v$ of $\TT_v$ identifies
$G_1$-equivariantly with $\PP_1(K_v)$. The Nagao lattice $\Ga_1$ acts
transitively on the subset $\PP_1(K)$ of $\PP_1(K_v)$.

We denote by $\G\TT_v$ the space of {\it geodesic lines} in $\TT_v$
(that is, the set of isometric maps $\ell:\ZZ\ra V\TT_v$ endowed with
the compact-open topology), endowed with the action by
post-composition of $G_1$ defined by $(g,\ell)\mapsto
\{g\,\ell:k\mapsto g\,\ell(k)\}$. Let
\[
\G_{\rm even}\TT_v=\{\ell\in \G\TT_v:\ell(0)\in V_{\rm even} \TT_v\}\,,
\]
which is invariant by $G_1$. Let $\ell_*\in\G\TT_v$ be the unique
geodesic line with $\ell_*(-\infty) =\infty\in\partial_\infty\TT_v$,
$\ell_*(+\infty) = 0\in\partial_\infty\TT_v$, and $\ell_*(0)=*_v$, so
that $\ell_*(2n)= \aaa^{n}\,*_v$ for every $n\in\ZZ$.

A geodesic line $\ell$ in $\G\TT_v$, as well as its image in
$\Ga_1\bs\G\TT_v$, is called {\it birational} if its two points at
infinity $\ell(\pm\infty)$ belong to $\PP_1(K)$.  For instance,
$\ell_*$ belongs to $\G_{\rm even}\TT_v$ and is birational.

We denote by $(\phi^n)_{n\in\ZZ}$ the (discrete-time) {\it geodesic
  flow} on the space $\G\TT_v$, defined by $\phi^n\ell:k\mapsto
\ell(k+n)$ for all $n\in\ZZ$ and $\ell \in \G\TT_v$, which commutes
with the action of $G_1$, as well as its quotient flow on
$\Ga_1\bs\G\TT_v$.  The stabilizer of $\ell_*$ for the transitive
action of $G_1$ on $\G_{\rm even} \TT_v$ is exactly $A_1(\OOO_v)$.
Hence the map $A_1(\OOO_v)g\mapsto g^{-1}\,\ell_*$ is a homeomorphism
$\wt\Xi$ from $A_1(\OOO_v)\bs G_1$ to $\G_{\rm even}\TT_v$, which is
(anti-)equivariant with respect to the actions of $\Ga_1$ on the right
on $A_1(\OOO_v)\bs G_1$ and on the left on $\G_{\rm even}\TT_v$~:
\[
\forall\; g\in G_1,\;\;\forall\;\ga\in\Ga_1,\quad
\wt\Xi (A_1(\OOO_v)g\ga)=\ga^{-1}\wt\Xi (A_1(\OOO_v)g)\,.
\]
We denote by $\Xi:A_1(\OOO_v)\bs G_1/\Ga_1\ra\Ga_1\bs\G_{\rm even}
\TT_v$ the homeomorphism induced by $\wt\Xi$. We have the following
crucial property relating the right action of $A_1$ (or the commuting
right actions of $D$ and $A_1(\OOO_v)$) on $\X_1=G_1/\Ga_1$ and the
even-time geodesic flow on $\Ga_1\bs\G_{\rm even}\TT_v$: for every
$n\in\ZZ$, we have $\wt\Xi\circ\aaa^{-n}= \phi^{2n}\circ \wt\Xi$, or
equivalently
\begin{equation}\label{eq:commuta0geodflow}
  \forall\;g\in G_1,\quad\phi^{2n}(g\,\ell_*)= g\,\aaa^{n}\,\ell_*\,.
\end{equation}
Refering to \cite{Serre83} for background, the quotient graph of
groups $\Ga_1\dbs\TT_v$ is the following {\it modular ray}

\begin{center}
\begin{picture}(0,0)%
\includegraphics{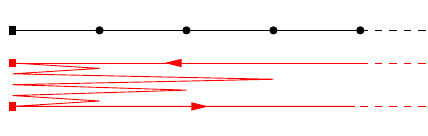}%
\end{picture}%
\setlength{\unitlength}{3812sp}%
\begingroup\makeatletter\ifx\SetFigFont\undefined%
\gdef\SetFigFont#1#2#3#4#5{%
  \reset@font\fontsize{#1}{#2pt}%
  \fontfamily{#3}\fontseries{#4}\fontshape{#5}%
  \selectfont}%
\fi\endgroup%
\begin{picture}(3537,1107)(256,-481)
\put(3061,-421){\makebox(0,0)[lb]{\smash{{\SetFigFont{11}{13.2}{\rmdefault}{\mddefault}{\updefault}{\color[rgb]{0,0,0}$\Ga_1  \ell$}%
}}}}
\put(271,479){\makebox(0,0)[lb]{\smash{{\SetFigFont{11}{13.2}{\rmdefault}{\mddefault}{\updefault}{\color[rgb]{0,0,0}$H_{-1}$}%
}}}}
\put(991,479){\makebox(0,0)[lb]{\smash{{\SetFigFont{11}{13.2}{\rmdefault}{\mddefault}{\updefault}{\color[rgb]{0,0,0}$H_0$}%
}}}}
\put(631,209){\makebox(0,0)[lb]{\smash{{\SetFigFont{11}{13.2}{\rmdefault}{\mddefault}{\updefault}{\color[rgb]{0,0,0}$H'_0$}%
}}}}
\put(1711,479){\makebox(0,0)[lb]{\smash{{\SetFigFont{11}{13.2}{\rmdefault}{\mddefault}{\updefault}{\color[rgb]{0,0,0}$H_1$}%
}}}}
\put(1351,209){\makebox(0,0)[lb]{\smash{{\SetFigFont{11}{13.2}{\rmdefault}{\mddefault}{\updefault}{\color[rgb]{0,0,0}$H_0$}%
}}}}
\put(2431,479){\makebox(0,0)[lb]{\smash{{\SetFigFont{11}{13.2}{\rmdefault}{\mddefault}{\updefault}{\color[rgb]{0,0,0}$H_2$}%
}}}}
\put(2071,209){\makebox(0,0)[lb]{\smash{{\SetFigFont{11}{13.2}{\rmdefault}{\mddefault}{\updefault}{\color[rgb]{0,0,0}$H_1$}%
}}}}
\put(2791,209){\makebox(0,0)[lb]{\smash{{\SetFigFont{11}{13.2}{\rmdefault}{\mddefault}{\updefault}{\color[rgb]{0,0,0}$H_2$}%
}}}}
\end{picture}%

\end{center}

\noindent 
where $H_{-1}=\Ga_1\cap\GL_2(\FF_q)$, $H'_0=H_0\cap H_{-1}$
and, for every $n\in\NN$,
\[
H_n=\bigg\{\begin{pmatrix}a&b\\0&d\end{pmatrix}
\in\Ga_1\;:\; a,d\in \FF_q^\times, b\in \FF_q[Y], \deg b\leq n+1\bigg\}\,.
\]
%The stabilizer of a vertex $s$ in $\Ga_1\dbs\TT_v$, that is different
%from the origin $\Ga_1\,*$ of the ray $\Ga_1\bs\TT_v$, is equal to the
%stabilizer of the edge starting from $s$ and pointing towards infinity
%in the ray $\Ga_1\bs\TT_v$. Hence a geodesic line in $\Ga_1\bs\G\TT_v$
%that starts to go down towards the origin $\Ga_1\,*$ goes down all the
%way to $\Ga_1\,*$.
A birational geodesic line in $\Ga_1\bs\G\TT_v$ starts from the point
at infinity of the ray $\Ga_1\bs\TT_v$, goes down to the origin
$\Ga_1\,*_v$, then makes some back-and-forth to the origin for some
(even, possibly zero) finite time, then goes up all the way to the
point at infinity of the ray (see \cite[page 116]{Serre83},
\cite[\S 6.1]{Paulin02} and the above picture).

Half the (even) length of the birational geodesic line $\Ga_1\ell$ in
$\Ga_1\bs\G_{\rm even}\TT_v$ between the first and last time of
passage through the origin $\Ga_1\,*_v$ is called the {\it zigzag
  length} of $\Ga_1\ell$, and denoted by $\operatorname{zz}(\Ga_1\ell)
\in \NN$.  It is invariant under the action of the geodesic flow. For
instance, $\operatorname{zz}(\Ga_1\ell_*)=0$.

Any element $f\in K_v=\FF_q((Y^{-1}))$ may be uniquely written as a
sum $f=[f]+\{f\}$ with $[f]\in R_v=\FF_q[Y]$ (called the {\it integral
  part} of $f$) and $\{f\}\in \pi_v\OOO_v$ (called the {\it fractional
  part} of $f$).  The {\it Artin map} $\Psi: \pi_v\OOO_v\ssm
\{0\}\ra \pi_v\OOO_v$ is defined by $f\mapsto \big\{\frac{1}{f}
\big\}$.  Any $f\in K=\FF_q(Y)$ has a unique finite {\it continued
  fraction expansion}
\[
f=[a_0;a_1,\ldots, a_n]=a_0 + \cfrac{1}{a_1+\cfrac{1}{a_2+
    \cfrac{1}{\cdots+\cfrac{1}{a_n}}}}=\frac{P_n}{Q_n}\;,
\]
with $a_0=a_0(f)=[f]\in R_v$ and if $f\neq a_0$ then $n=n(f) \in\NN
\ssm\{0\}$ is such that we have $\Psi^{n}(f-a_0)= 0$ and $a_i= a_i(f)=
\big[\frac{1}{\Psi^{i-1}(f-a_0)} \big]$ is a nonconstant polynomial
for $1\leq i\leq n$. The elements $a_0,a_1,\ldots, a_n\in R_v$ are
called the {\it coefficients} of the continued fraction expansion of
$f$. The fraction $\frac{P_i}{Q_i}= [a_0; a_1,\ldots, a_i]\in K$ is
called the {\it $i$-th convergent} of $f$ and is uniquely defined by
induction by
\begin{equation}\label{eq:recurconvergent}
\begin{array}{cccc}
P_{-1}=1 & P_{0}=a_{0}, &   & P_{i}=P_{i-1}a_{i}+P_{i-2}\\
Q_{-1}=0 & Q_{0}=1, &   & Q_{i}=Q_{i-1}a_{i}+Q_{i-2}
\end{array}
\end{equation}
for $1\leq i\leq n$. We refer for instance to \cite{Lasjaunias00,
  Schmidt00, Paulin02} for background on the above notions.

The stabilizer of $\infty\in\PP_1(K_v)$ for the projective action of
$\Ga_1$ is its upper triangular subgroup $H_\infty =\bigcup_{n\in\NN}
H_n$. For every $f\in K_v=\PP_1(K_v)\ssm \{\infty\}$, there exists
$g\in H_\infty$ such that $gf\in \pi_v\OOO_v$ (for instance $g=
\big(\begin{smallmatrix}1&-[f]\\0&1\end{smallmatrix}\big)$ ) and if
$g'\in H_\infty$ also satisfies that $g'f\in \pi_v\OOO_v$, then there
exists $u\in \FF_q^{\,\times}$ such that $g'f= u(gf)$. Hence every
birational geodesic line $\ell$ in $\G\TT_v$ has a representative $\wt
\ell$ in its class $\phi^\ZZ\Ga_1 \ell$ modulo the (commuting) actions
of the geodesic flow and of $\Ga_1$ which starts at time $-\infty$
from $\infty\in\partial_\infty \TT_v$, passes at time $t=0$ through
$*\in V\TT_v$, and ends at a point in
$\pi_v\OOO_v\subset\partial_\infty \TT_v$, unique up to multiplication
by an element of $\FF_q^\times$. Note that $t=0$ is the time when $\wt
\ell$ starts its zigzag part (see Proposition \ref{prop:treecase} for
the computation of the time $\wt \ell$ ends its zigzag part). We
define the {\it continued fraction total length}
$\operatorname{cf}(\Ga_1\ell)$ of $\Ga_1\ell$ as the sum of the
degrees of the coefficients of the continued fraction expansion
$[0;a_1,a_2\ldots, a_n]$ of $\wt \ell(+\infty)$~:
\[
\operatorname{cf}(\Ga_1\ell)=\sum_{i=1}^n\deg(a_i)\,.
\]
This does not depend on the choice of $\wt \ell$, since for all $u\in
\FF_q^\times$ and $a_1,\ldots, a_n\in R_v\ssm \FF_q$, we have
$u[0;a_1,a_2\ldots, a_n]= [0;u^{-1}a_1,ua_2,\ldots, u^{(-1)^n}a_n]$.  We
define the {\it height} $\operatorname{ht}(\Ga_1\ell)$ of $\Ga_1\ell$
as the degree of the denominator of the last convergent
$\frac{P_n}{Q_n}$ of $\wt \ell(+\infty)$~:
\[
\operatorname{ht}(\Ga_1\ell)=\deg(Q_n)\,.
\]
The following result says that the truncated covolume of a divergent
orbit in $G_1/\Ga_1$ coincides with the zigzag length, with the
continued fraction total length and with the height of the
corresponding orbit of the even-time geodesic flow in $\Ga_1\bs\G_{\rm
  even} \TT_v$.

\bprop\label{prop:treecase} For every $g\in G_1$, the $A_1$-orbit
$A_1gR_v^{\,2}$ of the $R_v$-lattice $gR_v^{\,2}$ is divergent in
$\X_1$ if and only if the geodesic line $\Xi(A_1(\OOO_v)g\Ga_1)=
\Ga_1g^{-1}\ell_* \in\Ga_1\bs\G_{\rm even}\TT_v$ is birational, and we
then have
\begin{equation}\label{eq:quadrupleq}
\tau_{A_1g R_v^{\,2}}=\operatorname{ht}(\Ga_1g^{-1}\ell_*)
=\operatorname{cf}(\Ga_1g^{-1}\ell_*)
=\operatorname{zz}(\Ga_1g^{-1}\ell_*)\,.
\end{equation}
\eprop

\dem Let $g\in G_1$. By Corollary \ref{coro:classdivunimod}, we know
that $A_1gR_v^{\,2}$ is divergent if and only if $g\in A_1\GL_2^1(K)$.

The group $\GL_2^1(K)$ acts transitively on the set of ordered pairs
of distinct points of $\PP_1(K)$, since for all $x,y\in K$, the
element $\begin{psmallmatrix} 0 & 1\\ 1 & -x \end{psmallmatrix}\in
\GL_2^1(K)$ maps $x$ to $\infty$ and the element $\begin{psmallmatrix}
  1 & -y\\ 0 & 1 \end{psmallmatrix}\in\GL_2^1(K)$ maps $y$ to $0$
while fixing $\infty$.  Hence if the geodesic line $g^{-1}\ell_*\in
\G_{\rm even}\TT_v$ is birational, then there exist $h\in \GL_2^1(K)$
and $n\in\NN$ such that $hg^{-1}\ell_*=\phi^{2n}\ell_* =\aaa^{n}
\ell_*$, using Equation \eqref{eq:commuta0geodflow} for the last
equality.  Hence $\aaa^{-n}hg^{-1}$ fixes $\ell_*$, whose stabilizer
is $A_1(\OOO_v)$.  Therefore $g\in A_1\GL_2^1(K)$.

Conversely, assume that $g\in A_1\GL_2^1(K)$. Since $A_1=
D A_1(\OOO_v)$, there exist $n\in\ZZ$, $h'\in A_1(\OOO_v)$ and $h\in
\GL_2^1(K)$ such that $g=\aaa^{n}h'h$. The points at infinity of
$g^{-1}\ell_*$ are hence equal to the points at infinity of
$h^{-1}\ell_*$, which are both in $\PP_1(K)$, hence $g^{-1}\ell_*$ is
birational. This proves the first claim.

Let us prove the first equality of Equation \eqref{eq:quadrupleq}. If
$A_1g R_v^{\,2}$ is divergent, we may assume that $g\in \GL_2^1(K)$ by
Corollary \ref{coro:classdivunimod}.  By the transitivity properties
of $\Ga_1$, up to multiplying $g$ on the left by an element of
$A_1\cap\GL_2^1(K)$ and on the right by an element of $\Ga_1$, we may
assume that $g^{-1}*_v=*_v$ and that the projective action of $g^{-1}$
fixes $\infty$ and sends $0$ to the last convergent $\frac{P_n}{Q_n}$
of $g^{-1}\ell_*(+\infty)$. Thus $g$ has the form
$\begin{psmallmatrix} a & b\\ 0 & d \end{psmallmatrix}$ with $a,b,d\in
K$ with $|ad|=1$. In particular, we have $g^{-1}\ell_*=
\wt{g^{-1}\ell_*}$ with the above notation. Since multiplying $g$ on
the left by an element of $A_1\cap \Ga_1$ does not change
$g^{-1}\ell_*(+\infty)$, we may assume that $a=d=1$. Then $b=
-\frac{P_n}{Q_n}$. Now, we have $g R_v^{\,2}= R_ve_1+ (be_1+e_2) R_v$.
Hence $g R_v^{\,2}\cap K_v e_1=R_v e_1$ and $g R_v^{\,2}\cap K_v e_2=
R_v Q_n e_2$. Thus by Equation \eqref{eq:inversiRv} and the definition
of the directional systoles, we have
\[
\sys_1(g R_v^{\,2})=\log_{q}1=0\quad\text{and} \quad
\sys_2(g R_v^{\,2})=\log_{q}|Q_n|=\deg Q_n\;,
\]
so that $\tau_{A_1g R_v^{\,2}}= \deg Q_n= \operatorname{ht}(\Ga_1
g^{-1} \ell_*)$, as wanted. For use in the following remark, note that
if $m=\big\lceil \frac{\deg Q_n}{2}\big\rceil$, then by Proposition
\ref{prop:propriinvar} \eqref{item2:propriinvar}, we have
\begin{equation}\label{eq:pourremsuiv}
\sys_1(\aaa^{-m}g R_v^{\,2})=-v(\pi_v^{-m})=m\quad\text{and} \quad
\sys_2(\aaa^{-m}g R_v^{\,2})=\deg Q_n-v(\pi_v^{-m})=\Big\lfloor\frac{\deg
  Q_n}{2}\Big\rfloor\;.
\end{equation}

The middle equality of Equation \eqref{eq:quadrupleq} follows by
induction from Equation \eqref{eq:recurconvergent}, noting that $\deg
a_i\geq 1$ if $i\geq 1$.  The last equality follows from \cite[\S 6.3,
  Remarque 2.]{Paulin02}.
\cqfd

\medskip
\rem The above proof also gives that the compact core of a divergent
$A_1$-orbit $A_1g R_v^{\,2}$ corresponds to the part of the geodesic flow
orbit of a birational geodesic line $\Ga_1\ell$ in the space
$\Ga_1\bs\G_{\rm even}\TT_v$ where the base point $\Ga_1\ell(0)$ lies
exactly in the zigzag part: more precisely
\[
\Xi(A_1(\OOO_v)C_{A_1g\Ga_1})=
\big\{\Ga_1\phi^{2k}\big(\,\wt{g^{-1}\ell_*}\,\big):
0\leq k\leq \operatorname{zz}(\Ga_1g^{-1}\ell_*)\big\}\,.
\]
Note that a quasicenter $\wh x$ of a divergent $A_1$-orbit $x$ is well
defined up to the action of $A_1(\OOO_v)$, and $\Xi$ is a
homeomorphism from $A_1(\OOO_v)\bs \X_1$ to $\Ga_1\bs\G_{\rm even}
\TT_v$, hence looking at the image $\Xi(A_1(\OOO_v)\wh x)$ of the set
of quasicenters is well defined. Equation \eqref{eq:pourremsuiv} in
the above proof gives besides that the quasicenter of a divergent
$A_1$-orbit $A_1g R_v^{\,2}$ (defined by Equation
\eqref{eq:defiquasicenter}) corresponds to the geodesic flow orbit
point of a birational geodesic line in the space $\Ga_1\bs\G_{\rm
  even}\TT_v$ where the base point is almost at the midpoint of the
zigzag part: more precisely
%%
%\todo{vérifier $\lceil\;\rceil$ ou $\lfloor\;\rfloor$...}
%%
\[
\Xi(A_1(\OOO_v)\wh{A_1g R_v^{\,2}})=
\Ga_1\phi^{2m}\big(\,\wt{g^{-1}\ell_*}\,\big)\quad\text{where}\quad
m=\Big\lceil\frac{\operatorname{zz}(\Ga_1g^{-1}\ell_*)}{2}\Big\rceil\,.
\]

\subsection{Type and discriminant of the divergent diagonal orbits}
\label{subsec:typdiscdivdiaorb}

In this section, we introduce two new invariants of the divergent
diagonal orbits in $\X_1$, we gather the technical notation that will
be used in the following sections, and we give a precise description
of the divergent orbits whose equidistribution we will study.

Let $x\in P\!\X$ be the homothety class of an $R_v$-lattice whose
$P\!A$-orbit is divergent in $P\!\X$. By Theorem \ref{theo:defidiv}
\eqref{item4:defidiv}, we know that the orbit $P\!A\,x$ contains at least
one homothety class $[L]$ of an integral $R_v$-lattice $L$. Since the
normalized covolume $\frac{\covol(L)} {\covol(R_v^{\;n})}\in
q_v^{\,\ZZ}$ of $L$ is at least $1$ as $L\subset R_v^{\;n}$, and since
any nonempty subset of $\NN$ has a lower bound, there exists at least
one integral $R_v$-lattice $L_x$ whose homothety class belongs to $P\!A\,x$
and whose covolume is minimal. We define

$\bullet$~ the {\it discriminant}
%%
%\todo{\tiny sans $\log_{q_v}$ ?}
%%
of the divergent $P\!A$-orbit $P\!A\,x$ as $\operatorname{disc}
(P\!A\,x)= \log_{q_v}\frac{\covol(L_x)}{\covol(R_v^{\;n})}\in \NN$

\noindent and

$\bullet$~ the {\it type} of the divergent $P\!A$-orbit $P\!A\,x$ as
the set of types (see Subsection \ref{subsec:lattices}) of the
finitely many integral lattices $L_x$ minimizing the covolume among
the integral $R_v$-lattices whose homothety class belongs to
$P\!A\,x$.
%%
%\todo{\tiny do all minimizing integral lattices have the same type ?}
%%

We refer to the proof of Proposition \ref{prop:descripxt}
\eqref{item3:descripxt} for examples of divergent $P\!A$-orbits having
nonunique homothety classes of covolume-minimizing integral
$R_v$-lattices.

\medskip
We endow the infinite sets $(R_v\ssm \{0\})^{n-1}$ and $R_v\ssm \{0\}$
with the Fréchet filter of the complementary sets of their finite
subsets, and we will denote as usual by ${\displaystyle
  \lim_{+\infty}}$ the limits along this filter.

Let us introduce the notation that will be used in the remainder of
this paper. In this section, we fix an element $\mb{s}= (s_2, \ldots,
s_n) \in (R_v\ssm \{0\})^{n-1}$ such that 
\begin{equation}\label{eq:proprisbold}
  \exists\;\sigma\in \operatorname{Bij}(\llbracket 2,n \rrbracket),
  \qquad s_{\sigma^{-1}(2)}\mid s_{\sigma^{-1}(3)}\mid\ldots\mid
  s_{\sigma^{-1}(n)} \quad\text{ and } s_{\sigma^{-1}(n)} \in n\,\ZZ\,.
\end{equation}
We then define $\mb{s}_*= s_{\sigma^{-1}(n)}\in n\,\ZZ$. Since the
valuation $v$ is nonpositive on $R_v\ssm\{0\}$ by Equation
\eqref{eq:inversiRv}, we have $v(\mb{s}_*)\leq 0$. Note that
$\mb{s}\ra+\infty$ if and only if $\mb{s}_*\ra+\infty$, hence if and
only if $-v(s)\ra+\infty$.  We define (independently of the choice of
such a permutation $\sigma$)
%%
%\todo{\tiny voir si inégalités strictes dans $\wt\loz_{\mb{s}}$ ?}
%%
%%
%\todo{\tiny dans $\Lambda_{\mb{s}}$, remplacer $=R_v$ par $=I_i$ où
%  $I_i\in\I_v^+$, peut-être avec $I_i$ constant ?}
%%
\begin{align*}
\Lambda_{\mb{s}}&=\Big\{\big(\frac{r_2}{s_2},\frac{r_3}{s_3},
  \ldots, \frac{r_n}{s_n}\big):
\forall\;i\in\intbra{2,n},\;r_i\in R_v\;\text{and}\;
r_iR_v+s_iR_v=R_v\Big\} \mod R_v^{\;n-1}\,,
\\
\Delta_{\mb{s}}&=\big\{\mb{k}=(k_1,\ldots, k_n)\in\ZZ_0^n:
k_1\geq v(\mb{s}_*)\;\text{and}\;
\forall\;i\in\intbra{2,n},\;k_i\geq 0\big\}\,,
\\
\loz_{\mb{s}}&=\big\{\mb{k}=(k_1,\ldots, k_n)\in\ZZ_0^n:
\forall\;i\in\intbra{2,n},\;0\leq k_i\leq k_1-v(\mb{s}_*)\big\}\,,
\\
\mb{k}_{\mb{s}}&=\Big(\frac{n-1}{n}\,v(\mb{s}_*),
-\frac{1}{n}\,v(\mb{s}_*),\ldots,
-\frac{1}{n}\,v(\mb{s}_*)\Big)\in\ZZ_0^n\,,
\\
\wt\Delta_{\mb{s}}&=\Delta_{\mb{s}}-\mb{k}_{\mb{s}}=\Big\{\mb{k}
=(k_1,\ldots, k_n)\in\ZZ_0^n:
\forall\;i\in\intbra{1,n},\;k_i\geq \frac{v(\mb{s}_*)}{n}\Big\}\,,
\\
\wt\loz_{\mb{s}}&=\loz_{\mb{s}}-\mb{k}_{\mb{s}}=\big\{\mb{k}
=(k_1,\ldots, k_n)\in\ZZ_0^n:
\forall\;i\in\intbra{2,n},\;\frac{v(\mb{s}_*)}{n}\leq k_i\leq k_1\big\}\,.
\end{align*}
Let us make some comments on this notation. We have $\mb{k}_{\mb{s}}
\in \loz_{\mb{s}}\subset \Delta_{\mb{s}}$ (see the following picture
when $n=3$). We have
\begin{equation}\label{eq:cardLambdasss}
\card(\Lambda_{\mb{s}})=\prod_{i=2}^n\;\varphi_v(s_i)\;,
\quad\text{hence}\quad
  \card(\Lambda_{(s,\ldots,s)})=(\varphi_v(s))^{n-1}\;,
\end{equation}
where $\varphi_v$ is the Euler function of $R_v$ defined in Equation
\eqref{eq:defiEulerfunct}.

\begin{center}
  \!\!\!\!\!\!\!\!\!\!\!\!\!\!\!\!\!\!\!\!\!\!\!\!\!\!\!
  \begin{picture}(0,0)%
\includegraphics{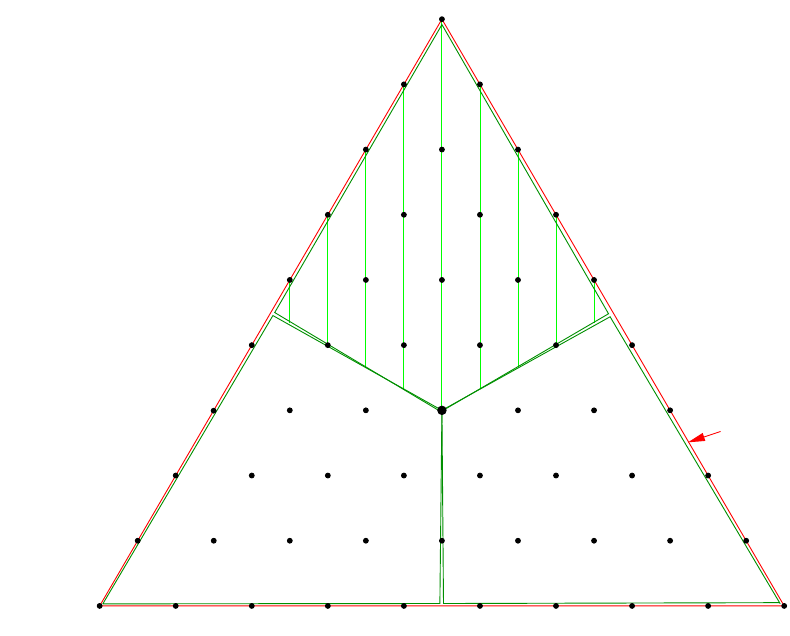}%
\end{picture}%
\setlength{\unitlength}{3812sp}%
\begingroup\makeatletter\ifx\SetFigFont\undefined%
\gdef\SetFigFont#1#2#3#4#5{%
  \reset@font\fontsize{#1}{#2pt}%
  \fontfamily{#3}\fontseries{#4}\fontshape{#5}%
  \selectfont}%
\fi\endgroup%
\begin{picture}(6523,5251)(-1184,-2555)
\put(2211,119){\makebox(0,0)[lb]{\smash{{\SetFigFont{11}{13.2}{\rmdefault}{\mddefault}{\updefault}{\color[rgb]{0,.56,0}{\Large$\loz_{\mb{s}}$}}%
}}}}
\put(2521,-871){\makebox(0,0)[lb]{\smash{{\SetFigFont{11}{13.2}{\rmdefault}{\mddefault}{\updefault}{\color[rgb]{0,0,0}$k_{\mb{s}}$}%
}}}}
\put(4861,-871){\makebox(0,0)[lb]{\smash{{\SetFigFont{11}{13.2}{\rmdefault}{\mddefault}{\updefault}{\color[rgb]{1,0,0}$\Delta_{\mb{s}}$}%
}}}}
\put(3781,389){\makebox(0,0)[lb]{\smash{{\SetFigFont{11}{13.2}{\rmdefault}{\mddefault}{\updefault}{\color[rgb]{0,0,0}$(\lfloor\frac{v(s_3)}{2}\rfloor,-\lfloor\frac{v(s_3)}{2}\rfloor,0)$}%
}}}}
\put(-404,389){\makebox(0,0)[lb]{\smash{{\SetFigFont{11}{13.2}{\rmdefault}{\mddefault}{\updefault}{\color[rgb]{0,0,0}$(\lfloor\frac{v(s_3)}{2}\rfloor,0,-\lfloor\frac{v(s_3)}{2}\rfloor)$}%
}}}}
\put(4816,-2491){\makebox(0,0)[lb]{\smash{{\SetFigFont{11}{13.2}{\rmdefault}{\mddefault}{\updefault}{\color[rgb]{0,0,0}$(v(s_3),-v(s_3),0)$}%
}}}}
\put(-1169,-2491){\makebox(0,0)[lb]{\smash{{\SetFigFont{11}{13.2}{\rmdefault}{\mddefault}{\updefault}{\color[rgb]{0,0,0}$(v(s_3),0,-v(s_3))$}%
}}}}
\put(3151,1469){\makebox(0,0)[lb]{\smash{{\SetFigFont{11}{13.2}{\rmdefault}{\mddefault}{\updefault}{\color[rgb]{0,0,0}$(-2,2,0)$}%
}}}}
\put(2836,2009){\makebox(0,0)[lb]{\smash{{\SetFigFont{11}{13.2}{\rmdefault}{\mddefault}{\updefault}{\color[rgb]{0,0,0}$(-1,1,0)$}%
}}}}
\put(2521,2549){\makebox(0,0)[lb]{\smash{{\SetFigFont{11}{13.2}{\rmdefault}{\mddefault}{\updefault}{\color[rgb]{0,0,0}$(0,0,0)$}%
}}}}
\put(1036,1469){\makebox(0,0)[lb]{\smash{{\SetFigFont{11}{13.2}{\rmdefault}{\mddefault}{\updefault}{\color[rgb]{0,0,0}$(-2,0,2)$}%
}}}}
\put(1351,2009){\makebox(0,0)[lb]{\smash{{\SetFigFont{11}{13.2}{\rmdefault}{\mddefault}{\updefault}{\color[rgb]{0,0,0}$(-1,0,1)$}%
}}}}
\put(2251,1289){\makebox(0,0)[lb]{\smash{{\SetFigFont{11}{13.2}{\rmdefault}{\mddefault}{\updefault}{\color[rgb]{0,0,0}$(-2,1,1)$}%
}}}}
\end{picture}%

\end{center}

\noindent We denote by $\delta_y$ the unit Dirac mass at $y$. With
$\uuu_{\mb{t}}$ for $\mb{t}\in K_v^{\;n-1}$ defined in Equation
\eqref{eq:defuuut}, we define two probability measures on $\X$ by
%%
%\todo{\tiny vérifier le signe de $\exp(\pm\mb{k})$ vs celui de
%  [DS] bas page 1222}
%%
\begin{equation}\label{eq:definuDelta}
\nu_{\mb{s}}^\Delta=
\frac{1}{\card\;\Lambda_{\mb{s}}\;\card\;\Delta_{\mb{s}}}\;
\sum_{\mb{t}\in\Lambda_{\mb{s}},\;\mb{k}\in\Delta_{\mb{s}}}\int_{a\in A(\OOO_v)}
\;\delta_{a\exp(\mb{k})\,\uuu_{\mb{t}}R_v^{\;n}}\;da\;,
\end{equation}
and similarly, replacing $\Delta_{\mb{s}}$ by $\loz_{\mb{s}}$,
\begin{equation}\label{eq:definuloz}
\nu_{\mb{s}}^\loz=
\frac{1}{\card\;\Lambda_{\mb{s}}\;\card\;\loz_{\mb{s}}}\;
\sum_{\mb{t}\in\Lambda_{\mb{s}},\;\mb{k}\in\loz_{\mb{s}}}\int_{a\in A(\OOO_v)}
\;\delta_{a\exp(\mb{k})\,\uuu_{\mb{t}}R_v^{\;n}}\;da\,.
\end{equation}
We define similarly two probability measures $\ov\nu_{\mb{s}}^\Delta$
and $\ov\nu_{\mb{s}}^\loz$ on $\X_1$ by replacing $A$ by $A_1$ in
Equations \eqref{eq:definuDelta} and \eqref{eq:definuloz}, so that
by Equation \eqref{eq:decompdmA1} we have
\begin{equation}\label{eq:relatnudeltnuloz}
  \ov\nu_{\mb{s}}^\Delta=\frac{q_v}{q_v-1}\int_{\lambda\in\OOO_v^{\,\times}}
  \xi(\lambda)_*\nu_{\mb{s}}^\Delta\;d\vol_v(\lambda)
%\quad\text{and}\quad
%\ov\nu_{\mb{s}}^\loz=\frac{q_v}{q_v-1}
%\int_{\lambda\in\OOO_v^{\,\times}}
%\xi(\lambda)_*\nu_{\mb{s}}^\loz\;d\vol_v(\lambda)
\,.
\end{equation}

%\noindent For every $s\in R_v$ with $v(s)\in n\ZZ$ (so that
%$s\neq 0$), setting $\mb{s}_s=(s,s,\ldots,s)\in (R_v\ssm
%\{0\})^{\;n-1}$, we define
%\[
%\Lambda_{s}=\Lambda_{\mb{s}_s}, \;\Delta_{s}=\Delta_{\mb{s}_s}, \;
%\loz_{s}=\loz_{\mb{s}_s}, \;\mb{k}_{s}=\mb{k}_{\mb{s}_s}, \;
%\wt\Delta_{s}=\wt\Delta_{\mb{s}_s}, \;\wt\loz_{s}=\wt\loz_{\mb{s}_s},\;
%\nu_{s}^\Delta=\nu_{\mb{s}_s}^\Delta\;\text{and}\;
%\nu_{s}^\loz=\nu_{\mb{s}_s}^\loz\,.
%\]

\medskip
We denote by $\Scal_n$ the permutation group of $\llbracket1,
n\rrbracket$, and by $\Scal_{n-1}$ the stabilizer of $\{1\}$ in
$\Scal_n$.  The group $\Scal_{n-1}$ acts $K$-linearly on $K^{n-1}$
(preserving the subset $(R_v\ssm\{0\})^{n-1}$) by $\sigma\cdot \mb{t}=
(t_{\sigma^{-1}(2)}, \ldots, t_{\sigma^{-1}(n)})$ for all $\sigma\in
\Scal_{n-1}$ and $\mb{t}=(t_2,\ldots,t_n) \in K^{n-1}$. By
construction, for all $\sigma\in\Scal_{n-1}$ and $\mb{s}\in
(R_v\ssm\{0\})^{n-1}$, the element $\mb{s}$ satisfies Equation
\eqref{eq:proprisbold} if and only if $\sigma\cdot \mb{s}$ does, and
we have
\begin{equation}\label{eq:equivarsigma}
  \Lambda_{\sigma\cdot\mb{s}}=\sigma\cdot\Lambda_{\mb{s}},\quad
  v((\sigma\cdot \mb{s})_*)=v( \mb{s}_*),\quad
\Delta_{\sigma\cdot\mb{s}}=\Delta_{\mb{s}},\quad
\loz_{\sigma\cdot\mb{s}}=\loz_{\mb{s}},\quad\text{and}\quad 
\mb{k}_{\sigma\cdot\mb{s}}=\mb{k}_{\mb{s}}\,.
\end{equation}

The group $\Scal_n$ acts $\ZZ$-linearly on $\ZZ_0^n$ by $\sigma\cdot
\mb{k}=(k_{\sigma^{-1}(1)},\ldots,k_{\sigma^{-1}(n)})$ for all $\sigma
\in\Scal_{n}$ and $\mb{k}=(k_1,\ldots,k_n)\in\ZZ_0^n$. Note that the
subset $\wt\Delta_{\mb{s}}$ of $\ZZ_0^n$ is invariant under the action
of $\Scal_n$. As said above, the subsets $\Delta_{\mb{s}}$ and
$\loz_{\mb{s}}$ of $\ZZ_0^n$, as well as their point
$\mb{k}_{\mb{s}}$, are invariant under the action of $\Scal_{n-1}=
\operatorname{Stab}_{\Scal_n}\{1\}$. Let $\sigma_n=(1\ldots n)$ be the
standard $n$-cycle in $\Scal_n$. The cyclic group $\sigma_n^{\;\ZZ}$
of order $n$ generated by $\sigma_n$ acts freely on
$\wt\Delta_{\mb{s}} \ssm\{\mb{0}\}$ where $\mb{0}=(0,\ldots,0)$. The
subset $\wt\loz_{\mb{s}}$ is a (weak) fundamental domain for the
action of $\sigma_n^{\;\ZZ}$ on $\wt\Delta_{\mb{s}}$~: we have (with
nondisjoint union)
\[
\wt\Delta_{\mb{s}}=
\bigcup_{j=0}^{n-1}\;\sigma_n^{\; j}\cdot\wt\loz_{\mb{s}}\,.
\]

For every $\sigma\in \Scal_n$, we denote by $P_\sigma\in\GL_n(K_v)=
\GL(K_v^{\;n})$ the corresponding permutation matrix of the canonical
basis $(e_1,\ldots, e_n)$ of $K_v^{\;n}$, so that $P_\sigma(e_i)=
e_{\sigma(i)}$ for every $i\in \intbra{1,n}$. The map $\sigma\mapsto
P_\sigma$ is a group morphism from $\Scal_n$ to $\GL_n(K_v)$.  This
\mbox{$K_v$-linear} action of $\Scal_n$ on $K_v^{\;n}$ gives an action of
$\Scal_n$ on the set of $R_v$-lattices $x$ of $K_v^{\;n}$ by $x\mapsto
\sigma(x)= \{P_\sigma(w):w\in x\}$, that preserves $\X_1$ since the
determinant of $P_\sigma$ is the signature $\varepsilon(\sigma)
\in\{\pm1\}$ of $\sigma$ for every $\sigma\in\Scal_n$.
%The group $\Scal_n$  also acts on $\M_n(K_v)$ by conjugation
%$(\sigma,M)\mapsto M^\sigma= P_\sigma MP_\sigma^{-1}$.

\blemm \label{lem:equivarnusbold} For every $\mb{s}\in
(R_v\ssm\{0\})^{n-1}$ satisfying Equation \eqref{eq:proprisbold} and
for every $\sigma\in \Scal_{n-1}$, we have $\nu^\Delta_{\sigma\cdot\mb{s}}
=\sigma_*(\nu^\Delta_{\mb{s}})$, $\nu^\loz_{\sigma\cdot\mb{s}}
=\sigma_*(\nu^\loz_{\mb{s}})$, $\ov\nu^\Delta_{\sigma\cdot\mb{s}}
=\sigma_*(\,\ov\nu^\Delta_{\mb{s}})$ and
$\ov\nu^\loz_{\sigma\cdot\mb{s}} =\sigma_*(\,\ov\nu^\loz_{\mb{s}})$.
\elemm

\dem We only prove the first equality. The proofs of the other ones
are similar.  Let $\sigma,\mb{s}$ be as in the statement. Let $a\in
A(\OOO_v)$, $\mb{k} \in\ZZ_0^n$ and $\mb{t}\in K^{n-1}$. Since
$\uuu_{\sigma\cdot \mb{t}}=P_\sigma \uuu_{\mb{t}} P_\sigma^{-1}$ and
$R_v^{\;n}$ is invariant by $P_\sigma^{-1}$, we have
\begin{align*}
\sigma(a\,\exp(\mb{k})\;\uuu_{\mb{t}}\;R_v^{\;n})&=
P_\sigma\,a\,\,\exp(\mb{k})\;\uuu_{\mb{t}}\;R_v^{\;n}=
P_\sigma\,a\,P_\sigma^{-1}P_\sigma\,\exp(\mb{k})\,P_\sigma^{-1}P_\sigma
\uuu_{\mb{t}}\,P_\sigma^{-1}\,R_v^{\;n}\\&=
(P_\sigma\,a\,P_\sigma^{-1})\exp(\sigma\cdot\mb{k})\;\uuu_{\sigma\cdot\mb{t}}\;
R_v^{\;n}\,.
\end{align*}
For every $\sigma\in\Scal_n$, the conjugation by $P_\sigma$ in
$\GL_n(K_v)$ preserves $A(\OOO_v)$ and its Haar measure.  Hence, by
Equation \eqref{eq:equivarsigma} and by changes of variables in the
sums and integral, we have
\begin{align*}
\nu_{\sigma\cdot\mb{s}}^\Delta&=
\frac{1}{\card\;\Lambda_{\sigma\cdot\mb{s}}\;\card\;\Delta_{\sigma\cdot\mb{s}}}\;
\sum_{\mb{t}\in\Lambda_{\sigma\cdot\mb{s}},\;\mb{k}\in\Delta_{\sigma\cdot\mb{s}}}\int_{a\in
  A(\OOO_v)} \;\delta_{a\exp(\mb{k})\,\uuu_{\mb{t}}R_v^{\;n}}\;da\\&=
\frac{1}{\card\;\Lambda_{\mb{s}}\;\card\;\Delta_{\mb{s}}}\;
\sum_{\mb{t}\in\Lambda_{\mb{s}},\;\mb{k}\in\Delta_{\mb{s}}}\int_{a\in A(\OOO_v)}
\;\delta_{(P_\sigma a P_\sigma^{\,-1})\exp(\sigma\cdot\mb{k})\,\uuu_{\sigma\cdot\mb{t}}
R_v^{\;n}}\;da=\sigma_*\big(\nu_{\mb{s}}^\Delta\big)\,.\;\;\;\Box
\end{align*}

\medskip
Note for future use that for all $\sigma\in \Scal_n$, $a\in
A_1(\OOO_v)$, $\mb{k} \in\ZZ_0^n$ and $\mb{t}\in K^{n-1}$, we have
\begin{equation}\label{eq:sigmapermut}
\sigma\big(a\,\exp(\mb{k})\;\uuu_{\mb{t}}\;R_v^{\;n}\big)=
\exp(\sigma\cdot\mb{k})\;\sigma(a\;\uuu_{\mb{t}}\;R_v^{\;n})\,.
\end{equation}

\medskip
For every $\mb{t}=\big(t_2,\ldots, t_n\big)
\in K_v^{\;n-1}$, let
\begin{equation}\label{eq:defixt}
x_{\mb{t}}=\uuu_{\mb{t}}R_v^{\;n}=R_v\Big(e_1+\sum_{i=2}^n\;t_i\;e_i\Big)
+R_v\,e_2+\ldots+R_v\,e_n\,.
\end{equation}
Note that $x_{\mb{t}}\in\X$ since $\uuu_{\mb{t}}\in G=\SL_n(K_v)$. We
have $x_{\mb{0}}=R_v^{\;n}$, and if $\mb{t}'\in\mb{t}+R_v^{\;n-1}$,
then $x_{\mb{t}'}=x_{\mb{t}}$. Note that $x_{\mb{t}}$ is a unimodular
\mbox{$R_v$-lattice,} which is rational if and only if $\mb{t}\in
K^{n-1}$ and integral if and only if $\mb{t}\in R_v^{\;n-1}$. For
every $i\in\intbra{2,n}$, we have $x_{\mb{t}}\cap (K_v e_i)=R_v\,e_i$.
Hence by the definition of the directional systoles and by Equation
\eqref{eq:inversiRv}, we have
\begin{equation}\label{eq:sysigeq2xt}
  \forall\;i\in\intbra{2,n},\qquad
  \sys_i(x_{\mb{t}})=\log_{q_v}\min_{s\in R_v\ssm\{0\}}|s|=0\,.
\end{equation}
Since
\begin{equation}\label{eq:interxtpremaxe}
  x_{\mb{t}}\cap (K_v e_1)=\{\lambda_1 e_1:\lambda_1\in R_v
  \;\;\text{and}\;\;\forall\;i\in\intbra{2,n},\;\;
  \lambda_1 t_i\in R_v\}\;,
\end{equation}
the $R_v$-lattice $x_{\mb{t}}$ is axial if and only if it is rational,
hence if and only if $\mb{t}\in K^{n-1}$.

\bprop \label{prop:descripxt} Let $\mb{s}\in (R_v\ssm \{0\})^{n-1}$
satisfying Equation \eqref{eq:proprisbold}. Let $\mb{t}=\big(
\frac{r_2}{s_2},\ldots, \frac{r_n}{s_n} \big) \in \Lambda_{\mb{s}}$.
\begin{enumerate}
\item\label{item1:descripxt} The first directional systole of the
  $R_v$-lattice $x_{\mb{t}}$ is $\sys_1(x_{\mb{t}})=
  -v(\mb{s}_*)$. The truncated covolume of the $A_1$-orbit of
  $x_{\mb{t}}$ is $\tau_{x_{\mb{t}}} =-v(\mb{s}_*)$. The coordinate
  sublattice of $x_{\mb{t}}$ is $(x_{\mb{t}})^{\rm coo}= R_v\,\mb{s}_*
  \,e_1+ R_v\,e_2 +\cdots +R_v\,e_n$. The compact core of the
  $A_1$-orbit of $x_{\mb{t}}$ is $C_{x_{\mb{t}}}=A_1(\OOO_v)
  \exp(\Delta_{\mb{s}}) x_{\mb{t}}$. A quasicenter of the $A_1$-orbit
  of $x_{\mb{t}}$ is $\widehat{x_{\mb{t}}}= \exp(\mb{k}_{\mb{s}})
  x_{\mb{t}}$.
\item\label{item1_5:descripxt} As $\mb{s}\ra+\infty$, we have
\begin{equation}\label{eq:initem1_5:descripxt}
\card\; \loz_{\mb{s}}=\frac{1}{n!}(-v(\mb{s}_*))^{n-1}
+\bigO((-v(\mb{s}_*))^{n-2}) =\frac{1}{n}\;\card\; \Delta_{\mb{s}}
+\bigO((-v(\mb{s}_*))^{n-2})\,.
\end{equation}
\item\label{item2:descripxt} We have $x_{\mb{t}}\cap R_v^{\;n}=
  R_v\,\mb{s}_*\,e_1+R_v\,e_2 +\cdots +R_v\,e_n=(x_{\mb{t}})^{\rm coo}$.
  Hence the type of the integral lattice $x_{\mb{t}}\cap R_v^{\;n}$ is
  $(1,\ldots, 1,\mb{s}_*)$.
\item\label{item3:descripxt} The $P\!A$-orbit of the homothety class
  of $x_{\mb{t}}$ has discriminant $\prod_{i=2}^n|s_i|$ and has type
  $\{(1,s_{\sigma^{-1}(2)},\ldots, s_{\sigma^{-1}(n)})\}$ where $\sigma
  \in\Scal_{n-1}$ is such that $s_{\sigma^{-1}(2)}\mid\ldots \mid
  s_{\sigma^{-1}(n)}$.
\item\label{item4:descripxt} With $\nu_y$ for $y\in\X$ the probability
  measure given by Equation \eqref{eq:calcpratnuy}, we have
  \[
  \nu_{\mb{s}}^\Delta=\frac{1}{\card\;\Lambda_{\mb{s}}}\;
  \sum_{\mb{t}\in\Lambda_{\mb{s}}}\;\nu_{x_{\mb{t}}}\,.
  %\quad\text{and} \quad
  %\ov\nu_{\mb{s}}^\Delta=\frac{1}{\card\;\Lambda_{\mb{s}}}\;
  %\sum_{\mb{t}\in\Lambda_{\mb{s}}}\;\ov\nu_{x_{\mb{t}}}
  \]
  As $s\ra+\infty$ in $R_v\ssm\{0\}$ with $v(s) \in n\ZZ$, defining
  ${\mathfrak s}=(s,\ldots,s)\in (R_v\ssm\{0\})^{\,n-1}$, for every
  $f\in {\rm C}_c(\X)$, we have
%%
%\todo{\tiny analogue lem 2.5 [DS]}
%%
\[
\nu_{{\mathfrak s}}^\Delta(f)=\frac{1}{n}\;\sum_{j=0}^{n-1}\;(\sigma_n^{\; j})_*
\nu_{{\mathfrak s}}^\loz(f)+\bigO\Big(\frac{\|f\|_\infty}{-v(s)}\Big)
%\]
%\[\text{and} \quad
%\ov\nu_{s}^\Delta(f)=\frac{1}{n}\;\sum_{j=0}^{n-1}\;(\sigma_n^{\; j})_*
%\ov \nu_{s}^\loz(f)+\bigO\Big(\frac{\|f\|_\infty}{-v(s)}\Big)
\,.
\]
\end{enumerate}
\eprop

\medskip
\dem If $r_2,\ldots, r_n, r'_2,\ldots,
r'_n\in R_v$ satisfy $r'_2\equiv r_2\mod s_2,\ldots, r'_n\equiv
r_n\mod s_n$, and if $\mb{t}=\big(\frac{r_2}{s_2},\ldots,
\frac{r_n}{s_n} \big)$, $\mb{t}'=\big(\frac{r'_2}{s_2},\ldots,
\frac{r'_n}{s_n} \big)$, then we have $\mb{t}'-\mb{t}\in R_v^{\;n-1}$
and $x_{\mb{t}}= x_{\mb{t}'}$.  In particular, the $R_v$-lattice
$x_{\mb{t}}$, as well as the measures $\nu_{x_{\mb{t}}}$,
$\nu_{\mb{s}}^\Delta$ and $\nu_{\mb{s}}^\loz$, do not depend on the
choice of representatives of the elements $\mb{t}$ in the index set
$\Lambda_{\mb{s}}\subset K^{n-1}/R_v^{n-1}$.

%For every $\sigma\in\Scal_{n-1}$, we have
%$x_{\sigma\cdot\mb{t}}=P_\sigma x_{\mb{t}}$ and $P_\sigma$ fixes
%$K_ve_1$, so that $\tau_{x_{\sigma\cdot\mb{t}}}=\tau_{x_{\mb{t}}}$ by
%Equation \eqref{eq:defitruncovol} and \eqref{eq:deficoosublat} By
%Equation \eqref{eq:equivarsigma} and Lemma \ref{lem:equivarnusbold},
%since for every $\sigma\in\Scal_{n-1}$ , we may assume that $s_2\mid
%s_3\mid\dots\mid s_n$, so that $\mb{s}_*=s_n$.

\medskip
\eqref{item1:descripxt} By Equation \eqref{eq:interxtpremaxe}, the set
$x_{\mb{t}}\cap (K_v e_1)$ consists in the elements $\lambda_1 e_1$
where $\lambda_1\in R_v$ is such that, for every $i\in\intbra{2,n}$,
we have $\lambda_1 \frac{r_i}{s_i}\in R_v$. Since $r_i$ is invertible
modulo $s_i$ and by Equation \eqref{eq:proprisbold}, this occurs if
and only if $\lambda_1\in \bigcap_{i=2}^n s_iR_v= \mb{s}_*\,R_v$.
Hence we have $x_{\mb{t}}\cap (K_v e_1)=\mb{s}_*\,R_v$. By the
definition of the first directional systole, this proves the first
claim. By the definition \eqref{eq:defitruncovol} of the truncated
covolume and by Equation \eqref{eq:sysigeq2xt}, we hence have
$\tau_{x_{\mb{t}}} =\sum_{i=1}^n\sys_i(x_{\mb{t}})=-v(\mb{s}_*)$.  By
the definition \eqref{eq:deficoosublat} of the coordinate sublattice,
we have
\[
(x_{\mb{t}})^{\rm coo}= (x_{\mb{t}}\cap K_v e_1)+\ldots +
(x_{\mb{t}}\cap K_v e_n)= R_v\,\mb{s}_*\,e_1+R_v\,e_2 +\cdots
+R_v\,e_n\,.
\]
By Equation \eqref{eq:defiDeltax}, by the above computation of the
directional systoles of $x_{\mb{t}}$ and by the definition of
$\Delta_{\mb{s}}$, we have $\Delta^{x_{\mb{t}}}=\Delta_{\mb{s}}$. This
gives the description of the compact core of the $A_1$-orbit of
$x_{\mb{t}}$ by Equation \eqref{eq:deficompactcore}.

Recall that $-v(\mb{s}_*)\in n\NN$ and that $\mb{k}_{\mb{s}}=\big(k_1=
\frac{n-1}{n} v(\mb{s}_*),k_2 =-\frac{v(\mb{s}_*)}{n}, \ldots,k_n=-
\frac{v(\mb{s}_*)}{n}\big)$. By the first claim of Proposition
\ref{prop:propriinvar} \eqref{item2:propriinvar}, we have that
\[\sys_i(\exp(\mb{k}_{\mb{s}})x_{\mb{t}})=\sys_i(x_{\mb{t}})-
v(\pi_v^{-k_i})=\sys_i(x_{\mb{t}})+k_i
\]
is equal to $-\frac{v(\mb{s}_*)}{n}$ if $i\in\llbracket 2,n\rrbracket$
by Equation \eqref{eq:sysigeq2xt} and to $-v(\mb{s}_*)+\frac{n-1}{n}
v(\mb{s}_*) = -\frac{v(\mb{s}_*)}{n}$ if $i=1$. The description of the
quasicenters of the $A_1$-orbit of $x_{\mb{t}}$ (well-defined up to
left translation by an element of $A_1(\OOO_v)$) then follows from
their definition that requires, by Equation \eqref{eq:defiquasicenter}
and since $\tau_{x_{\mb{t}}}=-v(\mb{s}_*)\in\ZZ$, that we have
$\sys_i(\wh{x_{\mb{t}}})= \frac{\tau_{x_{\mb{t}}}}{n}
=-\frac{v(\mb{s}_*)}{n}$ for every $i\in\llbracket 1,n\rrbracket$.

\medskip
\eqref{item1_5:descripxt} Using the facts that $\Delta^{x_{\mb{t}}}=
\Delta_{\mb{s}}$ and $\tau_{x_{\mb{t}}} =-v(\mb{s}_*)$ seen in Assertion
\eqref{item1:descripxt}, the first equality in Equation
\eqref{eq:initem1_5:descripxt} follows by the same proof as the one of
Proposition \ref{prop:propriinvar} \eqref{item5:propriinvar}, by
comparing with the Euclidean volume of the polytop
\[
\wh{\loz}'(\alpha)
=\big\{\mb{t}=(t_1,\dots,t_n)\in \RR_0^n:\forall \;i\in\intbra{2,n},\;
0\leq t_i\leq t_1+\alpha\big\}
\]
for $\alpha=-v(\mb{s}_*)$ whose images under the powers
$\sigma_n^{\;0}, \sigma_n^{\;1},\dots,\sigma_n^{\;n-1}$ of $\sigma_n$
have pairwise disjoint interior and cover $\Delta'(\alpha)$. The
second equality in Equation \eqref{eq:initem1_5:descripxt} then
follows from Proposition \ref{prop:propriinvar}
\eqref{item5:propriinvar}.

\medskip
\eqref{item2:descripxt} Let $\lambda_1,\ldots,\lambda_n\in R_v$. We
have $\lambda_1\big(e_1+\sum_{i=2}^n\;\frac{r_i}{s_i}\;e_i\big)
+\sum_{i=2}^n \lambda_i\,e_i\in R_v^{\;n}$ if and only if for every
$i\in\intbra{2,n}$ we have $\lambda_1\frac{r_i}{s_i}\in
R_v$, hence if and only if $\lambda_1\in \mb{s}_*\,R_v$ as seen in
Assertion \eqref{item1:descripxt}.  The result follows by Equation
\eqref{eq:defixt}.

\medskip
\eqref{item3:descripxt} Let $a\in\wt A$ be a diagonal matrix with
diagonal coefficients $a_1,\ldots, a_n\in K_v^{\,\times}$. The
$R_v$-lattice $a\, x_{\mb{t}} =R_v\big(a_1\,e_1+\sum_{i=2}^n\;
\frac{r_i}{s_i}\;a_i\,e_i\big) +R_v\,a_2\,e_2+\ldots+R_v\,a_n\,e_n$ is
integral, that is contained in $R_v^{\;n}$, if and only if $a_i\in
R_v$ for every $i\in\intbra{2,n}$, $a_1\in R_v$ and
$\frac{r_i}{s_i}\;a_i\in R_v$ for every $i\in\intbra{2,n}$, hence if
and only if $a_1\in R_v$ and $a_i\in s_i\,R_v$ for every
$i\in\intbra{2,n}$ since $r_i$ is invertible modulo $s_i$,

Recall that by Equation \eqref{eq:inversiRv}, for every $z\in
R_v\ssm\{0\}$, we have $|z|\geq 1$, with equality if and only if $z\in
R_v^{\,\times}$. Recall that $x_{\mb{t}}$ is unimodular. Hence by
Equation \eqref{eq:changevarcovol}, we have
\[
\frac{\covol(a\, x_{\mb{t}})}{\covol(R_v^{\;n})}=|\det a|\;
\frac{\covol(x_{\mb{t}})}{\covol(R_v^{\;n})}=\prod_{i=1}^n|a_i|
\geq \prod_{i=2}^n|s_i|
\]
with equality if and only if $a_1\in R_v^{\,\times}$ and $a_i\in s_i\,
R_v^{\,\times}$ for every $i\in\intbra{2,n}$. Therefore the integral
$R_v$-lattices contained in elements of the orbit $\wt A \,x_{\mb{t}}$
with minimal covolume are exactly the $R_v$-lattices $L=R_v\big(e_1+
\sum_{i=2}^n\; r_i\;a'_i\,e_i\big) +R_v\,s_2\,e_2 +\ldots+
R_v\,s_n\,e_n$ where $a'_2, \ldots,a'_n\in R_v^{\,\times}$. Note that
there is no uniqueness of such an $R_v$-lattice $L$ in general.

Since $R_v^{\;n}=R_v\big(e_1+\sum_{i=2}^n\; r_i\;a'_i\,e_i\big) +
R_v\,e_2 +\ldots+R_v\,e_n$ by an immediate change of $R_v$-basis, the
$R_v$-module $R_v^{\;n}/L$ is isomorphic to $\prod_{i=2}^n R_v/
(s_iR_v)$ for every such $L$. All such integral lattices $L$ hence
have the same type, and this proves the result by Equation
\eqref{eq:proprisbold}.

\medskip
\eqref{item4:descripxt}
%We prove the two claims in $\X$, the corresponding claims in
%$\X_1$ follow by averaging, using Equation \eqref{eq:relatnudeltnuloz}.
%
Since $\Delta^{x_{\mb{t}}}=\Delta_{\mb{s}}$ for every $\mb{t}\in
\Lambda_{\mb{s}}$ by Assertion \eqref{item1:descripxt}, the first
claim of Assertion \eqref{item4:descripxt} follows from the
definitions of the measure $\nu_{\mb{s}}^\Delta$ in Equation
\eqref{eq:definuDelta} and of the measure $\nu_y$ for any
$R_v$-lattice $y\in\X$ in Equation \eqref{eq:calcpratnuy}.

\medskip 
Let us prove the second claim of Assertion \eqref{item4:descripxt},
which uses the symmetry of the $(n-1)$-uples $(s,\ldots, s)$ for $s\in
R_v$. We start with a computational lemma.

We fix $s\in R_v$ with $v(s)\in n\ZZ$ and we consider the $(n-1)$-uple
${\mathfrak s}= (s,\ldots, s)$. Let $a\in A(\OOO_v)$ with diagonal
coefficients $a_1,\ldots, a_n$ and let $\mb{t}= \big(\frac{r_2}{s},
\ldots, \frac{r_n}{s} \big)\!\!\mod R_v^{\;n-1}\in \Lambda_{\mathfrak s}$.
Let $s'\in \OOO_v^{\,\times}$ be such that $s=\pi_v^{\;v(s)} \,s'$,
and let $\overline{r_n} \in R_v$ be an inverse of $r_n$ modulo
$s$. Let $a'$ be the element of $A(\OOO_v)$ with diagonal coefficients
\[
a'_1=(s')^{-1} a_n, a'_2=s'a_1,a'_3=a_2,\ldots, a'_n=a_{n-1}\,.
\]
Note that the map $a\mapsto a'$ is a homeomorphism of $A(\OOO_v)$
preserving its Haar measure. Let $\mb{t}'=\big(\frac{\overline{r_n}}
{s},\frac{\overline{r_n}\,r_2}{s}, \ldots, \frac{\overline{r_n}
  \,r_{n-1}}{s} \big)\!\!\mod R_v^{\;n-1}$, and note that the map
$\mb{t} \mapsto \mb{t}'$ is a bijection of $\Lambda_{\mathfrak s}$
with inverse $\big(\frac{r'_2} {s}, \ldots, \frac{r'_n}{s} \big)
\!\!\mod R_v^{\;n-1} \mapsto \big(\frac{\overline{r'_2}\,r'_3}{s}, \ldots,
\frac{\overline{r'_2}\,r'_n}{s}, \frac{\overline{r'_2}}{s} \big)\!\!\mod
R_v^{\;n-1}$.

\blemm \label{lem:tournemanege} Denoting by $y_{a,\mb{t}}$ the
$R_v$-lattice $a\,\exp(\mb{k}_{\mathfrak s})\,x_{\mb{t}}\in \X$, we
have $\sigma_n(y_{a,\mb{t}})=y_{a',\mb{t}'}$.
\elemm

\dem Let $w=s^{-1}\big(a'_2\; e_1+\sum_{j=2}^n\;r_j\;a_j\;e_j\big)\in
K_v^{\;n}$. Since $\pi_v^{\;-\frac{n-1}{n}v(s)} a_1=
\pi_v^{\;\frac{v(s)}{n}} s^{-1}a'_2$, by Equation \eqref{eq:defixt}
and by the diagonal action, we have
\[
y_{a,\mb{t}}= a\,\exp(\mb{k}_{\mathfrak s})\,x_{\mb{t}}=
\pi_v^{\;\frac{v(s)}{n}}\big(R_v\,w+\sum_{j=2}^nR_v\,a_j\;e_j\big)\,.
\]
In particular, $\pi_v^{\;\frac{v(s)}{n}}w,\pi_v^{\;\frac{v(s)}{n}}
a_2e_2, \dots, \pi_v^{\;\frac{v(s)}{n}}a_ne_n$ belong to
$y_{a,\mb{t}}$. Hence $\pi_v^{\;\frac{v(s)}{n}}a'_2\;e_1$ belongs to
$y_{a,\mb{t}}$ since $a'_2\;e_1=s\;w- \sum_{j=2}^n r_j\; a_j\;e_j$ and
$s,r_2,\ldots, r_n\in R_v$.  Let $w'=\overline{r_n} \;w$, so that
$\pi_v^{\;\frac{v(s)}{n}}w'= \overline{r_n}\;(\pi_v^{\;\frac{v(s)}{n}}
w)$ belongs to the $R_v$-lattice $y_{a,\mb{t}}$, and
\[
y_{a,\mb{t}}=\pi_v^{\;\frac{v(s)}{n}}\big(R_v\,w+R_v\,w'+R_v\,a'_2\;e_1+
\sum_{j=2}^nR_v\,a_j\;e_j\big)\,.
\]
Let $s_0\in R_v$ be such that $1=\overline{r_n}\;r_n+s\,s_0$, so that
$w= r_nw'+s_0\big(a'_2\; e_1+\sum_{j=2}^n\;r_j\;a_j\;e_j\big)$ and we
can remove $R_v\,w$ in the above expression of $y_{a,\mb{t}}$.
Therefore, plugging in the expressions of $w'=\overline{r_n}
\;s^{-1}\big(a'_2\; e_1+\sum_{j=2}^n\;r_j\;a_j\;e_j\big)$ and
$a'_2=s'a_1$, we have
\begin{align*}
  y_{a,\mb{t}}=\;&
  \pi_v^{\;\frac{v(s)}{n}}\Big(R_v\,s^{-1}\big(\overline{r_n}
  \,(s'\,a_1)\, e_1+ (\overline{r_n} \;r_2)\,a_2\, e_2+\ldots+
  (\overline{r_n} \;r_{n-1})\,a_{n-1} e_{n-1}+s'(s'^{-1}\,a_{n}) e_n\big)
  \\&\quad+R_v\,(s'\,a_1)\;e_1+R_v\,a_2\,e_2+\ldots+R_v\, a_{n-1}\, e_{n-1}
  +R_v \,s'(s'^{-1}\,a_{n}) \,e_n\Big)\,.
\end{align*}
The result follows by the action $P_{\sigma_n}: e_i \mapsto e_{i+1}$
(for $i\in \intbra{1,n}$ modulo $n$) of $\sigma_n$ on the canonical
basis, and by the definition of $a'$ and $\mb{t}'$.
\cqfd

\medskip
Recall that ${\mathfrak s}= (s,\ldots, s)$ till the end of the proof
of the second claim of Assertion \eqref{item4:descripxt}. Since we
have $a\exp(\mb{k})\,\uuu_{\mb{t}}R_v^{\;n}= \exp(\mb{k}
-\mb{k}_{\mathfrak s})\,y_{a,\mb{t}}$ and by using the change of
variable $\mb{k}\in\loz_{\mathfrak s}\mapsto\mb{k}-\mb{k}_{\mathfrak
  s} \in\wt\loz_{{\mathfrak s}}$ in Equation \eqref{eq:definuloz}, we
have
\begin{equation}\label{eq:reducnuloz1}
\nu_{\mathfrak s}^\loz=
\frac{1}{\card\;\Lambda_{{\mathfrak s}}\;\card\;\loz_{{\mathfrak s}}}\;
\sum_{\mb{t}\in\Lambda_{{\mathfrak s}},\;\mb{k}\in\wt\loz_{\mathfrak s}}\int_{a\in A(\OOO_v)}
\;\delta_{\exp(\mb{k})\,y_{a,\mb{t}}}\;da\,.
\end{equation}
Similarly, by Equation \eqref{eq:definuDelta}, we have
\begin{equation}\label{eq:reducnuDelta}
\nu_{{\mathfrak s}}^\Delta=
\frac{1}{\card\;\Lambda_{{\mathfrak s}}\;\card\;\Delta_{{\mathfrak s}}}\;
\sum_{\mb{t}\in\Lambda_{\mathfrak s},\;\mb{k}\in\wt\Delta_{\mathfrak s}}\int_{a\in A(\OOO_v)}
\;\delta_{\exp(\mb{k})\,y_{a,\mb{t}}}\;da\,.
\end{equation}

Let $\wt\loz^\sharp_{\mathfrak s}$ be a strict fundamental domain for
the (free) action of $\sigma_n^{\,\ZZ}$ on $\wt\Delta_{\mathfrak s}
\ssm\{\mb{0}\}$, so that
\[
\big\{\,\mb{k}=(k_1,\ldots, k_n)\in\ZZ_0^n:
\forall\;j\in\intbra{2,n},\;\frac{v(s)}{n}\leq k_j< k_1\big\} \subset
\wt\loz^\sharp_{\mathfrak s}\subset \wt\loz_{\mathfrak s}
\]
and $\wt\Delta_{\mathfrak s}\ssm\{\mb{0}\}=\bigsqcup_{j=0}^{n-1}\;
\sigma_n^{\;j}\cdot\wt\loz^\sharp_{\mathfrak s}$ (see the picture at
the beginning of this Section for an illustration when $n=3$ of this
partition, after translating by $-\mb{k}_{\mathfrak s}$). By the standard
Gauss counting argument and by Assertion \eqref{item1_5:descripxt}
applied with $\mb{s}={\mathfrak s}= (s,\ldots, s)$, we have
\begin{equation}\label{eq:petitbord}
\frac{\card\big(\wt\loz_{\mathfrak s}\ssm\wt\loz^\sharp_{\mathfrak s}\big)}
{\card\;\loz_{\mathfrak s}}
=\bigO\Big(\frac{1}{-v(s)}\Big)\quad\text{and}\quad
\frac{\card\;\Delta_{\mathfrak s}}{\card\;\loz_{\mathfrak s}}=n
+\bigO\Big(\frac{1}{-v(s)}\Big)\,.
\end{equation}
Let
\[
\nu_{\mathfrak s}^\sharp=
\frac{1}{\card\;\Lambda_{\mathfrak s}\;\card\;\loz_{\mathfrak s}}\;
\sum_{\mb{t}\in\Lambda_{\mathfrak s},\;\mb{k}\in\wt\loz^\sharp_{\mathfrak s}}\int_{a\in A(\OOO_v)}
\;\delta_{\exp(\mb{k})\,\,y_{a,\mb{t}}}\;da\,,
\]
so that by Equation \eqref{eq:reducnuloz1} and by the left hand part
of Equation \eqref{eq:petitbord}, for every $f\in {\rm C}_c (\X)$,
we have
\begin{equation}\label{eq:passagesharp}
\nu_{\mathfrak s}^\loz(f)=\nu_{\mathfrak s}^\sharp(f)+ \bigO\Big(
\frac{\|f\|_\infty}{-v(s)}\Big)\,.
\end{equation}
By Equation \eqref{eq:sigmapermut}
and by Lemma \ref{lem:tournemanege}, we have
\[
(\sigma_n)_*\delta_{\exp(\mb{k})\,y_{a,\mb{t}}}
=\delta_{\sigma_n(\exp(\mb{k})\,y_{a,\mb{t}})}=
\delta_{\exp(\sigma_n\cdot\mb{k})\;\sigma_n(y_{a,\mb{t}})}=
\delta_{\exp(\sigma_n\cdot\mb{k})\;y_{a',\mb{t}'}}\,.
\]
By the changes of variable $a\mapsto a'$ in the integral, and
$\mb{k}\mapsto\sigma_n\cdot\mb{k}$ as well as $\mb{t}\mapsto \mb{t}'$
in the sum, we hence have
\[
(\sigma_n)_*\nu_{\mathfrak s}^\sharp=
\frac{1}{\card\;\Lambda_{\mathfrak s}\;\card\;\loz_{\mathfrak s}}\;
\;\sum_{\mb{t}\,\in\,\Lambda_{\mathfrak s},
  \;\mb{k}\,\in\,\sigma_n\cdot\,\wt\loz^\sharp_{\mathfrak s}}
\;\int_{a\in A(\OOO_v)}\;\delta_{\exp(\mb{k})\,\,y_{a,\mb{t}}}\;da\,.
\]
By iteration, we therefore have
\[
\frac{1}{n}\sum_{j=0}^{n-1}(\sigma_n^{\;j})_*\nu_{\mathfrak s}^\sharp =
\frac{1}{n\;\card\;\Lambda_{\mathfrak s}\;\card\;\loz_{\mathfrak s}}\;\;
\sum_{\mb{t}\,\in\,\Lambda_{\mathfrak s},\;\mb{k}\,\in\,\bigsqcup_{j=0}^{n-1}
  \sigma_n^{\;j}\cdot\,\wt\loz^\sharp_{\mathfrak s}}\;\int_{a\in A(\OOO_v)}
\;\delta_{\exp(\mb{k})\,\,y_{a,\mb{t}}}\;da\,.
\]
By Assertion \eqref{item1_5:descripxt}, for every $f\in {\rm C}_c
(\X)$, we have
\begin{align*}
\frac{1}{n\;\card\;\Lambda_{\mathfrak s}\;\card\;\loz_{\mathfrak s}}\;
\sum_{\mb{t}\in\Lambda_{\mathfrak s}}\int_{a\in A(\OOO_v)}\;f(y_{a,\mb{t}})\;da
=\bigO\Big(\frac{\|f\|_\infty}{(-v(s))^{n-1}}\Big)=
\bigO\Big(\frac{\|f\|_\infty}{-v(s)}\Big)\,.
\end{align*}
Since $\wt\Delta_{\mathfrak s}=\{\mb{0}\}\cup\bigsqcup_{j=0}^{n-1}\;
\sigma_n^{\;j}\cdot\wt\loz^\sharp_{\mathfrak s}$, by Equations
\eqref{eq:passagesharp} and \eqref{eq:reducnuDelta} and by the right
hand part of Equation \eqref{eq:petitbord}, the result follows.
\cqfd

\bigskip
The aim of the following sections will be to prove that the probability
measures 
%%
%\todo{\tiny see if we do manage this generality}
%%
$\nu_{\mb{s}}^\loz$ weak-star converge as $\mb{s}\ra +\infty$ (with
appropriate conditions the components of $\mb{s}$ that are satisfied
if they are all equal and have absolute values equal to a multiple of
$n$) to the $G$-homogeneous measure ${\tt m}_{\X}$ on $\X$,
renormalized to be a probability measure, see Theorem
\ref{theo:main}. In the particular case when $\mb{s}=(s,s,\ldots,s)$,
this will imply Theorem \ref{theo:mainuintro2} by the following
lemma. We denote by $c_{K,n}$ the constant defined in the statement of
Theorem \ref{theo:mainuintro2}.

\blemm\label{lem:deducmainuintro2} Assume that ${\displaystyle
  \lim_{s\ra\infty,\;v(s)\, \in\, n\ZZ}} \;\;\nu_{(s,\ldots,s)}^\loz
=\frac{{\tt m}_{\X}}{\|{\tt m}_{\X}\|}$.  Then
\[
\lim_{s\ra\infty,\;v(s)\, \in\, n\ZZ} \;\;
\frac{c_{K,n}}{(\varphi_v(s)\, \log_{q_v}|s|\,)^{n-1}}\;
\sum_{\mb{t}\in \Lambda_{(s,\ldots,s)}}\; \ov\mu_{\uuu_{\mb{t}} R_v^{\;n}}\;
=\;{\tt m}_{\X_1}\,.
\]
\elemm

\dem For every $s\in R_v$ such that $v(s)\, \in\, n\ZZ$, let
${\mathfrak s}= (s,s,\ldots,s)\in (R_v\ssm\{0\})^{\;n-1}$.  By the
assumption of the lemma and by a finite average using the last claim
of Proposition \ref{prop:descripxt} \eqref{item4:descripxt} (the only
point where we use the symmetry of ${\mathfrak s}$), we have
${\displaystyle \lim_{s\ra\infty, \; v(s)\, \in\, n\ZZ}} \;\;
\nu_{{\mathfrak s}}^\Delta = \frac{{\tt m}_{\X}} {\|{\tt m}_{\X}\|}$.
By the first claim of Proposition \ref{prop:descripxt}
\eqref{item4:descripxt} and by Equation \eqref{eq:cardLambdasss} on
the right, we have $\nu_{{\mathfrak s}}^\Delta =\frac{1}
{(\varphi_v(s))^{n-1}} \; \sum_{\mb{t}\in \Lambda_{{\mathfrak s}}}
\;\nu_{x_{\mb{t}}}$. By Proposition \ref{prop:descripxt}
\eqref{item1:descripxt}, we have $\tau_{x_{\mb{t}}} =-v(s)=
\log_{q_v}|s|$ for all $\mb{t}\in\Lambda_{{\mathfrak s}}$. Hence by
the last claim of Lemma \ref{lem:massmostoncompcore} (and the
definition of the weak-star convergence), we have
\begin{equation}\label{eq:maindansX}
  \lim_{s\ra\infty, \;v(s)\, \in\, n\ZZ}\;
  \frac{1}{c_n\,(\varphi_v(s)\,\log_{q_v}|s|\,)^{n-1}}
\sum_{\mb{t}\in \Lambda_{{\mathfrak s}}} \mu_{x_{\mb{t}}}=\frac{{\tt m}_{\X}}
    {\|{\tt m}_{\X}\|}\,.
\end{equation}
Again by averaging, this time over the compact probability space
$(\OOO_v^{\,\times}/R_v^{\,\times}, \frac{q_v(q-1)}{q_v-1} \vol'_v)$
defined at the end of Section \ref{subsec:homogspalat}, and using
Equations \eqref{eq:desintegovmumu} and \eqref{eq:mXmX1}, we have
\[
\lim_{s\ra\infty,\;v(s)\, \in\, n\ZZ}
\;\;\frac{1}{c_n\,(\varphi_v(s)\,\log_{q_v}|s|\,)^{n-1}}
\sum_{\mb{t}\in \Lambda_{{\mathfrak s}}} \ov\mu_{x_{\mb{t}}} = \frac{{\tt
    m}_{\X_1}}{\|{\tt m}_{\X_1}\|}\,.
\]
By Equation \eqref{eq:totmassmX_1} giving the value of $\|{\tt
  m}_{\X_1}\|$ and by Proposition \ref{prop:propriinvar}
\eqref{item5:propriinvar} giving the value of $c_n$, we have
\[
\frac{\|{\tt m}_{\X_1}\|}{c_n}=\frac{(n-1)!\,(q_v-1)}{q_v\,(q-1)}
\prod_{i=1}^{n-1}\frac{\zeta_v(-i)}{q_v^{\;i}-1}\,.
\]
This is exactly the value of the constant $c_{K,n}$ defined in the
statement of Theorem \ref{theo:mainuintro2}. The result follows. \cqfd

\medskip
\rem Equation \eqref{eq:maindansX} says that the distribution
statement in $\X$ analogous to the one of Theorem
\ref{theo:mainuintro2} in $\X_1$ will also follow from Theorem
\ref{theo:main}.

\section{The high entropy method for equidistribution problem}
\label{sec:entropytoequidistrib}

In this section, we explain how we are going to use entropy technics
in homogeneous dynamics of diagonal actions (that are currently more
and more used, see for instance \cite{EinLin10, LimSaxSha18,
  EinLinWar22, KimLimPau23}) in order to prove the equidistribution of
our families of measures. We start by a brief reminder of entropy
theory.

Let $(X',\mu')$ be a Borel probability space, $\phi:X'\ra X'$ a
measurable map, and $\P,\P'$ finite measurable partitions of $X'$.

We denote by $\phi^{-1}\P=\{\phi^{-1}(B):B\in \P, \;\phi^{-1}(B)\neq
\emptyset\}$ the pull-back partition and by $\P\vee\P''=\{B\cap B':
B\in\P,\;B'\in \P',\;B\cap B'\neq\emptyset\}$ the joint partition.
%and for every $M\in\NN$ by $\P^{M}
%=\bigvee_{i=0}^{M-1}\phi^{-i}\P$ the iterated pull-back joint.
Using the convention $0 \log_{q_v} 0 =0$, the {\it entropy} of the
partition $\P$ with respect to $\mu'$ is
\[
H_{\mu'}(\P) = - \sum_{P\in \P} \mu'(P) \log_{q_v} \mu'(P) \;\in [0, \infty[\;.
\]
The usual definition of the entropy of a partition uses the Neperian
logarithm $\ln$ instead of $\log_{q_v}$, but the above convention will
be technically easier in this paper.
%Recall that $H_\mu(\P)\leq \log_{q_v}(\card\P)$.
We have $H_{\phi_*\mu'}(\P)=H_{\mu'}(\phi^{-1}\P)$.  We have the
following concavity properties of the entropy of a partition as a
function of the measure.

\blemm\label{lem:entromino} (David-Shapira \cite[Lem~3.4]{DavSha20}) Let
$(X',\mu'),\phi,\P$ be as above.
\begin{enumerate}
\item\label{item1:entromino}
For all $M\leq N$ in $\NN\ssm\{0\}$, we have
\[
\frac{1}{M}
H_{\frac{1}{N}\sum_{i=0}^{N-1}(\phi^i)_*\mu'}\Big(\bigvee_{i=0}^{M-1}\phi^{-i}\P\Big)
\geq \frac{1}{N} H_{\mu'}\Big(\bigvee_{i=0}^{N-1}\phi^{-i}\P\Big)- \frac{M}{N}
\log_{q_v}\card\;\P\,. 
\]
\item\label{item2:entromino} Let $(\Omega, \omega)$ be a probability
  space and let $x\mapsto \mu'_x$ be a measurable map from $\Omega$ to
  the space of probability measures on $X'$ such that $\mu'=
  \int_{x\in\Omega} \mu'_x\; d\omega(x)$. Then we have $H_{\mu'}(\P)
  \geq \int_{x\in\Omega} H_{\mu'_x}(\P)\;d\omega(x)$.
  \cqfd
\end{enumerate}
\elemm

If $\phi$ preserves the measure $\mu'$, the {\it (dynamical) entropy}
of $\phi$ with respect to $\mu'$ is defined by $h_{\mu'}(\phi)=
\sup_{\P}\; h_{\mu'}(\phi,\P)$ where the least upper bound is taken
over all finite measurable partitions $\P$ of $X'$ and
\begin{equation}\label{eq:defientrorelpart}
h_{\mu'}(\phi,\P)=\lim_{M\ra+\infty}
\frac{1}{M}H_{\mu'}\Big(\bigvee_{i=0}^{M-1}\phi^{-i}\P\Big)\,.
\end{equation}

\medskip
The following result says that the homogeneous measure $m_{\X}$ is
after renormalisation the unique probability measure of maximal
entropy on the space $\X$ for the transformation
$\aaa=\big(\begin{smallmatrix} \pi_v^{\;n-1} & 0\\0 &
  \pi_v^{\;-1}I_{n-1} \end{smallmatrix}\big)$ given by Equation
\eqref{eq:defaaa}.

\btheo\label{theo:EL} {\bf (Einsiedler-Lindenstrauss)} Let $\nu$ be an
$\aaa$-invariant probability measure on $\X$. Then $h_\nu(\aaa) \leq
h_{\frac{m_{\X}}{\|m_{\X}\|}}(\aaa)=n(n-1)$ with equality if and only
if $\nu=\frac{m_{\X}}{\|m_{\X}\|}$.
%%
%\todo{\tiny version effective : autre projet}
%%
\etheo

We will apply this theorem in Section \ref{sec:lowboundentrop} to
every weak-star accumulation point $\nu$ of the measures
$\nu_{\mb{s}}^\loz$ as $\mb{s}$ tends appropriately to $+\infty$.
Since the space $\X$ is not compact, we will first need to prove that
$\nu$ is a probability measure (see the arguments in Section
\ref{sec:noescapmass}), and then that its entropy $h_\nu(\aaa)$ is
equal to $h_{\frac{m_{\X}}{\|m_{\X}\|}}(\aaa)$ (see Section
\ref{sec:lowboundentrop}).

\medskip
\dem Let ${\mathbb G}$ be the algebraic group $\SL_n$ over the local
field $K_v$ and let $G=\SL_n(K_v)$ be its locally compact group of
$K_v$-points, so that $\Ga=\SL_n(R_v)$ is a lattice in $G$.  Let
\[
G^-=\{g\in G:\lim_{i\ra-\infty }\aaa^ig\aaa^{-i}=I_n\}=\big\{\big(
\begin{smallmatrix}1 & 0\\b & I_{n-1}\end{smallmatrix}\big):
b\in K_v^{\;n-1}\big\}\]
and
\[
G^+=\;^tG^-=\{g\in G:\lim_{i\ra+\infty }
\aaa^ig\aaa^{-i}=I_n\}=\big\{\big(\begin{smallmatrix}1 & ^tb\\0 &
  I_{n-1} \end{smallmatrix}\big): b\in K_v^{\;n-1}\big\}
\]
be respectively the unstable and stable horospherical groups of $\aaa$
in $G$. By \cite[Prop~4.11]{BorTit65}, the groups $G^-$ and $G^+$
generate a normal subgroup $H$ of $G$. It is well known that $H=G$
when $n=2$. Hence $H$ contains the copies of $\SL_2(K_v)$ with upper
and lower unipotent subgroups contained in $G^-$ and $G^+$
respectively. Therefore $H$ contains the diagonal subgroup $A$ of $G$,
thus contains properly the center of $G$, hence is equal to $G$ since
$\PSL_n(K_v)$ is simple. By \cite[Th.~7.10]{EinLin10}, the normalized
Haar measure $\frac{m_{\X}}{\|m_{\X}\|}$ of the homogeneous space
$G/\Ga$ is hence the unique measure of maximal entropy on $G/\Ga$ for
the left action of $\aaa$.

The entropy of $\aaa$ with respect to the homogeneous measure of
$\frac{m_{\X}}{\|m_{\X}\|}$ is well-known (see for instance \cite[\S
  7.8]{EinLin10}) to be the logarithm (in basis $q_v$ for entropy
computations in nonarchimedian local fields with residual fields of
order $q_v$) of the unstable Jacobian of $\aaa$. That is, with
$\underline{\uuu}^-$ the strict lower triangular linear subspace of
the Lie algebra $\sss\lll_n(K_v)$ of $\SL_n(K_v)$, with basis the
family of elementary matrices $(E_{i,j})_{1\leq j<i\leq n}$, since we
have $\Ad \aaa\; (E_{i,j}) = \pi_v^{\;-n}E_{i,j}$ if $j=1$ and $\Ad
\aaa \;(E_{i,j})= E_{i,j}$ otherwise, we have
\[
h_{\frac{m_{\X}}{\|m_{\X}\|}}(\aaa)=
\log_{q_v}|\det (\Ad \aaa)_{\mid\,\underline{\uuu}^-}|=
\log_{q_v}\Big|\prod_{j=2}^n\pi_v^{\;-n}\,\Big|=n(n-1)
\,.\;\;\;\Box
\]

\section{Constructing high entropy partitions from dynamical
neighborhoods}
\label{sec:countdynball}

In this section, using the contraction and dilation properties of the
action of the diagonal element $\aaa=\big(\begin{smallmatrix}
\pi_v^{\;n-1} & 0\\0 & \pi_v^{\;-1}I_{n-1} \end{smallmatrix}\big)$ on its
unstable and stable horospherical subgroups $G^-$ and $G^+$, we give a
construction of good measurable partitions in the homogeneous space
$\X$, that will turn out in Section \ref{sec:lowboundentrop} to be
well adapted in order to obtain entropy lower bounds of
$\aaa$-invariant measures. This construction is essentially due to
\cite[Lem.~4.5]{EinLinMicVen12} in dimension $2$ (see
\cite[Lem.~2.9]{DavSha18} correcting a small inaccuracy in
\cite{EinLinMicVen12}) and to \cite[Lem.~3.7]{DavSha20} for any
dimension, see also \cite[\S 2.3]{LimSaxSha18}, \cite[\S
3.3]{KimKimLim21}, all these references in the real case, and
\cite[\S 6.1]{KimLimPau23} in the function fields case.

\subsection{Dynamical neighbourhoods in $\SL_n(K_v)$}
\label{subsec:dynneig}

We first define the dynamical neighbourhoods of the identity element
$I_n$ in $\SL_n(K_v)$ that we will consider.

We denote by $\|\;\|:\M_{n}(K_v)\ra [0,+\infty[$ the ultrametric norm
on $\M_{n}(K_v)$ defined by $(x_{ij})_{1\leq i,j\leq n}\mapsto
\max_{1\leq i,j\leq n} |\,x_{ij}\,|$, which is, since the absolute
value of $K_v$ is ultrametric, a submultiplicative norm on the
$K_v$-algebra $\M_{n}(K_v)$.

For all $\ell,N\in\ZZ$, let
\begin{align}
W_{\ell,N} & =
\{w=(w_{i,j})_{1\leq i,j\leq n}\in\M_n(K_v) : \|w\| \leq q_v^{\;-\ell}
\mbox{ and } \forall\;i\in\llbracket2,n\rrbracket,\;
|w_{i,1}| \leq q_v^{\;-(\ell+nN)}\},
\nonumber\\ & = \{w\in\M_n(\pi_v^{\;\ell} \OOO_v)  : \forall\;i\in
\llbracket2,n\rrbracket,\; w_{i,1} \in \pi_v^{\; \ell+nN} \OOO_v\},
\nonumber \\ B_{\ell,N} & = (I_n+W_{\ell,N}) \cap \SL_n(K_v)\,.
\label{eq:defiBellN}
\end{align}
We also define $W_{\ell}=W_{\ell,0}=\M_n(\pi_v^\ell\OOO_v)$ and
$B_\ell=B_{\ell,0}$. For all $\ell,\ell',N \in \ZZ$, by the
ultrametric inequalities, we have
\begin{equation}\label{eq:subaddsubmult}
W_{\ell,N} + W_{\ell',N} \subset W_{\min(\ell,\ell'), N}\quad\text{and}\quad
W_{\ell,N} W_{\ell',N} \subset W_{\ell+\ell',N}\,.
\end{equation}
We also have the following decreasing properties $W_{\ell,N+1}\subset
W_{\ell,N}$, $B_{\ell,N+1}\subset B_{\ell,N}$ in the parameter $N$
and, in the parameter $\ell$,
\begin{equation}\label{eq:voisdyncroiss}
W_{\ell+1,N}\subset W_{\ell,N},\quad B_{\ell+1,N}\subset B_{\ell,N},\quad
\bigcap_{\ell\in\NN}W_{\ell,N}=\{0\},\quad
\bigcap_{\ell\in\NN}B_{\ell,N}=\{I_n\}\,.
\end{equation}
%For all $\ell>0$ and $w \in W_{\ell,N}$, we also obtain the equality
%$|\det(I_n+w)|=1$, that is, $I_n+W_{\ell,N} \subset \GL_n^1(\OOO_v)$.
%This is not sufficient for $I_n+w$ to be in $\SL_n(K_v)$ but one can
%notice that the quotient topological group $\GL_n^1(K_v)/\SL_n(K_v)
%\simeq \OOO_v^\times$ is compact.
Since the multiplication by an element of $\OOO_v^\times$ preserves
the absolute value on $K_v$, for every $a\in A(\OOO_v)$, we have
\begin{equation}\label{eq:invarWellNpara}
  a\,W_{\ell,N}\,a^{-1}=W_{\ell,N}\qquad\text{and}\qquad
  a\,B_{\ell,N}\,a^{-1}=B_{\ell,N}
\end{equation}
The action by conjugation of the transformation $\aaa$ on these
dynamical balls $W_{\ell,N}$ and $B_{\ell,N}$ satisfy the following
contraction/dilation properties: For all $\ell,\ell',N \in \ZZ$, we
have
\begin{align}\label{eq:dynamic_a_on_balls}
\aaa^{\ell'} W_{\ell, N} \,\aaa^{-\ell'} &=  \Big\{w\in\M_n(K_v) :
\begin{array}{l} \forall\;i,j \in\llbracket2,n\rrbracket,\;
w_{1,1}, w_{i,j}\in\pi_v^{\; \ell}\OOO_v,\\
w_{i,1}\in\pi_v^{\;\ell+(N-\ell')n}\OOO_v, \;
w_{1,j}\in\pi_v^{\; \ell + n\ell'} \OOO_v\end{array}\Big\} \nonumber
\\ &\subset  W_{\min\{\ell, \,\ell+n\ell'\},\; N-\ell'} \nonumber
\\ \mbox{hence }\quad \aaa^{\ell'} B_{\ell, N} \,\aaa^{-\ell'}& \subset 
B_{\min\{\ell,\,\ell+n\ell'\},\; N-\ell'}\,.
\end{align}
The dynamical neighbourhoods $W_{\ell,N}$ and $B_{\ell,N}$ satisfy the
following three elementary lemmas.
%%
%\todo{\tiny Lem A.1-A.3 \cite{DavSha20} nonarchimédiens}
%%
The first one says that the balls
$B_{\ell,N}$ are invariant upon taking inverses.

\blemm\label{lem:inverse_dyn_balls} Let $N\in\NN$, $\ell\in\NN\ssm
\{0\}$ and $w \in W_{\ell,N}$. Then $(I_n+w)^{-1} \in I_n - w +
W_{2\ell,N}$. In particular, we have $(B_{\ell,N})^{-1} = B_{\ell,N}$.
\elemm

\dem By Equations \eqref{eq:subaddsubmult} and
\eqref{eq:voisdyncroiss}, we have $w^i \in W_{i\ell,N}$ and
$\lim_{i\ra+\infty}w^i=0$ since $\ell>0$. Hence $I_n+w$ is invertible
with $(I_n+w)^{-1} = I_n - w + \sum_{i=2}^{\infty} (-1)^i w^i$. By
Equation \eqref{eq:subaddsubmult} and since $W_{2\ell, N}$ is closed,
we have $\sum_{i=2}^{\infty} (-1)^i w^i \in W_{2\ell,N}$. In
particular we have the inclusion $(B_{\ell,N})^{-1} \subset
B_{\ell,N}$, and equality holds by taking the inverses.  \cqfd

\medskip
The next lemma is a version, for dynamical balls in $\X$ of the form
$B_{\ell,N}\,x$ centred at any point $x\in\X$, of the intersection
property of ultrametric balls.

\blemm\label{lem:reciprocal_inclusion_dyn_balls} For all $N\in\NN$,
$\ell\in\NN\ssm\{0\}$ and $x,y \in \X$ with $B_{\ell,N} \,x \cap
B_{\ell,N} \,y \neq \emptyset$, we have $B_{\ell,N} \,x = B_{\ell,N}\,
y$.
\elemm

\dem First notice the inclusion
\begin{align}
B_{\ell,N}B_{\ell,N} & \subset ((I_n+W_{\ell,N})(I_n+W_{\ell,N}))\cap\SL_n(K_v)
\nonumber\\ & \subset (I_n+W_{\ell,N}+W_{2\ell,N})\cap\SL_n(K_v)
\subset B_{\ell,N}\,.\label{eq:BellBelldansBell}
\end{align}
Let $g,h \in B_{\ell,N}$ be such that $gx=hy$. Using Lemma
\ref{lem:inverse_dyn_balls} and the latter inclusion, we have
\begin{align*}
B_{\ell,N} \,x & =  B_{\ell,N} \,g^{-1} h y \subset
B_{\ell,N} (B_{\ell,N})^{-1} B_{\ell,N}\, y =
B_{\ell,N} B_{\ell,N} B_{\ell,N} \, y \subset B_{\ell,N}\, y\,.
\end{align*}
By symmetry, the result follows.
\cqfd

\medskip
The final lemma gives a quantitative covering property for any
dynamical ball $B_{\ell,N}$ in $\SL_n(K_v)$ by smaller dynamical balls
$B_{\ell+\ell',N}g_i$.

\blemm\label{lem:cover_dynamic_balls_SLn} Let $N\in\NN$ and
$\ell,\ell'\in\NN\ssm\{0\}$ with $\ell' \leq \ell$. Let $S \subset
B_{\ell,N}$. Then, there exist an integer $C \leq (q_v^{\;
\ell'})^{n^2}$ and matrices $g_1, \ldots, g_C \in S$ such that
\[
S \subset \bigsqcup_{i=1}^C B_{\ell+\ell',N} \, g_i\,.
\]
\elemm

\dem We may assume that $S$ is nonempty, otherwise $C=0$ works. As a
preliminary remark, let us prove that for all integers $\ell,\ell'
>0$, there exist an integer $C \leq (q_v^{\; \ell'})^{n^2}$ and points
$w_1, \ldots, w_C \in W_{\ell,N}$ such that
\begin{equation}\label{eq:fact_cover_dynamic_balls_Mn}
W_{\ell,N} = \bigsqcup_{i=1}^{C} (w_i + W_{\ell+\ell', N})\,.
\end{equation}
Indeed, recall that Equation \eqref{eq:changevarHaar} when $n=1$ gives
$\vol_v(\pi_v^{\;\ell+\ell'}\OOO_v)=q_v^{\;-\ell'}\vol_v(\pi_v^{\;\ell}
\OOO_v)$. Let $\{x_i: i \in I\}$ be a set of representatives of the
classes in $\pi_v^{\;\ell} \OOO_v/ \pi_v^{\;\ell+\ell'} \OOO_v$, so
that we have a partition $\pi_v^{\; \ell}\OOO_v = \bigsqcup_{i\in I}
(x_i+\pi_v^{\; \ell+\ell'}\OOO_v)$. Furthermore, by the invariance of
$\vol_v$ under translations, we have $\card(I) \leq \vol_v
(\pi_v^{\;\ell}\OOO_v)/\vol_v(\pi_v^{\;\ell+\ell'}\OOO_v)=q_v^{\;\ell'}$.
The same argument replacing $\ell$ by $\ell+nN$ proves that
$\pi_v^{\;\ell+nN} \OOO_v$ can be covered by at most $q_v^{\;\ell'}$
pairwise disjoint translates of the ball $\pi_v^{\;\ell+\ell'+nN}
\OOO_v$. Equation \eqref{eq:fact_cover_dynamic_balls_Mn} follows by
applying this construction for each matrix entries.

Now, take $C \leq (q_v^{\; \ell'})^{n^2} $ and $w_1, \ldots, w_C \in
W_{\ell,N}$ as in Equation \eqref{eq:fact_cover_dynamic_balls_Mn}. We
obtain a partition
\[
S = \bigsqcup_{i=1}^C (I_n+w_i+W_{\ell+\ell', N})\cap\SL_n(K_v)\cap S\,.
\]
Up to decreasing $C$, we may assume that the set $S_i=(I_n + w_i +
W_{\ell+\ell', N}) \cap \SL_n(K_v) \cap S$ is nonempty for every
$i\in\llbracket1, C\rrbracket$. Let us fix an element $g_i \in S_i$
and let us prove that we have $(I_n+w_i+W_{\ell+\ell', N})\,g_i^{-1}
\cap \SL_n(K_v) \subset B_{\ell+\ell',N}$. By Equations
\eqref{eq:subaddsubmult} and \eqref{eq:voisdyncroiss}, since
$\ell'\leq \ell$ and $I_n\in W_{0,N}$, we have
\[
w_i^{\;2}\in W_{\ell,N}W_{\ell,N}\subset W_{2\ell,N}\subset
W_{\ell+\ell',N}\quad\text{and}\quad W_{\ell+\ell',N}W_{2\ell,N}
\subset W_{3\ell+\ell',N}\subset W_{\ell+\ell',N}\,,
\]
\[
(I_n+w_i)W_{2\ell,n}\subset (W_{0,N}+W_{\ell,N})W_{2\ell,N}\subset
W_{0,N}W_{2\ell,N}\subset W_{2\ell,N}\subset W_{\ell+\ell',N}\,.
\]
By Lemma \ref{lem:inverse_dyn_balls}, we know that $g_i^{-1} \in
I_n-w_i+W_{2\ell,N}$. We hence have
\begin{align*}
(I_n+w_i+W_{\ell+\ell', N})\,g_i^{-1} &
\subset (I_n+w_i+W_{\ell+\ell', N})(I_n-w_i+W_{2\ell,N})
\\ & \subset I_n-w_i^2 +(I_n + w_i) W_{2\ell,N} +
W_{\ell+\ell',N}(I_n - w_i) + W_{\ell+\ell', N}W_{2\ell,N}
\\ & \subset I_n + W_{\ell+\ell',N}\,.
\end{align*}
We obtain the inclusions $S_i\subset B_{\ell+\ell',N} \, g_i $ by
taking the intersection with $\SL_n(K_v)$, so that $S \subset
\bigcup_{i=1}^C B_{\ell+\ell',N}\,g_i$. Up to decreasing $C$, we may
assume that this intersection is disjoint by using Lemma
\ref{lem:reciprocal_inclusion_dyn_balls}. This concludes the proof.
\cqfd

\subsection{Dynamical partitions in $\X$}
\label{subsec:dynpart}

We now construct measurable partitions of $\X$, that will be useful
for entropy lower bounds computations. We start with a systole
minoration result for $R_v$-lattices that are close enough to the
standard ``cubic'' $R_v$-lattice. Recall that the map $\exp$ is
defined just before Subsection \ref{subsec:homogdiag} and the systole
function in Subsection \ref{subsec:systole}.

\blemm\label{lem:minosystole}
%%
%\todo{\tiny Lem 3.9 \cite{DavSha20} nonarchimédien}
%%
Let $\ell,d\in\NN\ssm\{0\}$, ${\bf k} =(k_1,\ldots,k_d)\in\ZZ^d$ and
$g\in \M_{d}(K_v)$ be such that $\|g\|\leq q_v^{\,-\ell}$.
Then if $L=(I_d+\exp ({\bf k})\, g \,\exp( -{\bf k})) R_v^{\;d}$, we
have $\sys(L)\geq 1-q_v^{\,-\ell}$.
\elemm

\dem By the equality case of the ultrametric triangular inequality and
since $\|g\|< 1$, note that $|\det (I_d+g)|= \max_{\sigma\in \Scal_d}
\prod_{1\leq i\leq d}|(I_d+g)_{i\sigma(i)}|= 1$, hence $I_d+g$ belongs
to $G_1$ and so does $\exp ({\bf k})(I_d+g)\exp(-{\bf k})$. Therefore
$L$ is indeed an $R_v$-lattice in $K_v^{\;d}$ and $\covol (L)=\covol
(R_v^{\;d})$ by Equation \eqref{eq:changevarcovol}.

By Equation \eqref{eq:defisys}, assume for a contradiction that
$\sys(L)=\min_{w\in L\,\ssm\,\{0\}}\|w\|<1-q_v^{\,-\ell}$. Let ${\bf x}
=(x_1,\ldots,x_d)\in R_v^{\;d}\ssm\{0\}$ be such that $\|(I_d+\exp
({\bf k})\, g \,\exp( -{\bf k})){\bf x}\,\| <1-q_v^{\,-\ell}$.  For
every $i_0\in\intbra{1,d}$ such that $x_{i_0}\neq 0$, by computing the
$i_0$-th coordinate of the vector $(I_d+\exp ({\bf k})\, g \,\exp(
-{\bf k})){\bf x}$, by the ultrametric triangular inequality and since
$\|g\|\leq q_v^{\,-\ell}$, we then have
\begin{align*}
1-q_v^{\,-\ell}&>\big|x_{i_0}+\sum_{j=1}^d x_j\,
\pi_v^{k_j-k_{i_0}}g_{i_0j}\big|\geq
|x_{i_0}|-\max_{1\leq j\leq d}\; |x_j|\,|g_{i_0j}|\;q_v^{-k_j+k_{i_0}}\\ & \geq
|x_{i_0}|-q_v^{\,-\ell}\max_{1\leq j\leq d} \;|x_j|\;q_v^{-k_j+k_{i_0}}\,.
\end{align*}
Noting that $|x_{i_0}|\geq 1$ since $x_{i_0}\in R_v\ssm \{0\}$, and
since $q_v^{\,\ell}-1\geq 0$, we have $|x_{i_0}|(q_v^{\,\ell}-1)\geq
q_v^{\,\ell}-1$, thus $|x_{i_0}|\leq q_v^{\,\ell}\big(|x_{i_0}|-1
+q_v^{\,-\ell}\big)$. Therefore
\[
|x_{i_0}|\;q_v^{-k_{i_0}}\leq
q_v^{\,\ell}\big(|x_{i_0}|-1+q_v^{\,-\ell}\big)\,q_v^{-k_{i_0}}
<\max_{1\leq j\leq d} |x_j|\;q_v^{-k_j}\,.
\]
Thus, there exists $i_1\neq i_0$ such that $|x_{i_0}|\;q_v^{-k_{i_0}}<
|x_{i_1}|\;q_v^{-k_{i_1}}$. By iteration, and since there is no
strictly increasing  sequence in the finite subset
$\{|x_i|\;q_v^{-k_i}:i\in\intbra{1,d}, x_i\neq 0\}$ of
$[0,+\infty[\,$, we obtain a contradiction.
\cqfd

\medskip
Let us fix $\mb{s}=(s_2,\ldots,s_n) \in (R_v\ssm \{0\})^{\;n-1}$
satisfying Equation \eqref{eq:proprisbold}: There exists a permutation
$\sigma$ of $\llbracket 2,n \rrbracket$ such that
$s_{\sigma^{-1}(2)}\mid s_{\sigma^{-1}(3)}\mid\ldots\mid \mb{s}_*
=s_{\sigma^{-1}(n)}$ and
\begin{equation}\label{eq:propmbsstar}
v(\mb{s}_*)=\min_{2\leq i\leq n}\;v(s_i)\;\in n\,\ZZ\,.
\end{equation}
%We denote by $\delta_v$ the diameter of $K_v / R_v$ for the quotient
%distance, which satisfies $\delta_v\geq 1$.
The next lemma gives a cardinality estimate for the number of
$R_v$-lattices whose equidistribution we want to study, that belong to
a small dynamical ball of $\X$.

\blemm\label{lem:counting_latpts_dynballs}
For every $x \in \X = G / \Ga$, every integer $\ell > \max \{0, -
\log_{q_v}(\sys(x))\}$, every ${\bf k}=(k_1,\ldots, k_n) \in
\Delta_{\mb{s}}$ and every $N\in\llbracket 0,
\frac{-v(\mb{s}_*)-\ell+k_1-\max_{2 \leq i \leq n} k_i}{n}\rrbracket$, we
have
%%
%\todo{\tiny Lem 3.8 \cite{DavSha20} nonarchimédien}
%%
\[
\card\big( \{ \mb{t} + R_v^{\; n-1} \in \Lambda_{\mb{s}}:\exp(\mb{k})
\,\uuu_{\mb{t}} \,\Ga \in B_{\ell, N} \;x\} \big) \leq
2^{n-1}\,q_v^{\,-\ell(n-1)-v(\mb{s}_*)(n-1)-n(n-1)N}\,.
\]
\elemm

\dem Fix $x$, $\ell$, $\mb{k}$ and $N$ as in the lemma. Let $g\in G$
be such that $x=g\Ga$. For simplicity, we fix a lift
$\widetilde{\Lambda_{\mb{s}}}$ of $\Lambda_{\mb{s}}$ in $K^{n-1}$,
hence having the same cardinality $\prod_{i=2}^n\varphi_v(s_i)$ as
$\Lambda_{\mb{s}}$.
%By the definition of $\delta_v$ and a homothety argument, for every $s
%\in R_v \ssm \{0\}$, the restriction to the closed ball $B(0,|s|
%\delta_v)$ in $K_v$ of the canonical projection $K_v \to K_v/s R_v$ is
%surjective. Hence $B(0,|s|\delta_v)\cap R_v$ maps onto $R_v/s R_v$ by
%the canonical projection $R_v \to R_v/s R_v$. Setting $\delta =
%{\displaystyle \max_{\mb{t}=(t_2,\dots,t_n) \,\in\,
%    \widetilde{\Lambda_{\mb{s}}}} \;\;\max_{2\leq i\leq
%    n}}\;|t_i|\;\in q_v^{\,\ZZ}$, we may then choose
%$\widetilde{\Lambda_{\mb{s}}}$ so that $\delta \leq \delta_v$.

We want to evaluate the number of points $\mb{t} \in
\wt{\Lambda_{\mb{s}}}$ such that $\exp(\mb{k})\uuu_{\mb{t}} \Ga \in
B_{\ell, N} g \Ga$, in other words such that there exist $\ga \in \Ga$
and $h \in B_{\ell,N}$ verifying $h^{-1}\exp(\mb{k}) = g \ga
\uuu_{-\mb{t}}$.  Since $B_{\ell,N}^{-1} = B_{\ell,N}$ by Lemma
\ref{lem:inverse_dyn_balls}, we have to bound from above the
nonnegative quantity
\[
\card\{\mb{t}\in\wt{\Lambda_{\mb{s}}} : \exists \ga\in\Ga, \;
g \ga\uuu_{-\mb{t}} \in B_{\ell,N} \exp(\mb{k})\}\,.
\]
We may assume that this quantity is nonzero.  Let us take $\mb{t} =
(\frac{r_2}{s_2}, \ldots, \frac{r_n}{s_n}) \in
\widetilde{\Lambda_{\mb{s}}}$ and $\gamma \in \Ga$ whose column
matrices are denoted by $\ga_1, \ldots, \ga_n$. Recall the notation
$(e_1,\ldots, e_n)$ for the canonical $K_v$-basis of $K_v^{\;n}$. Let
us consider the constraints on the columns in the condition $g \ga
\uuu_{-\mb{t}} \in (I_n+W_{\ell,N}) \exp(\mb{k})$, which is equivalent
to the condition $g \ga \uuu_{-\mb{t}} \in B_{\ell,N} \exp(\mb{k})$
since $\det(g\ga\uuu_{-\mb{t}}\exp(-\mb{k}))=1$. We obtain the
equivalent system of conditions
\begin{align}
g \ga_1 - \sum_{i=2}^n \frac{r_i}{s_i} g \ga_i &\; \in\; \pi_v^{-k_1} e_1+
(\pi_v^{\ell-k_1} \OOO_v) \times (\pi_v^{\ell + n N-k_1} \OOO_v)^{n-1}\,,
\label{eq:system_count_latpt_dynball_1} \\
\forall\; i \in \llbracket 2, n\rrbracket, \quad g \ga_i &\; \in\;
\pi_v^{-k_i}e_i+(\pi_v^{\ell-k_i} \OOO_v)^n\,.
\label{eq:system_count_latpt_dynball_2}
\end{align}
%(In order to obtain an equivalence with $g \gamma \uuu_{-\mb{t}} \in
%B_{\ell,N} \exp(\mb{k})$, one would have to add the condition
%$\det(g \gamma \uuu_{-\mb{t}} \exp{-\mb{k}})=1$).
Since $\ell > -\log_{q_v}(\sys(x))$, the lattice $x = g R_v^{\; n}$
contains at most one point in each translate of $(\pi_v^{\; \ell}
\OOO_v)^n$. As seen in the proof of Equation
\eqref{eq:fact_cover_dynamic_balls_Mn}, for every $i \in\llbracket
2,n\rrbracket$, since $k_i \geq 0$ as $\mb{k}\in\Delta_{\mb{s}}$, we
can cover (any translate of) $(\pi_v^{\; \ell-k_i} \OOO_v)^n$ by at
most $q_v^{\; k_i n}$ pairwise disjoint translates of $(\pi_v^\ell
\OOO_v)^n$. Hence for each fixed $i\in\intbra{2,n}$, the condition on
$\gamma_i \in R_v^{\; n}$ given in Equation
\eqref{eq:system_count_latpt_dynball_2} has at most $q_v^{\; k_i n}$
solutions, and this set of solutions is independent of $\mb{t}$.

%By the definition of $\delta$, for every $i\in\intbra{2,n}$, we have
%\[
%\frac{r_i}{s_i}\in \pi_v^{-\log_{q_v}|\frac{r_i}{s_i}|} \OOO_v^\times
%\subset \pi_v^{-\max\{0, \,\log_{q_v}\delta\}}\OOO_v\,.
%\]
%Let $\alpha=\max\{0, \log_{q_v}\delta\}+\max_{2\leq j\leq n} k_j$.  It
%follows from Equation \eqref{eq:system_count_latpt_dynball_2} that
%\[
%\sum_{i=2}^n \frac{r_i}{s_i} g \gamma_i \in \sum_{i=2}^n
%\frac{r_i}{s_i}\pi_v^{-k_i} e_i\;+ (\pi_v^{\; \ell-\alpha} \OOO_v)^n\,.
%\]
%Since $\alpha\geq k_1$ and $nN\geq 0$, we have $(\pi_v^{\ell-k_1}
%\OOO_v) \times (\pi_v^{\ell + n N-k_1} \OOO_v)^{n-1}\subset (\pi_v^{\;
%  \ell-\alpha} \OOO_v)^n$. Hence it follows from Equation
%\eqref{eq:system_count_latpt_dynball_1} that
%\[
%g \ga_1 \in \pi_v^{-k_1}
%e_1+\sum_{i=2}^n \frac{r_i}{s_i}\pi_v^{-k_i} e_i\;+ (\pi_v^{\;
%  \ell-\alpha} \OOO_v)^n\,.
%\]
%As above, this equation with unknown $\gamma_1$ has at most
%$q_v^{\;\alpha n}$ solutions, hence the set of solutions of Equation
%\eqref{eq:system_count_latpt_dynball_1} with unknown $\gamma_1$ is
%contained in a set (that depends on $\mb{t} \in
%\widetilde{\Lambda_{\mb{s}}}$) of cardinality at most $q_v^{\alpha n}$.

Let us fix a solution $(\ga_2,\ldots,\ga_n)$ of the system of
equations \eqref{eq:system_count_latpt_dynball_2}.  Given an element
$\mb{t} = (\frac{r_2}{s_2}, \ldots, \frac{r_n}{s_n})
\in\widetilde{\Lambda_{\mb{s}}}$, let us fix a solution
$\ga=(\ga_1\;\ga_2\;\dots\;\ga_n)\in\Ga$ of Equation
\eqref{eq:system_count_latpt_dynball_1} with prescribed last $n-1$
columns $\ga_2,\ldots,\ga_n$. Note that since $\ga\in\Ga=\SL_n(R_v)$,
the vectors of $K_v^{\;n}$ with column matrices $\ga_1,\ldots,\ga_n$
form an $R_v$-basis of $R_v^{\;n}$. Hence any other solution $\ga'=
(\ga'_1\;\ga_2\;\dots\;\ga_n)\in\Ga$ of this equation with these last
$n-1$ columns has a first column $\ga'_1$ such that there exists
$\lambda_1,\dots,\lambda_n\in R_v$ with $\ga'_1=\lambda_1\ga_1+
\lambda_2\ga_2 +\dots +\lambda_n\ga_n$.  Since the determinant of
$n$-tuples of elements of $K_v^{\;n}$ is multilinear and alternating,
and since $\det \ga=\det \ga'=1$, we have
\[
\lambda_1=\det(\lambda_1\ga_1,\ga_2,\dots,\ga_n)=
\det(\ga'_1,\ga_2,\dots,\ga_n)=1\,.
\]
Hence Equation \eqref{eq:system_count_latpt_dynball_1} for the matrix
$\ga'$ becomes 
\begin{equation}\label{eq:system_count_latpt_dynball_1prim}
g \ga_1 + \sum_{i=2}^n \big(\lambda_i-\frac{r_i}{s_i} \big)g \ga_i
\; \in\; \pi_v^{-k_1} e_1+ (\pi_v^{\ell-k_1} \OOO_v) \times
(\pi_v^{\ell + n N-k_1} \OOO_v)^{n-1}\,.
\end{equation}
%Let us fix $\ga \in \Ga$ whose colums $\ga_2,\ldots\ga_n$ verify
%Equation \eqref{eq:system_count_latpt_dynball_2}. Let us give an upper
%bound on the number of solutions $\mb{t}$ of Equation
%\eqref{eq:system_count_latpt_dynball_1} with $\ga_1$ thus fixed, when
%$\mb{t} = (\frac{r_2}{s_2}, \ldots, \frac{r_n}{s_n})$ varies in
%$\widetilde{\Lambda_{\mb{s}}}$.
We denote by $\operatorname{pr} : K_v^{\; n} \to K_v^{\; n-1}$ the
projection onto the last $n-1$ coordinates. Let $\ddd=\diag
(\pi_v^{\;k_2}, \ldots, \pi_v^{\; k_n})$. Let us multiply Equation
\eqref{eq:system_count_latpt_dynball_1prim} by the scalar $\mb{s}_*$
(defined above Equation \eqref{eq:propmbsstar}, so that $\frac{r_i
  \,\mb{s}_*}{s_i}\in R_v$ for every $i \in \llbracket 2,
n\rrbracket$). Let then project it to $K_v^{\; n-1}$ as well as
Equation \eqref{eq:system_count_latpt_dynball_2}. Let us then multiply
them by the matrix $\ddd=\diag(\pi_v^{\; k_2}, \ldots, \pi_v^{\;
  k_n})$ on the left.  The condition $g \gamma \uuu_{-\mb{t}} \in
B_{\ell,N} \exp(\mb{k})$ thus provides the system of conditions
\begin{align}
\sum_{i=2}^n \Big(\lambda_i\,\mb{s}_*- \frac{r_i \,\mb{s}_*}{s_i}\Big)
\ddd \operatorname{pr}(g\ga_i) &
\;\in\; \mb{s}_* \ddd \operatorname{pr}(g \ga_1) + \mb{s}_* \prod_{i=2}^n
(\pi_v^{\; \ell + nN + k_i - k_1} \OOO_v)\, ,
\label{eq:system_count_latpt_dynball_1_alt} \\
\forall\; i \in \llbracket 2, n\rrbracket, \quad
\ddd \operatorname{pr}(g \ga_i) &
\;\in\; e_i+ \ddd (\pi_v^{\; \ell-k_i} \OOO_v)^{n-1}
\label{eq:system_count_latpt_dynball_2_alt}\,.
\end{align}

We want to give an upper bound on the number of elements $\mb{t} \in
\wt{\Lambda_{\mb{s}}}$ such that there exists $\lambda_2,\dots,
\lambda_n\in R_v$ satisfying this system. By Equation
\eqref{eq:system_count_latpt_dynball_2_alt}, since $\ell>0$, there
exists a matrix $\wt g\in\M_{n-1}(K_v)$ with $\|\wt g\| \leq
q_v^{\;-\ell}$ such that the $R_v$-lattice $L=\oplus_{2\leq i\leq n}
R_v \ddd \operatorname{pr}(g \gamma_i)$ in $K_v^{\; n-1}$ is equal to
$(I_{n-1}+\ddd \,\wt{g}\,\ddd^{-1}) R_v^{\; n-1}$. By Lemma
\ref{lem:minosystole} applied with $d=n-1$, we have $\sys(L) \geq
1-q_v^{\,-\ell}\geq\frac{1}{2}$.  Note that for every
$i\in\intbra{2,n}$, the assumption on $N$ of Lemma
\ref{lem:counting_latpts_dynballs} gives the inequalities
$-v(\mb{s}_*)-\ell-nN-k_i+k_1 \geq 0$. Note that
$\mb{s}_*\in\pi^{v(\mb{s}_*)}\OOO_v^{\,\times}$.  Since each solution
$(\lambda_2\,\mb{s}_*+\frac{r_2\,\mb{s}_*}{s_2}, \ldots,
\lambda_n\,\mb{s}_*- \frac{r_n \,\mb{s}_*}{s_n}) \in R_v ^{\; n-1}$ of
Equation \eqref{eq:system_count_latpt_dynball_1_alt} corresponds to
one point of the lattice $L$, the associated number of solutions is
bounded from above by 
\begin{align*}
&\Big\lceil\frac{1}{\sys(L)^{n-1}}
  \vol_v^{n-1}\Big(\mb{s}_*\prod_{i=2}^n\pi_v^{\;\ell+nN+k_i-k_1}\OOO_v\Big)
  \Big\rceil\\\leq\; & \Big\lceil2^{n-1} \prod_{2 \leq i \leq n}
\vol_v\big(\pi_v^{\,v(\mb{s}_*)+\ell+nN+k_i-k_1}\OOO_v\big)\Big\rceil%\\&
=2^{n-1} \prod_{2 \leq i \leq n}q_v^{\,-v(\mb{s}_*)-\ell-nN-k_i+k_1}\\\leq\; &
2^{n-1}\,q_v^{\,-\ell(n-1)-v(\mb{s}_*)(n-1)-n(n-1)N+(n-1)k_1-\sum_{i=2}^n k_i}\,.
\end{align*}
Combining this with the previous counting results for the columns
$\ga_2,\ldots, \ga_n$, recalling that $\sum_{i=1}^n k_i =0$ (since
$\mb{k} \in \Delta_{\mb{s}}\subset \ZZ_0^n$), we finally obtain, as
wanted,
\begin{align*}
& \card\big(\{\exp(\mb{k}) \uuu_{\mb{t}} \Ga :
\mb{t}  \in \wt\Lambda_{\mb{s}}\} \cap (B_{\ell, N} \;x)\big)
\\ \leq \; & 2^{n-1}\,
q_v^{\,-\ell(n-1)-v(\mb{s}_*)(n-1)-n(n-1)N+(n-1)k_1-\sum_{i=2}^n k_i}
\,\prod_{i=2}^nq_v^{\;nk_i}
\\ \leq \; & 2^{n-1}\,q_v^{\,-\ell(n-1)-v(\mb{s}_*)(n-1)-n(n-1)N}\,.
\;\;\;\Box
\end{align*}
\medskip

Before stating the main result of Subsection \ref{subsec:dynpart}, let
us give some definitions.  For all $\ell,m\in\NN$, an {\it
  $(m,\ell)$-partition} of $\X$ is a finite measurable partition $\P=
\{P_1,P_2, \ldots, P_{|\P|} \}$ of $\X$ such that $P_1$ is equal to the
$q_v^{-m}$-thin part $\X^{<q_v^{-m}} = \{x\in\X : \sys(x)<q_v^{-m}\}$
of $\X$ (see Subsection \ref{subsec:systole}) and such that for every
$i \in \llbracket 2, |\P|\rrbracket$, there exists $x_i\in\X$ with
$P_i \subset B_\ell \,x_i$ with $B_\ell=B_{\ell,0}$ defined in
Equation \eqref{eq:defiBellN}. Note that for every
$(m,\ell)$-partition $\P$ of $\X$ and every $a\in A(\OOO_v)$, since
$\sys(ax)=\sys(x)$ for every $x\in \X$ and by Equation
\eqref{eq:invarWellNpara}, the partition $a^{-1}\P$ is also an
$(m,\ell)$-partition of $\X$.

For every $N \in\NN\ssm\{0\}$, the {\it $N$-th dynamical partition}
for $\aaa$ associated with a finite measurable partition $\P$ of $\X$
is the finite measurable partition
\begin{equation}\label{eq:defPN}
\P^N = \bigvee_{i=0}^{N-1}\aaa^{-i} \P\,.
\end{equation}
Note that for every $a\in A(\OOO_v)$, since $a$ commutes with $\aaa$,
we have $a^{-1}(\P^N)=(a^{-1}\P)^N$.  The {\it $N$-th Birkhoff
  average} for $\aaa$ of a Borel probability measure $\mu$ on $\X$ is
\begin{equation}\label{eq:defSNmu}
S_N \mu = \frac{1}{N} \sum_{i=0}^{N-1} \aaa^i_* \mu\,.
\end{equation}

The next lemma says that the thick part of the space $\X$ of special
unimodular $R_v$-lattices may be almost entirely covered by dynamical
balls $B_{\ell,N}\, x_j$ that are essentially finer than the partition
$\P^N$, with a good control of the cardinality of this cover.

\blemm\label{lem:cover_X_1_dynamic_balls}
%%
%\todo{\tiny Lem 3.7 \cite{DavSha20} nonarchimédien}
%%
For every $m\in\NN$, there exists $\ell_m\in\NN\ssm\{0\}$ such that
for every integer $\ell\geq\ell_m$, for every $(m,\ell)$-partition
$\P$ of $\X$, for every $\kappa \in \; ]0,1[\,$, and for every
$N\in\NN\ssm\{0\}$, the $q_v^{-m}$-thick part $\X^{\geq q_v^{-m}}$
of $\X$ contains a measurable subset $\X'=\X'_{\P,\,\kappa,\,N}$
satisfying the two following conditions.
\begin{enumerate}
\item\label{item1:cover_X_1_dynamic_balls} There exists a subset $\P'$
of $\P^N$ such that $\X' = \bigcup \P'$ and such that, for every $P
\in \P'$, there exists a finite subset $F_P$ of $P$ with cardinality
at most $q_v^{\; n^3\kappa N}$ such that $P\subset\bigcup_{x\in F_P}
B_{\ell,N-1}\, x$.
\item\label{item2:cover_X_1_dynamic_balls} For every Borel probability
measure $\mu$ on $\X$, we have $\mu(\X') \geq 1 - \frac{1}{\kappa}
\, S_N \mu(\X^{< q_v^{-m}})$.
\end{enumerate}
\elemm

\dem Let $m$, $\kappa$, $N$ and $\mu$ be as in the statement. Since
the action by left translations of $G$ on $G/\Ga$ is locally free and
since $\X^{\geq q_v^{-m}}$ is compact, there exists $\ell_m \in\NN
\ssm \{0\}$ such that for every $x \in \X^{\geq q_v^{-m}}$, the map $g
\mapsto gx$ is injective on the dynamical ball $B_{\ell_m-n}$. We may
assume that $\ell_m\geq n$ for future use.
%%
%\todo{Trouver un $\ell_m$ explicite ?}
%%
Let $\ell \geq \ell_m$ and let $\P = \{P_1, \ldots, P_{|\P|}\}$ be an
$(m,\ell)$-partition of $\X$ so that $P_1=\X^{<q_v^{-m}}$ and for
every $k \in \llbracket2, |\P|\rrbracket$, there exists $x_k \in \X$
such that $P_k \subset B_\ell\, x_k$. We define a function $f_N:\X\ra
[0,+\infty[$ counting in average the excursions before time $N$ of the
    diagonal orbits under $\aaa$ into the $q_v^{-m}$-thin part of $\X$
    by
\[
f_N:x\mapsto \frac{1}{N} \;\sum_{j=0}^{N-1} \mathbbm{1}_{\X^{< q_v^{-m}}}
(\aaa^j x)\,.
\]
We define $\X' = \{x\in\X : f_N(x) \leq \kappa \}$. By Markov's
inequality applied to the nonnegative random variable $f_N$, we have
\begin{align*}
1-\mu(\X')&=\mu\big(\{x\in\X : f_N(x) > \kappa \}\big)
\leq \frac{1}{\kappa}\int_\X f_N\;d\mu=
\frac{1}{\kappa N}\sum_{j=0}^{N-1}\int_\X \mathbbm{1}_{\X^{< q_v^{-m}}}
(\aaa^j x)\;d\mu(x)\\&=\frac{1}{\kappa N}\sum_{j=0}^{N-1}\int_\X
\mathbbm{1}_{\X^{< q_v^{-m}}}\;d({\aaa^j}_*\mu)=\frac{1}{\kappa}S_N\mu
(\mathbbm{1}_{\X^{< q_v^{-m}}})\,.
\end{align*}
Hence the set $\X'$ satisfies Assertion
\eqref{item2:cover_X_1_dynamic_balls}.

In order to prove Assertion \eqref{item1:cover_X_1_dynamic_balls},
first notice that $f_N$ can be described on the partition $\P^N$ as
follows. For all $P = \bigcap_{j=0}^{N-1} \aaa^{-j} P_{k_j} \in \P^N$,
where $k_0, \ldots, k_{N-1} \in \llbracket 1, |\P|\rrbracket$, since $P_1
= \X^{< q_v^{-m}}$, we have, for all $x \in P$,
\begin{equation}\label{eq:caracfN}
f_N(x)=\frac{1}{N}\;\card\{j \in \llbracket 0, N-1 \rrbracket :
k_j = 1\}\,.
\end{equation}
In particular, $f_N$ is constant on every $P \in \P^N$. The definition
of $\X'$ then implies that there exists a subset $\P'$ of the
partition $\P^N$ such that $\X' = \bigcup \P'$.

Now let $P \in \P'$, so that we have ${f_N}_{\mid P}\leq \kappa<1$. By
the definition of the partition $\P^N$ and since $\P$ is an
$(m,\ell)$-partition of $\X$, for every $j \in \llbracket 0, N-1
\rrbracket$, we have either $\aaa^j P \subset \X^{\geq q_v^{-m}}$ or
$\aaa^j P \subset \X^{< q_v^{-m}}$. The set $\{j \in \intbra{0, N-1} :
\aaa^j P \subset \X^{\geq q_v^{-m}} \}$ is nonempty, since otherwise
we would have $\aaa^j P \subset \X^{< q_v^{-m}}$ for every $j \in
\intbra{0, N-1}$, hence ${f_N}_{\mid P}=1$, a contradiction.  Let
$j_0\in\intbra{0, N-1}$ be the minimum of this nonempty subset of
$\NN$. By the definition of the partition $\P^N$ and since $\P$ is an
$(m,\ell)$-partition of $\X$, there exists $k_0 \in \llbracket 2, |\P|
\rrbracket$ such that $\aaa^{j_0} P \subset P_{k_0} \subset B_\ell
\,x_{k_0}$. This inclusion $\aaa^{j_0} P \subset B_\ell \,x_{k_0}$ is
the starting point in order to prove Assertion
\eqref{item1:cover_X_1_dynamic_balls} by using iterations of Lemma
\ref{lem:cover_dynamic_balls_SLn} and of Equation
\eqref{eq:dynamic_a_on_balls}. We define, for every $j \in \llbracket
0, N-1 \rrbracket$,
\[
V_j = \{i\in \llbracket 0, j \rrbracket  :  \aaa^i P \subset
\X^{<q_v^{-m}} \}\,,
\]
and we denote by $|V_j|$ its cardinality.

Let $C=q_v^{\; n^2}$ be the constant satisfying Lemma
\ref{lem:cover_dynamic_balls_SLn} for $\ell'=1$ (allowing
multiplicities), so that $C^n=q_v^{\; n^3}$ is the constant satisfying
Lemma \ref{lem:cover_dynamic_balls_SLn} for $\ell'=n$.

\medskip
\noindent{\bf Claim. } For every $j \in \llbracket 0, N-1 \rrbracket$
such that $|V_j| \neq j+1$, there exist $R_v$-lattices $y_{1,\,j},
\ldots, y_{C^{n|V_j|},\,j} \in P$ such that
\begin{equation}\label{fact:last_fact_lemma_cover_dynamic_balls}
P \subset \bigcup_{i=1}^{C^{n|V_j|}} B_{\ell,\,j} \,y_{i,\,j}\,.
\end{equation}

We have $|V_{N-1}| \leq \kappa N< N$ (hence $C^{n |V_{N-1}|} \leq
q_v^{\; n^3 \kappa N}$) since ${f_N}_{\mid P}=\frac{|V_{N-1}|}{N}$ by
Equation \eqref{eq:caracfN} and since $P \subset \X'$ so that
${f_N}_{\mid P}\leq \kappa$.  Therefore the case $j=N-1$ of this claim
implies Assertion \eqref{item1:cover_X_1_dynamic_balls}.

\medskip\noindent{\bf Proof of the claim. }  We proceed by induction
on $j\in \llbracket 0, N-1 \rrbracket$. By definition, we have $j_0 =
\min\{j\in\llbracket 0, N-1 \rrbracket : |V_j| \neq j+1 \}$, hence we
begin the induction at the step $j_0$, the previous cases being
empty. If $j_0=0$, we have $P \subset B_\ell\, x_{k_0}$ and by Lemma
\ref{lem:reciprocal_inclusion_dyn_balls} applied with $N=0$, we can
assume that $x_{k_0} \in P$, which proves the Claim at the $j_0$-th
step (since then $|V_0|=0$ or equivalently $C^{n|V_0|}=1$). If $j_0
\geq 1$, then we apply $n j_0$ times Lemma
\ref{lem:cover_dynamic_balls_SLn} with $N=0$, $\ell'=1$ and $S$
successively equal to $B_\ell, B_{\ell+1}, B_{\ell+2},
\ldots,B_{\ell+nj_0-1}$. This gives the existence of $g_{1}, \ldots,
g_{C^{n j_0}} \in G$ such that $B_\ell\subset \bigcup_{i=1}^{C^{n
    j_0}} B_{\ell+nj_0}\,g_i$. Hence by the inclusion $\aaa^{j_0} P
\subset B_\ell \,x_{k_0}$, we have $\aaa^{j_0} P \subset
\bigcup_{i=1}^{C^{n j_0}} B_{\ell+nj_0}\,g_ix_{k_0}$. Up to allowing
multiplicities, we may assume that for every $i\in\llbracket1,C^{n
  j_0}\rrbracket$, the intersection $(\aaa^{j_0} P)\cap
(B_{\ell+nj_0}\,g_ix_{k_0})$ is nonempty, hence contains an element
$\aaa^{j_0}y_{i,\,j_0}$ with $y_{i,\,j_0} \in P$. By Lemma
\ref{lem:reciprocal_inclusion_dyn_balls} and since $I_n\in
B_{\ell+nj_0}$, we have $B_{\ell+nj_0}\,g_ix_{k_0}=B_{\ell+nj_0}
\,\aaa^{j_0} y_{i,\,j_0}$. Therefore $\aaa^{j_0} P \subset
\bigcup_{i=1}^{C^{n j_0}} B_{\ell+nj_0} \,\aaa^{j_0} y_{i,\,j_0}$.
Hence by Equation \eqref{eq:dynamic_a_on_balls} applied with $N=0$,
$\ell'=-j_0$ and $\ell$ replaced by $\ell+nj_0$, we have
\[
P \subset\bigcup_{i=1}^{C^{nj_0}} \aaa^{-j_0} B_{\ell + n j_0} \,\aaa^{j_0}
y_{i, \,j_0} \subset \bigcup_{i=1}^{C^{n j_0}} B_{\ell, \,j_0} \,y_{i,\,j_0}\,.
\]
Since $|V_{j_0}|=j_0$, this proves the $j_0$-th step of the Claim. If
$j_0= N-1$, there is nothing more to be proved, hence we assume that
$j_0\leq N-2$.

Now let $j \in \llbracket j_0, N-2 \rrbracket$ and assume that the
$j$-th step of the Claim is satisfied, so that $P \subset
\bigcup_{i=1}^{C^{n|V_j|}} B_{\ell,\,j} \,y_{i,\,j}$ where $y_{1,\,j},
\ldots, y_{C^{n|V_j|},\,j} \in P$.

\begin{itemize}
\item First assume that $\aaa^{j+1} P \subset \X^{\geq q_v^{-m}}$ or
equivalently that $|V_{j+1}|=|V_j|$. Let us fix an element $i \in
\llbracket1, C^{n|V_j|}\rrbracket$ and let us prove that
$B_{\ell,\,j+1} \,y_{i,\,j} \cap P = B_{\ell,\,j} \,y_{i,\,j} \cap P$.
This will imply the $(j+1)$-th step of the Claim by setting
$y_{i,\,j+1}=y_{i,\,j}$. The inclusion $B_{\ell,\,j+1} \,y_{i,\,j}
\cap P \subset B_{\ell,\,j}\,y_{i,\,j} \cap P$ is clear by the
inclusion just above Equation \eqref{eq:voisdyncroiss}.  For the
converse one, let \mbox{$g \in B_{\ell,j}$} be such that $g y_{k,j}
\in P$. Since we have $\aaa^{j+1} P \subset \X^{\geq q_v^{-m}}$ and
since $\P$ is an $(m,\ell)$-partition of $\X$, there exists $k \in
\llbracket2, |\P|\rrbracket$ such that $\aaa^{j+1} P \subset P_k
\subset B_{\ell}\, x_k$. Let us define $\wt x= \aaa^{j+1}y_{i,\,j}$
and $\wt g=\aaa^{j+1} g\, \aaa^{-(j+1)}$. Since $y_{i,\,j} \in P$, we
have
\[
\wt x= \aaa^{j+1}y_{i,\,j}\in\aaa^{j+1}P\subset \X^{\geq
q_v^{-m}}\cap (B_{\ell}\, x_k)\,.
\]
Similarly, since $gy_{i,\,j} \in P$, we have $\wt g\,\wt x=
\aaa^{j+1}(gy_{k,\,j})\in B_{\ell} \,x_k$. Therefore we have $\wt
g\,\wt x\in B_{\ell} (B_{\ell})^{-1}\wt x\subset B_{\ell}\, \wt x$ by
Lemma \ref{lem:inverse_dyn_balls} and Equation
\eqref{eq:BellBelldansBell} (both with $N=0$ therein).  By Equations
\eqref{eq:dynamic_a_on_balls} and \eqref{eq:voisdyncroiss}, we have
$\wt g=\aaa^{j+1} g\, \aaa^{-(j+1)} \in B_{\ell, -1} \subset
B_{\ell-n} \subset B_{\ell_m -n}$.  We have $B_{\ell}\subset
B_{\ell_m}\subset B_{\ell_m -n}$ since $\ell\geq\ell_m\geq \ell_m-n$,
again by Equation \eqref{eq:voisdyncroiss}. Since $\wt g\,\wt x\in
B_{\ell}\, \wt x$, since $\wt x\in\X^{\geq q_v^{-m}}$ and by the
definition of $\ell_m$, we have $\wt g \in B_{\ell}$. Therefore, by
Equation \eqref{eq:dynamic_a_on_balls} again, we finally obtain
\[
g =\aaa^{-(j+1)} \wt g \,\aaa^{j+1}\in \aaa^{-(j+1)} B_\ell \,\aaa^{j+1}
\cap B_{\ell,\,j} \subset B_{\ell-n(j+1),\,j+1} \cap B_{\ell,\,j}
\subset B_{\ell,\,j+1}\,,
\]
so that $gy_{k,\,j}\in B_{\ell,\, j+1}\,y_{k,\,j}\cap P$, thus proving
the wanted converse inclusion.
%$B_{\ell,j+1} \,y_{k,j} \cap P \supset B_{\ell,j}\,y_{k,j} \cap P$
\item
Now assume that $\aaa^{j+1} P \subset \X^{< q_v^{-m}}$ or equivalently
that $|V_{j+1}|=|V_j|+1$. The proof of the $(j+1)$-th step of the
Claim is then straightforward by applying Lemma
\ref{lem:cover_dynamic_balls_SLn} with $\ell'=n$ (so that $\ell\geq
\ell_m\geq n\geq\ell'$) and $N=j$ to $S=B_{\ell,j}$ in order to cover
each $B_{\ell,\,j} \,y_{i,\,j} \cap P$ for $i \in \llbracket1,
C^{n|V_j|}\rrbracket$ by $C^n$ subsets of $\X$ of the form
\[
B_{\ell+n,\,j} \,y^{(i)}_{i',\,j+1} \cap P \subset
B_{\ell,\,j+1} \,y^{(i)}_{i',\,j+1} \cap P\,,
\]
where $i'\in \intbra{1,C^n}$ and $y^{(i)}_{i',\,j+1}\in\X$, thus
covering $P$ by $C^{n|V_j|}C^n=C^{n|V_{j+1}|}$ subsets $B_{\ell,\,j+1}
\,y^{(i)}_{i',\,j+1}$. As in the $j_0$-th step, by Lemma
\ref{lem:reciprocal_inclusion_dyn_balls}, we may assume that
$y^{(i)}_{i',\,j+1}\in P$.  \cqfd
\end{itemize}

\section{Non-escape of mass in the thin part}
\label{sec:noescapmass}

In this section, according to the first step of the program announced
after the statement of Theorem \ref{theo:EL}, we provide the material
that will be used in Section \ref{sec:lowboundentrop} in order to
prove that every weak-star accumulation point $\mu$ of the measures
$\nu_{\mb{s}}^\loz$ as $\mb{s}$ tends appropriately to $+\infty$ is a
probability measure on $\X$.
%
%\textcolor{blue}{{\bf Idée} : Pour $J$ idéal donné, compter les idéaux
%  $I$ tels que $\redn(I)\leq \epsilon\redn(J)$. Les principaux
%  sont en gros en proportion $\frac{1}{h_K}$, et il suffit de diviser par
%  $|R_v^\times|$ pour repasser aux éléments de $R_v$ engendrant les
%  idéaux principaux. Pour obtenir ceux premiers avec $J$, utiliser un
%  crible, ce qui revient parfois par diviser par la fonction
%  arithmétique adaptée évaluée en $J$. \`A faire plus tard, quand le
% reste aura avancé}

For every fixed $J \in \I_v^+$, we first estimate the number of
nonzero ideals that are coprime to $J$ and whose norm is comparatively
small with respect to the one of $J$. Recall (see Subsection
\ref{subsec:functionfield}) that $\varpi_v(J)$ is the number of prime
factors of $J$. For every $\epsilon>0$, let
\[
E_{J,\,{\rm prim}}(\epsilon)=
\big\{I\in \I_v^+ : (I,J) = 1,\;\redn(I)\leq \epsilon\redn(J)\big\}\,.
\]

\blemm\label{lem:cardEJprim}
There exists $c_1\geq 0$ such that for all $J \in \I_v^+$ and
$\epsilon\in q^{\ZZ} \,\cap \,]\,0,1[\,$, we have
%%
%\todo{\tiny Lem 3.16 \cite{DavSha20}}
%%
\[
\Big| \; \card (E_{J,\,{\rm prim}}(\epsilon)) - \epsilon\;
\frac{h_K\;q^{2-\ggg}\,(q_v-1)}{(q-1)^2\,q_v}\; \varphi_v(J) \;\Big|
\leq c_1\;2^{\,\varpi_v(J)}\,.
\]
\elemm

\dem Let $c_K=\frac{h_K\; q^{2-\ggg} \,(q_v-1)}{(q-1)^2\,q_v}>0$. Let
$E_J(\epsilon)= \big\{I\in \I_v^+ : \redn(I)\leq \epsilon\redn(J)
\big\}$.  By a standard sieving argument, with $\mu_v$ the Möbius
function defined in Subsection \ref{subsec:functionfield}, by Lemma
\ref{lem:gauscounting} since $\epsilon\in q^\ZZ$ and $N(I')\in q^\NN$
for every $I'\in \I_v^+$, and by Equation \eqref{eq:defiEulerfunct},
we have
\begin{align*}
  \card (E_{J,\,{\rm prim}}(\epsilon)) &= \sum_{I\in \I_v^+,\;I\,\mid\, J}
  \mu_v(I) \;\card (E_{JI^{-1}}(\epsilon))\\ &= \sum_{I\in \I_v^+,\;I\,\mid\, J}
  \mu_v(I)\big(c_K\;\epsilon\,N(JI^{-1})+\bigO(1)\big) \\ &
  = \epsilon\;c_K\,N(J)\sum_{I\in \I_v^+,\;I\,\mid\, J}
  \frac{\mu_v(I)}{N(I)}+\bigO\Big(\sum_{I\in \I_v^+,\;I\,\mid\, J}|\mu_v(I)|
  \;\Big)\\ & =  \epsilon\;c_K\;\varphi_v(J) +\bigO(2^{\,\varpi_v(J)})\,.
\end{align*}
This proves the result. \cqfd

\blemm \label{lem:decaythin} Assume that $R_v$ is principal.
%%
%\todo{\tiny Snif}
%%
There exists a constant $c_2\geq 1$ such that for all $\mb{s}= (s_2,
\ldots,s_n) \in (R_v\ssm \{0\})^{n-1} $ and $\mb{k} =(k_1,\ldots,k_n)
\in \Delta_{\mb{s}}$ with
\begin{equation}\label{eq:hyplemfuitemasse}
  \forall\;i\in\llbracket2,n \rrbracket,\quad
  2^{\varpi_v(s_i)} q_v^{\;k_i-k_1} \leq \frac{|s_i|}{\max\{1,\ln\ln |s_i|\}}\,,
\end{equation}
and for every $\epsilon\in q_v^{\,\ZZ}\cap\,]0,1[\,$, we have
%%
%\todo{\tiny analogue Lem 3.17 \cite{DavSha20}}
%%
\[\Big(\frac{1}{\card\; \Lambda_{\mb{s}}}\sum_{\mb{t}\in\Lambda_{\mb{s}}}
\delta_{\exp(\mb{k})x_{\mb{t}}}\Big) (\X^{\leq \epsilon}) \leq
c_2\;\epsilon^{n}\,.
\]
\elemm

\dem Let $\mb{s}$, $\mb{k}$, $\epsilon$ be fixed as in the statement.
Let $\mb{t} = (\frac{r_2}{s_2}, \ldots, \frac{r_n}{s_n}) \mod
R_v^{\,n-1}$ that varies in $\Lambda_{\mb{s}}$. Recall that
$x_{\mb{t}}=\uuu_{\mb{t}}R_v^{\,n}$. By the definition in Subsection
\ref{subsec:systole} of the $\epsilon$-thin part $\X^{\leq \epsilon}$
of $\X$, we have $\exp(\mb{k})\,x_{\mb{t}} \in \X^{\leq \epsilon}$ if
and only if there exists a nonzero element $\lambda \in R_v^{\,n}$
such that $\|\exp(\mb{k})\, \uuu_{\mb{t}}\,\lambda\| \leq \epsilon$,
or equivalently by an easy computation if and only if the following
joint system of inequalities with unknown $(\lambda_1,
\lambda_2,\ldots, \lambda_n)$ in $R_v^{\,n}$ has a nonzero solution
\begin{eqnarray}
\label{eq:lemfuitemasse1}
|\lambda_1| & \leq & \epsilon\, q_v^{\,-k_1}  \\  
\label{eq:lemfuitemasse2n}
\forall \;i\in\llbracket 2, n\rrbracket,\quad
\Big|\lambda_1 \frac{r_i}{s_i} +\lambda_i\Big| &\leq&
\epsilon\, q_v^{\,-k_i}\,.
\end{eqnarray}

Note that if $(\lambda_1,\ldots,\lambda_n)\in R_v^{\,n}$ is a nonzero
solution to the joint system \eqref{eq:lemfuitemasse1} and
\eqref{eq:lemfuitemasse2n}, then $\lambda_1\neq 0$. Indeed, by
Equation \eqref{eq:inversiRv}, the only element of $R_v$ contained in
the closed ball $B(0,q_v^{\,-1})$ of center $0$ and radius
$q_v^{\,-1}$ is $0$. Hence if $\lambda_1= 0$, then for every $i\in
\llbracket 2, n\rrbracket$, since $\epsilon<1$ and $k_i\geq 0$ as
$\mb{k} \in \Delta_{\mb{s}}$, we have $\lambda_i\in R_v\cap
B(0,q_v^{\,-k_i}\epsilon)=\{0\}$, which contradicts the fact that
$(\lambda_1,\ldots,\lambda_n)\neq 0$.

By the ultrametric triangle inequality, if $\lambda, \lambda'$ are
distinct elements of $R_v$, then the closed balls $B(\lambda,
q_v^{\,-1})$ and $B(\lambda',q_v^{\,-1})$ are disjoint.  Again by the
ultrametric triangle inequality, for every $\rho\geq q_v^{\,-1}$, the
closed ball $B(0,\rho)$ contains $\bigcup_{\lambda\in R_v\cap
  B(0,\,\rho)} B(\lambda, q_v^{\,-1})$, and this union is a disjoint
union.  Recall that $\vol_v(B(0,q_v^{\,-1}))= q_v^{\,-1}$ by the
normalisation of the Haar measure $\vol_v$ of $K_v$. Separating the
cases when $\epsilon\,q_v^{\,-k_1}<q_v^{\,-1}$ or the contrary, and
since $\epsilon\in q_v^{\,\ZZ}$, the number of nonzero solutions
$\lambda_1$ of Equation \eqref{eq:lemfuitemasse1} hence satisfies
\begin{equation}\label{eq:controllambdaone}
\card((R_v\ssm\{0\}) \cap B(0,\epsilon\,q_v^{\,-k_1}))\leq
\frac{\vol_v(B(0,\epsilon\,q_v^{\,-k_1}))}{\vol_v(B(0,q_v^{\,-1}))}
= \epsilon\,q_v^{\,-k_1+1}\,.
\end{equation}

Let us now fix $\lambda_1\in(R_v\ssm\{0\})\cap B(0,q_v^{-k_1}
\epsilon)$ and $i\in\llbracket 2,n\rrbracket$. In this proof, the
$\bigO(\;)$ functions do not depend on $\lambda_1$, $\mb{s}$,
$\mb{k}$, $\epsilon$.  Let
\[
\N_i(\lambda_1)=\card\Big\{\frac{r_i}{s_i}\!\!\mod R_v: (r_i,s_i)=1,\;
\exists\;\lambda_i\in R_v,\;
\big|\lambda_1\frac{r_i}{s_i}+\lambda_i\big|
\leq \epsilon\, q_v^{\,-k_i}\Big\}\,.
\]

\noindent{\bf Claim 1 : } We have $\N_i(\lambda_1)=
\bigO\big(\epsilon\,q_v^{\,-k_i}\,\varphi_v(s_i)\big)$.

\medskip
By the discussion above Equations \eqref{eq:lemfuitemasse1},
\eqref{eq:lemfuitemasse2n}, by Equation \eqref{eq:controllambdaone}
and this claim, since $k_1=-\sum_{i=2}^n k_i$ as $\mb{k}
\in\Delta_{\mb{s}}\subset \ZZ^n_0$, and by Equation
\eqref{eq:cardLambdasss} on the left, this will imply the inequality
\begin{align*}
  \card\big\{\mb{t}\in\Lambda_{\mb{s}}:
  \exp(\mb{k})x_{\mb{t}}\in \X^{\leq \epsilon}\big\} &\leq
  \sum_{\lambda_1\in(R_v\,\ssm\,\{0\})\cap B(0,\,q_v^{-k_1}\epsilon)}
  \prod_{i=2}^n\N_i(\lambda_1)\\&=\bigO\Big(\epsilon\,q_v^{\,-k_1+1}
  \prod_{i=2}^n\,\epsilon\,q_v^{\,-k_i}\,\varphi_v(s_i)\Big)\\&=
  \bigO\big(\epsilon^{n}\,\card\; \Lambda_{\mb{s}}\big)\,.
\end{align*}
This estimate will prove Lemma \ref{lem:decaythin}.

\medskip
\noindent{\bf Proof of Claim 1. } Let $J_i=s_iR_v \in\I_v^+$. Since
$R_v$ is assumed to be principal, let $d_i\in R_v $ be such that
$\lambda_1R_v+J_i=d_iR_v$. Let $\wt \lambda_1= \frac{\lambda_1}{d_i}$,
$\wt s_i=\frac{s_i}{d_i}$ and $\wt J_i=\wt s_iR_v\in\I_v^+$. By
dividing by $d_i$ and since the fibers of the canonical morphism
$(R_v/J_i)^\times \ra (R_v/\,\wt J_i)^\times$ have order
$\frac{\varphi_v(J_i)}{\varphi_v(\wt J_i)}$, we have
\begin{align*}
 \N_i(\lambda_1)&=\card\,\big\{r_i+J_i\in (R_v/J_i)^\times:
\exists\;\lambda_i\in R_v,\;\big|\,\lambda_1 \,r_i+
\lambda_i\;s_i\big| \leq \epsilon\,q_v^{\,-k_i}\,|s_i|\big\}
\\&=\frac{\varphi_v(J_i)}{\varphi_v(\wt J_i)}\;
\card\,\big\{r_i+\wt J_i\in (R_v/\wt J_i)^\times:
\exists\;\lambda_i\in R_v,\;\big|\,\wt \lambda_1 \,r_i+
\lambda_i\;\wt s_i\big| \leq \epsilon\,q_v^{\,-k_i}\,|\wt s_i|\big\}\,.
\end{align*}
Since $\wt \lambda_1$ and $\wt s_i$ are coprime, let $\wt\lambda_1^{-}
\in R_v$ be such that $\lambda_1\,\wt\lambda_1^{-}-1\in \wt J_i$. The
multiplication by $\wt\lambda_1^{-}$ is a bijective map from $(R_v/\wt
J_i)^\times$ to itself. Hence by Lemma \ref{lem:cardEJprim}, we have
\begin{align*}
\N_i(\lambda_1)&=\frac{\varphi_v(J_i)}{\varphi_v(\wt J_i)}\;
\card\,\big\{r_i+\wt J_i\in (R_v/\wt J_i)^\times:
\exists\;\lambda_i\in R_v,\;\big|\,r_i+
\lambda_i\;\wt s_i\big| \leq \epsilon\,q_v^{\,-k_i}\,|\wt s_i|\big\}\,.
\\&\leq\frac{\varphi_v(J_i)}{\varphi_v(\wt J_i)}\;
\card\,\big\{I\in\I_v^+: (I,\wt J_i)=1,\; \redn(I) \leq
\epsilon\,q_v^{\,-k_i}\,\redn(\wt J_i)\big\}
\\&=\frac{\varphi_v(J_i)}{\varphi_v(\wt J_i)}\;\bigO
\big(\epsilon\,q_v^{\,-k_i}\,\varphi_v(\wt J_i)+2^{\,\varpi_v(\wt J_i)}\big)
=\bigO\Big(\epsilon\,q_v^{\,-k_i}\,\varphi_v(J_i)+
\frac{2^{\,\varpi_v(\wt J_i)}}{\varphi_v(\wt J_i)}\,\varphi_v(J_i)\Big)\,.
\end{align*}

\noindent{\bf Claim 2 : }
We have $\frac{2^{\,\varpi_v(\wt J_i)}}{\varphi_v(\wt J_i)}=
\bigO(\epsilon\,q_v^{\,-k_i})$.

\medskip
With the previous formula, this implies Claim 1, hence concludes Lemma
\ref{lem:decaythin}

\medskip
\noindent{\bf Proof of Claim 2. } By Lemma \ref{lem:PoonenRosen},
since $\varpi_v(\wt J_i)\leq \varpi_v(J_i)$ as $\wt J_i$ divides
$J_i$, since $\redn(\wt J_i)\leq\redn(J_i)=\redn(\wt J_i)\,
\redn(\lambda_1R_v+\wt J_i)$, and since $\redn(\lambda_1R_v+\wt J_i)
\leq \redn(\lambda_1R_v)=|\lambda_1|\leq \epsilon\,q_v^{-k_1}$ by
Equation \eqref{eq:lemfuitemasse1}, we have
\begin{align*}
\frac{2^{\,\varpi_v(\wt J_i)}}{\varphi_v(\wt J_i)}&=\bigO\Big(
2^{\,\varpi_v(\wt J_i)}\frac{\ln\ln\redn(\wt J_i)}{\redn(\wt J_i)}\Big)
=\bigO\Big(2^{\,\varpi_v(J_i)}\frac{\ln\ln\redn(J_i)}{\redn(J_i)}
\redn(\lambda_1R_v+\wt J_i)\Big)
\\&=\bigO\Big(\epsilon\,q_v^{-k_1}\,2^{\,\varpi_v(J_i)}
\frac{\ln\ln\redn(J_i)}{\redn(J_i)}\Big)\,.
\end{align*}
Claim 2 hence follows by the technical Assumption
\eqref{eq:hyplemfuitemasse} of Lemma
\ref{lem:decaythin}. \hfill$\Box$ $\Box$ $\Box$

\medskip
Let $\mb{s}= (s_2, \ldots,s_n) \in (R_v\ssm \{0\})^{n-1}$ satisfying
Equation \eqref{eq:proprisbold}, so that there exists a permutation
$\sigma$ of $\llbracket 2,n \rrbracket$ with $s_{\sigma^{-1}(2)}\mid
s_{\sigma^{-1}(3)}\ldots\mid \mb{s}_*= s_{\sigma^{-1}(n)}$ and
$v(\mb{s}_*)\in n\,\ZZ$.  Let
\[
\mb{w}=(1-n,1,\dots,1)\in\ZZ_0^n, \quad\text{so
that}\quad
\mb{k}_{\mb{s}}=
-\frac{v(\mb{s}_*)}{n}\;\mb{w}\quad\text{and}\quad
\aaa
%=\begin{pmatrix} \pi_v^{n-1} &
%0\\0 & \pi_v^{\;-1}I_{n-1} \end{pmatrix}
=\exp(\mb{w})\,.
\]
For all $\mb{k}\in\loz_{\mb{s}}$ and $N\in\NN\ssm\{0\}$, we denote by
$[\mb{k},N]$ (see picture below) the discrete interval in $\ZZ_0^n$
defined by
\[
[\mb{k},N]=\{\mb{k}+\ell'\mb{w}:\ell'\in\llbracket 0,N-1\rrbracket\}\,.
\]
By the convexity of $\loz_{\mb{s}}$, we have $[\mb{k},N]\subset
\loz_{\mb{s}}$ if and only if $\mb{k}+(N-1)\mb{w}\in\loz_{\mb{s}}$. Let
\begin{equation}\label{eq:definuskN}
\nu_{\mb{s},[\mb{k},N]}=\frac{1}{N\;\card\; \Lambda_{\mb{s}}}
\sum_{\mb{t}\in\Lambda_{\mb{s}},\;\ell'\in\llbracket 0,\,N-1\rrbracket}
\delta_{\exp({\mb{k}+\ell'\mb{w}})x_{\mb{t}}}\,,
\end{equation}
which is a probability measure on $\X$.
\begin{center}
  \!\!\!\!\!\!\!\!\!\!\!\!\!\!\!\!\!\!\!\!\!\!\!\!\!\!\!
  \begin{picture}(0,0)%
\includegraphics{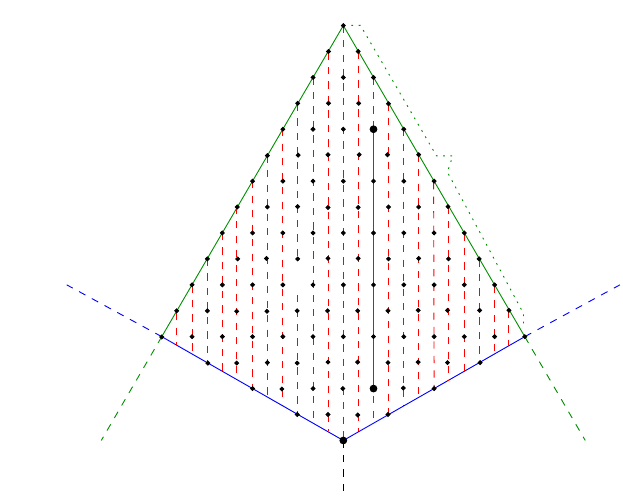}%
\end{picture}%
\setlength{\unitlength}{3812sp}%
\begingroup\makeatletter\ifx\SetFigFont\undefined%
\gdef\SetFigFont#1#2#3#4#5{%
  \reset@font\fontsize{#1}{#2pt}%
  \fontfamily{#3}\fontseries{#4}\fontshape{#5}%
  \selectfont}%
\fi\endgroup%
\begin{picture}(5146,4129)(-1094,-1337)
\put(1888,763){\makebox(0,0)[lb]{\smash{{\SetFigFont{9}{10.8}{\rmdefault}{\mddefault}{\updefault}{\color[rgb]{0,0,0}$\mb{k}+2\mb{w}$}%
}}}}
\put(1884,322){\makebox(0,0)[lb]{\smash{{\SetFigFont{9}{10.8}{\rmdefault}{\mddefault}{\updefault}{\color[rgb]{0,0,0}$\mb{k}+3\mb{w}$}%
}}}}
\put(1884,1622){\makebox(0,0)[lb]{\smash{{\SetFigFont{9}{10.8}{\rmdefault}{\mddefault}{\updefault}{\color[rgb]{0,0,0}$\mb{k}$}%
}}}}
\put(1888,1201){\makebox(0,0)[lb]{\smash{{\SetFigFont{9}{10.8}{\rmdefault}{\mddefault}{\updefault}{\color[rgb]{0,0,0}$\mb{k}+\mb{w}$}%
}}}}
\put(1787,-543){\makebox(0,0)[lb]{\smash{{\SetFigFont{9}{10.8}{\rmdefault}{\mddefault}{\updefault}{\color[rgb]{0,0,0}$\mb{k}+(N-1)\mb{w}$}%
}}}}
\put(-184,-844){\makebox(0,0)[lb]{\smash{{\SetFigFont{9}{10.8}{\rmdefault}{\mddefault}{\updefault}{\color[rgb]{0,.69,0}$k_2=0$}%
}}}}
\put(3789,-844){\makebox(0,0)[lb]{\smash{{\SetFigFont{9}{10.8}{\rmdefault}{\mddefault}{\updefault}{\color[rgb]{0,.69,0}$k_3=0$}%
}}}}
\put(-793,481){\makebox(0,0)[lb]{\smash{{\SetFigFont{9}{10.8}{\rmdefault}{\mddefault}{\updefault}{\color[rgb]{0,0,1}$k_3=k_1-v(\mb{s_*})$}%
}}}}
\put(3467,445){\makebox(0,0)[lb]{\smash{{\SetFigFont{9}{10.8}{\rmdefault}{\mddefault}{\updefault}{\color[rgb]{0,0,1}$k_2=k_1-v(\mb{s_*})$}%
}}}}
\put(2701,1514){\makebox(0,0)[lb]{\smash{{\SetFigFont{11}{13.2}{\rmdefault}{\mddefault}{\updefault}{\color[rgb]{0,.56,0}$\partial^+\loz_{\mb{s}}$}%
}}}}
\put(1263,494){\makebox(0,0)[lb]{\smash{{\SetFigFont{11}{13.2}{\rmdefault}{\mddefault}{\updefault}{\color[rgb]{1,0,0}{\large $\loz_{\mb{s}}$}}%
}}}}
\put(1785,-987){\makebox(0,0)[lb]{\smash{{\SetFigFont{9}{10.8}{\rmdefault}{\mddefault}{\updefault}{\color[rgb]{0,0,0}$k_{\mb{s}}$}%
}}}}
\put(3324,-56){\makebox(0,0)[lb]{\smash{{\SetFigFont{9}{10.8}{\rmdefault}{\mddefault}{\updefault}{\color[rgb]{0,0,0}$(\lfloor\frac{v(s_3)}{2}\rfloor,-\lfloor\frac{v(s_3)}{2}\rfloor,0)$}%
}}}}
\put(-1079,-56){\makebox(0,0)[lb]{\smash{{\SetFigFont{9}{10.8}{\rmdefault}{\mddefault}{\updefault}{\color[rgb]{0,0,0}$(\lfloor\frac{v(s_3)}{2}\rfloor,0,-\lfloor\frac{v(s_3)}{2}\rfloor)$}%
}}}}
\put(1781,2645){\makebox(0,0)[lb]{\smash{{\SetFigFont{9}{10.8}{\rmdefault}{\mddefault}{\updefault}{\color[rgb]{0,0,0}$(0,0,0)$}%
}}}}
\put(1001,2378){\makebox(0,0)[lb]{\smash{{\SetFigFont{9}{10.8}{\rmdefault}{\mddefault}{\updefault}{\color[rgb]{0,0,0}$(-1,0,1)$}%
}}}}
\put(1916,2374){\makebox(0,0)[lb]{\smash{{\SetFigFont{9}{10.8}{\rmdefault}{\mddefault}{\updefault}{\color[rgb]{0,0,0}$(-1,1,0)$}%
}}}}
\put(854,2159){\makebox(0,0)[lb]{\smash{{\SetFigFont{9}{10.8}{\rmdefault}{\mddefault}{\updefault}{\color[rgb]{0,0,0}$(-2,0,2)$}%
}}}}
\put(2059,2159){\makebox(0,0)[lb]{\smash{{\SetFigFont{9}{10.8}{\rmdefault}{\mddefault}{\updefault}{\color[rgb]{0,0,0}$(-2,2,0)$}%
}}}}
\put(1195,-1277){\makebox(0,0)[lb]{\smash{{\SetFigFont{9}{10.8}{\rmdefault}{\mddefault}{\updefault}{\color[rgb]{0,0,0}$k_2=k_3$}%
}}}}
\end{picture}%

\end{center}
The next corollary proves the non-escape of mass at infinity property 
for averages of measures parametrized by the discrete interval
$[\mb{k},N]$. Let us recall the positive constant $c_{\varpi_v}$
introduced in Equation \eqref{eq:majoomegav}. For every $\mb{s}= (s_2,
\ldots, s_n) \in(R_v\ssm \{0\})^{n-1}$ satisfying Equation
\eqref{eq:proprisbold}, let
\begin{equation}\label{eq:defikappaprims}
\kappa'(\mb{s})=\frac{1}{n}\Big(-v(\mb{s}_*)+\max_{i\in\llbracket
  2,n\rrbracket}\log_{q_v} \frac{2^{\varpi_v(s_i)}\,\max\{1,\ln\ln
  |s_i|\}}{|s_i|}\;\Big)\,.
\end{equation}

\brema\label{rem:controlkappaprims} If there exists $c_0\geq 0$ such
that $\max_{i\in\llbracket 2,n \rrbracket}\,v(s_i)-v(\mb{s}_*)\leq
c_0\frac{-v(\mb{s}_*)}{\max\{1,\, \ln(-v(\mb{s}_*))\}}$, then
$\kappa'(\mb{s})\leq \frac{1}{n} (c_{\varpi_v}+c_0+1)
\frac{-v(\mb{s}_*)}{\max\{1,\,\ln(-v(\mb{s}_*))\}}$. In particular,
$\kappa'(\mb{s})$ is negligible with respect to $-v(\mb{s}_*)$ as
$-v(\mb{s}_*)\ra+\infty$.
\erema

\dem Since $|s_i|=q_v^{-v(s_i)}$, since $\frac{1}{2}\leq\ln 2\leq \ln
q_v$ and by Equation \eqref{eq:majoomegav}, since the maps $f_1: t
\mapsto \frac{t}{\max\{1,\,\ln t\}}$ and $f_2:t\mapsto 2\ln(\max\{1,
\ln t\})$ on $[0,+\infty[$ are nondecreasing with $f_2\leq f_1$, and
by the assumption of the remark, we have
\begin{align*}
  n\;\kappa'(\mb{s})&=\max_{i\in\llbracket
    2,n\rrbracket}\Big(\log_{q_v}
  2^{\varpi_v(s_i)}+\log_{q_v}(\max\{1,\ln(-v(s_i))\}) +
  v(s_i)-v(\mb{s}_*)\Big)\\&\leq\max_{i\in\llbracket
    2,n\rrbracket} \Big(c_{\varpi_v}\frac{-v(s_i)}{\max\{1,\ln(-v(s_i))\}}
  +2\ln(\max\{1,\ln(-v(s_i))\})+ v(s_i)-v(\mb{s}_*)\Big)
  \\&\leq (c_{\varpi_v}+1)\frac{-v(\mb{s}_*)}{\max\{1,\ln(-v(\mb{s}_*))\}}
  +c_0\frac{-v(\mb{s}_*)}{\max\{1,\,\ln(-v(\mb{s}_*))\}}\,,
\end{align*}
which proves the result. \cqfd 

\bcoro \label{coro:decaythin2} Assume that $R_v$ is principal. For all
$\mb{s} \in (R_v\ssm \{0\})^{n-1}$ satisfying Equation
\eqref{eq:proprisbold}, $\mb{k} \in \loz_{\mb{s}}$, $N\in\NN\ssm\{0\}$
and $\epsilon \in q_v^{\,\ZZ}\cap\,]0,1[\,$ such that
\begin{equation}\label{eq:hyplemfuitemasse2}
  \mb{k}+(N-1)\mb{w}\in\loz_{\mb{s}}\quad\text{and}\quad
  N\geq\frac{\kappa'(\mb{s})}{c_2\,\epsilon^n}\,,
  %\kappa'(\mb{s})\max\big\{1,\frac{1}{c_2\,\epsilon^n}\big\}\,,
\end{equation}
we have
%%
%\todo{\tiny Lem 3.18 \cite{DavSha20}}
%%
\[
\nu_{\mb{s},[\mb{k},N]} (\X^{\leq \epsilon}) \leq 2\,c_2\;\epsilon^{n}\,.
\]
\ecoro

\dem Let $\mb{s}= (s_2, \ldots,s_n)$, $\mb{k} =(k_1,\ldots,k_n)$, $N$,
$\epsilon$ be fixed as in the statement. Let $\ell'$ be an integer
that will vary in $\llbracket 0,N-1 \rrbracket$. First assume that $\ell'
\leq N-1 -\kappa'(\mb{s})$. Note that $\kappa'(\mb{s})\geq \frac{1}{n}
\log_{q_v}(2^{\varpi_v(\mb{s}_*)}\,\max\{1,\ln\ln|\mb{s}_*|\})\geq 0$.

By the definition of $\loz_{\mb{s}}$, since $\mb{k}+ (N-1)\mb{w}\in
\loz_{\mb{s}}$ by the left hand side of Assumption
\eqref{eq:hyplemfuitemasse2}, we have $\max_{2\leq i\leq n}
(\mb{k}+(N-1)\mb{w})_i \leq (\mb{k}+(N-1)\mb{w})_1-v(\mb{s}_*)$.
Hence for every $i\in \llbracket 2,n\rrbracket$, since $\ell'\leq N-1
-\kappa'(\mb{s})$ and by the definition of $\kappa'(\mb{s})$, we have
\begin{align*}
  (\mb{k}+\ell'\mb{w})_i-(\mb{k}+\ell'\mb{w})_1
  &=(k_i+\ell')-(k_1+(1-n)\ell')
  \\&=(\mb{k}+(N-1)\mb{w})_i-(\mb{k}+(N-1)\mb{w})_1 +n(\ell'-N+1)
  \\&\leq -\,v(\mb{s}_*) -n\,\kappa'(\mb{s})
  \leq \log_{q_v}\frac{|s_i|}{2^{\varpi_v(s_i)}\,\max\{1,\ln\ln |s_i|\}}\,.
\end{align*}
Therefore the element $\mb{k}+\ell'\mb{w}$ of $\loz_{\mb{s}}$
satisfies the technical Assumption \eqref{eq:hyplemfuitemasse} of Lemma
\ref{lem:decaythin}, and we have
\[\Big(\frac{1}{\card\; \Lambda_{\mb{s}}}\sum_{\mb{t}\in\Lambda_{\mb{s}}}
\delta_{\exp(\mb{k}+\ell'\mb{w})x_{\mb{t}}}\Big) (\X^{\leq \epsilon}) \leq
c_2\;\epsilon^{n}\,.
\]
There are at most $N$ (resp.~$\kappa' (\mb{s})$) integral elements in
the real interval $[0,N-1-\kappa' (\mb{s})]$ (resp.~$]N-1-
  \kappa'(\mb{s}), N-1]$). Therefore separating, in Equation
\eqref{eq:definuskN} that defines the measure $\nu_{\mb{s},[\mb{k},N]}$,
the summation over $\ell'\in \llbracket0, N-1\rrbracket$ in firstly
$\ell'\in[0,N -1-\kappa'(\mb{s})]$ and secondly $\ell'\in\; ]N-1-
\kappa'(\mb{s}),N-1]$, we have
\[
\nu_{\mb{s},[\mb{k},N]} (\X^{\leq \epsilon}) \leq
\frac{1}{N}(N\;c_2\;\epsilon^{n})+
\frac{\kappa'(\mb{s})}{N}\,.
\]
By the right hand side of Assumption \eqref{eq:hyplemfuitemasse2},
this proves Corollary \ref{coro:decaythin2}.
\cqfd

\section{Optimal entropy lower bound}
\label{sec:lowboundentrop}

In this final section, we prove the main equidistribution result of
this paper, in the space $\X=\SL_n(K_v)/\SL_n(R_v)$ of special
unimodular $R_v$-lattices in $K_v^{\;n}$ towards its homogeneous
measure ${\tt m}_{\X}$, of the measures supported on large subsets of
divergent orbits of type $(1, s_2 ,\ldots , s_n)$ (up to permutation)
as $\mb{s}=(s_2 , \ldots , s_n) \in (R_v\ssm\{0\})^{n-1}$ tends to
infinity (for the Fréchet filter or equivalently when
$\min_{i\in\llbracket 2,n \rrbracket} v(s_i)$ tends to $-\infty$). We
will actually require some uniform convergence to $-\infty$ of the
valuations of the components $s_2 ,\ldots , s_n$ of $\mb{s}$, and
precisely
\begin{equation}\label{eq:unifbehavsi}
  \exists\,c_0\geq 0,\qquad
  \max_{i,j\in\llbracket 2,n \rrbracket}\,|\,v(s_i)-v(s_j)\,|=
  \max_{i\in\llbracket 2,n \rrbracket}\,v(s_i)-v(\mb{s}_*)
\leq c_0\frac{-v(\mb{s}_*)}{\max\{1,\, \ln(-v(\mb{s}_*))\}}\,.
\end{equation}
Note that this is for instance satisfied if $s_2 =\ldots = s_n$ as in
Theorem \ref{theo:mainuintro2} in the Introduction, and that this
assumption is optimal by Remark \ref{rem:controlkappaprims}.

\btheo\label{theo:main} Assume that $R_v$ is principal. As $\mb{s}\in
(R_v\ssm \{0\})^{n-1}$ satisfying Equations \eqref{eq:proprisbold} and
\eqref{eq:unifbehavsi} tends to infinity, the measures
$\nu_{\mb{s}}^\loz$ weak-star converge to $\frac{{\tt m}_{\X}}{\|{\tt
    m}_{\X}\|}$ on $\X$.
%%
%\todo{\tiny analogue Theorem 1.2 (3) \cite{DavSha20}}
%%
\etheo

\dem Let us fix a weak-star accumulation point $\nu$ of the measures
$\nu_{\mb{s}}^\loz$ as $\mb{s}\in (R_v\ssm \{0\})^{n-1}$ satisfying
Equations \eqref{eq:proprisbold} and \eqref{eq:unifbehavsi} tends to
infinity. We will prove that $\nu$ is a probability measure using the
work of Section \ref{sec:noescapmass} and that $\nu= \frac{{\tt
    m}_{\X}}{\|{\tt m}_{\X}\|}$ using the entropy method described in
Section \ref{sec:entropytoequidistrib}, which will conclude using the
Banach-Alaoglu theorem.

\blemm\label{lem:nuaaainvar}
The measure $\nu$ is $\aaa$-invariant.
\elemm

\dem Recall that $\mb{w}=(1-n,1, \ldots,1)\in\ZZ^n_0$. Using the
definitions \eqref{eq:definuloz} and \eqref{eq:defaaa}, since
$\aaa=\exp(\mb{w})$ commutes with $A(\OOO_v)$, we have
\[
\aaa_* \nu_{\mb{s}}^\loz =
\frac{1}{\card\;\Lambda_{\mb{s}}\;\card\;\loz_{\mb{s}}}\;
\sum_{\mb{t}\in\Lambda_{\mb{s}},\;\mb{k}\in\loz_{\mb{s}}}\int_{a\in A(\OOO_v)}
\;\delta_{a\exp(\mb{k}+\mb{w})\,x_{\mb{t}}}\;da\,.
\]
In order to compare $\nu_{\mb{s}}^\loz$ and $\aaa_*\nu_{\mb{s}}^\loz$,
let us give an upper estimate on the cardinality of the symmetric
difference between $\loz_{\mb{s}}$ and $\loz_{\mb{s}}+\mb{w}$.
\begin{center}
  \begin{picture}(0,0)%
\includegraphics{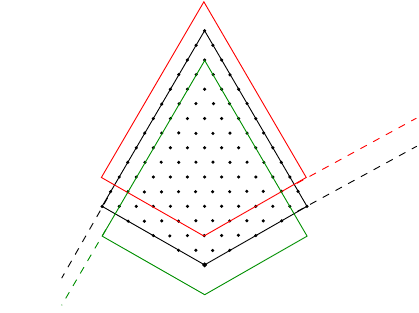}%
\end{picture}%
\setlength{\unitlength}{3812sp}%
\begingroup\makeatletter\ifx\SetFigFont\undefined%
\gdef\SetFigFont#1#2#3#4#5{%
  \reset@font\fontsize{#1}{#2pt}%
  \fontfamily{#3}\fontseries{#4}\fontshape{#5}%
  \selectfont}%
\fi\endgroup%
\begin{picture}(3464,2586)(616,-1876)
\put(1216,-1816){\makebox(0,0)[lb]{\smash{{\SetFigFont{11}{13.2}{\rmdefault}{\mddefault}{\updefault}{\color[rgb]{0,.69,0}$k_2=1$}%
}}}}
\put(2206,-1636){\makebox(0,0)[lb]{\smash{{\SetFigFont{11}{13.2}{\rmdefault}{\mddefault}{\updefault}{\color[rgb]{0,0,0}$k_{\mb{s}}$}%
}}}}
\put(2306,499){\makebox(0,0)[lb]{\smash{{\SetFigFont{11}{13.2}{\rmdefault}{\mddefault}{\updefault}{\color[rgb]{0,0,0}$(0,0,0)$}%
}}}}
\put(3871,-241){\makebox(0,0)[lb]{\smash{{\SetFigFont{11}{13.2}{\rmdefault}{\mddefault}{\updefault}{\color[rgb]{1,0,0}$k_2=k_1-v(\mb{s}_*)-n$}%
}}}}
\put(3961,-691){\makebox(0,0)[lb]{\smash{{\SetFigFont{11}{13.2}{\rmdefault}{\mddefault}{\updefault}{\color[rgb]{0,0,0}$k_2=k_1-v(\mb{s}_*)$}%
}}}}
\put(631,-1591){\makebox(0,0)[lb]{\smash{{\SetFigFont{11}{13.2}{\rmdefault}{\mddefault}{\updefault}{\color[rgb]{0,0,0}$k_2=0$}%
}}}}
\put(1216,-916){\makebox(0,0)[lb]{\smash{{\SetFigFont{11}{13.2}{\rmdefault}{\mddefault}{\updefault}{\color[rgb]{0,0,0}$\loz_{\mb{s}}$}%
}}}}
\put(2701, 74){\makebox(0,0)[lb]{\smash{{\SetFigFont{11}{13.2}{\rmdefault}{\mddefault}{\updefault}{\color[rgb]{1,0,0}$\loz_{\mb{s}}-\mb{w}$}%
}}}}
\put(2881,-1546){\makebox(0,0)[lb]{\smash{{\SetFigFont{11}{13.2}{\rmdefault}{\mddefault}{\updefault}{\color[rgb]{0,.56,0}$\loz_{\mb{s}}+\mb{w}$}%
}}}}
\end{picture}%

\end{center}
By construction, the boundary of $\loz_{\mb{s}}$ is contained in the
hyperplanes with equations $k_i=0$ and $k_j=k_1-v(\mb{s}_*)$ for
$i,j\in\llbracket2, n\rrbracket$. Hence $\loz_{\mb{s}}\ssm
(\loz_{\mb{s}}+\mb{w})= \bigcup_{i\in\llbracket2,n\rrbracket}\{\mb{k}\in
\loz_{\mb{s}} : 0\leq k_i<1\}$ (see the above picture). By Proposition
\ref{prop:descripxt} \eqref{item1_5:descripxt} in dimension $n-1$, for
every $i\in\llbracket2,n\rrbracket$, we have $\card\{\mb{k}\in
\loz_{\mb{s}}: k_i=0\}=\bigO((-v(\mb{s}_*))^{n-2})$. Therefore
\[\card\big(\loz_{\mb{s}}\ssm
(\loz_{\mb{s}}+\mb{w})\big)=\bigO((-v(\mb{s}_*))^{n-2})\,.
\]
Similarly, we have $\loz_{\mb{s}}\ssm (\loz_{\mb{s}}-\mb{w})=
\bigcup_{i\in\llbracket2,n\rrbracket}\{\mb{k}\in \loz_{\mb{s}} :k_1-
v(\mb{s}_*)-n<k_i\leq k_1-v(\mb{s}_*)\}$ (see the above picture) and
$\card\big((\loz_{\mb{s}}+\mb{w})\ssm \loz_{\mb{s}}\big)= \card(
\loz_{\mb{s}}\ssm (\loz_{\mb{s}}-\mb{w}))= \bigO((-v(\mb{s}_*))^{n-2}
)$.  Therefore by Proposition \ref{prop:descripxt}
\eqref{item1_5:descripxt}, the cardinality of the symmetric difference
between $\loz_{\mb{s}}$ and $\loz_{\mb{s}}+\mb{w}$ is negligible with
respect to the cardinality of $\loz_{\mb{s}}$.

This implies the weak-star convergence $\nu_{\mb{s}} - \aaa_*
\nu_{\mb{s}}\; \weakstar\; 0$ as $\mb{s} \to +\infty$. Finally, since
the transformation $\aaa:\X \to \X$ is a homeomorphism (in particular,
it is continuous and proper), we have $\aaa_* \nu = \nu$.
\cqfd

\medskip Let us recall the notation $\ell_m\in\NN\ssm\{0\}$ introduced
in Lemma \ref{lem:cover_X_1_dynamic_balls} for every $m\in\NN$ and
$\kappa'(\mb{s})$ introduced in Equation \eqref{eq:defikappaprims} for
every $\mb{s} \in(R_v\ssm \{0\})^{n-1}$ satisfying Equation
\eqref{eq:proprisbold}. Let us recall the notation $\P^N =
\bigvee_{i=0}^{N-1}\aaa^{-i} \P$ introduced in Equation
\eqref{eq:defPN} for every $N\in\NN\ssm\{0\}$ and every finite
measurable partition $\P$ of $\X$.

\blemm\label{lem:minoentronuskN}
%%
%\todo{\tiny analogue Lem 3.10 \cite{DavSha20}}
%%
Assume that $R_v$ is principal. For every $\eta\in\;]0,1[\,$, there
exists $m=m(\eta)\in\NN$ such that with $\ell= \max\{\ell_m,m+1\}$,
for every $(m,\ell)$-partition $\P$ of $\X$ and for every
$M\in\NN\ssm\{0\}$, there exists $N_0=N_0(\eta,\P,M)\in\NN\ssm\{0\}$
such that for every $\mb{s} \in (R_v\ssm \{0\})^{n-1}$ satisfying
Equations \eqref{eq:proprisbold} and \eqref{eq:unifbehavsi}, for every
$\mb{k}= (k_1,\ldots, k_n) \in \loz_{\mb{s}}$ and for every $N\in\NN$
satisfying the three assumptions 
\[
\mb{k}+(N-1)\mb{w}\in\loz_{\mb{s}}\,,
\]
\[
\max\Big\{N_0,\,\frac{4(1-\eta)}{\eta}(n+1)(c_0+1)
\frac{-v(\mb{s}_*)}{\max\{1,\, \ln(-v(\mb{s}_*))\}},\,
\frac{\kappa'(\mb{s})}{c_2\, q_v^{-mn}}\Big\}\leq N\]
\[
\text{and}\quad N\leq 
\frac{-v(\mb{s}_*) -\ell+k_1- \max_{2 \leq i \leq n} k_i}{n}\,,
\]
we have
\[
\frac{1}{M}\,H_{\nu_{\mb{s},[\mb{k},N]}}(\P^M)\geq(1-\eta)^2\;n(n-1)-\eta\,.
\]
\elemm

\dem For every $\eta\in\;]0,1[$, let $\kappa=\frac{n(n-1)\,\eta}{n^3}
\in\;]0,1[$, let $m=\big\lceil-\frac{1}{n}\log_{q_v}
\frac{\eta^2n(n-1)} {2c_2n^3} \big\rceil$ which belongs to $\NN$ since
$\eta^2< 1\leq 2c_2$, and let $\epsilon= q_v^{-m}\in
q_v^\ZZ\cap\,]0,1[\,$. Let $\ell,\P,M$ be as in the statement.  With
$c_{\varphi_v}$ the constant introduced in Equation
    \eqref{eq:minovarphiv}, let
\[
N_0=\max\Big\{M,\frac{2M}{\eta}\,\log_{q_v}\card\, \P,
\frac{4(1-\eta)}{\eta}(n-1)(n+1-\log_{q_v}c_{\varphi_v})\Big\}
\in\NN\ssm\{0\}
\]
and let $N\in\NN$ with $N\geq N_0$. Let $\mb{s}, \mb{k}$ be as in the
statement.  Let
\[
\nu_{\mb{s},\mb{k}}=\frac{1}{\card\; \Lambda_{\mb{s}}}\;
\sum_{\mb{t}\in\Lambda_{\mb{s}}}\;\delta_{\exp({\mb{k}})x_{\mb{t}}}\,.
\]
Note that by the definition of $\nu_{\mb{s},[\mb{k},N]}$ in Equation
\eqref{eq:definuskN}, since $\aaa=\exp(\mb{w})$ and by the definition
in Equation \eqref{eq:defSNmu} of the $N$-th Birkhoff average of
measures for $\aaa$, we have
\[
\nu_{\mb{s},[\mb{k},N]}=\frac{1}{N} \,\sum_{i=0}^{N-1}
(\aaa^i)_*\nu_{\mb{s},\mb{k}}=S_N\nu_{\mb{s},\mb{k}}\,.
\]
By Lemma \ref{lem:entromino} \eqref{item1:entromino} applied with
$\mu'=\nu_{\mb{s},\mb{k}}$ and $\phi=\aaa$ since $N\geq N_0\geq M\geq
1$, we have
\begin{equation}\label{eq:ramifpullback}
\frac{1}{M} H_{\nu_{\mb{s},[\mb{k},N]}}(\P^M) \geq \frac{1}{N}
H_{\nu_{\mb{s},\mb{k}}}(\P^N)- \frac{M}{N} \log_{q_v}\card\;\P\,.
\end{equation}
Since $N\geq N_0$, we have $\frac{M}{N} \log_{q_v}\card\;\P\leq
\frac{\eta}{2}$.

As in Lemma \ref{lem:cover_X_1_dynamic_balls}
\eqref{item1:cover_X_1_dynamic_balls}, let $\X'=\X'_{\P,\,\kappa,\,N}$
be a measurable subset of $\X^{\geq q_v^{-m}}$, let $\P'$ be a subset
of the partition $\P^N$ and for every $P \in \P'$, let $F_P$ be a
finite subset of $P$ with cardinality at most $q_v^{\; n^3\kappa N}
=q_v^{\; n(n-1)\eta N}$ such that $\X' = \bigcup \P'$ and $P\subset
\bigcup_{x\in F_P} B_{\ell,N-1}\, x$. Since $F_P \subset P\subset \X'
\subset \X^{\geq q_v^{-m}}$, for every $x\in F_P$, we have $\sys(x)
\geq q_v^{-m}$. Hence $\max\{0,-\log_{q_v} (\sys(x))\} \leq m < \ell$
by the definition of $\ell$. Therefore by the assumptions on $\mb{s}$
and $\mb{k}$, by Lemma \ref{lem:counting_latpts_dynballs} and by
Equation \eqref{eq:cardLambdasss} on the left, for every $x\in F_P$,
we have
\begin{align*}
\nu_{\mb{s},\mb{k}}(B_{\ell, N-1} \;x)&=
\frac{1}{\card(\Lambda_{\mb{s}})}\;\card\big( \{\mb{t} + R_v^{\; n-1}
\in \Lambda_{\mb{s}}: \exp(\mb{k}) \,x_{\mb{t}} 
\in B_{\ell, N-1} \;x\} \big)\nonumber\\&
\leq \frac{2^{n-1}}{\prod_{i=2}^n\;\varphi_v(s_i)}\;
q_v^{\,-\ell(n-1)-v(\mb{s}_*)(n-1)-n(n-1)(N-1)}\,.%\label{eq:uppercontrolnusN}
\end{align*}
Thus, since $\P'\subset \P^N$, since we have $P\subset\bigcup_{x\in
  F_P} B_{\ell,N-1}\, x$ and $\card\; F_P\leq q_v^{\; n(n-1)\,\eta N}$ for
every $P\in \P'$, and since $\X'=\bigsqcup_{P\in\P'}P$, we have
\begin{align}
  &H_{\nu_{\mb{s},\mb{k}}}(\P^N)=
  - \sum_{P\in \P^N} \nu_{\mb{s},\mb{k}}(P) \log_{q_v} \nu_{\mb{s},\mb{k}}(P)\geq
  - \sum_{P\in \P'} \nu_{\mb{s},\mb{k}}(P) \log_{q_v} \nu_{\mb{s},\mb{k}}(P)
  \nonumber\\\geq\;&
  - \sum_{P\in \P'} \nu_{\mb{s},\mb{k}}(P)\log_{q_v}\Big(q_v^{\; n(n-1)\,\eta N}\;
  \frac{2^{n-1}}{\prod_{i=2}^n\;\varphi_v(s_i)}\;
  q_v^{\,(n-\ell)\ell(n-1)-v(\mb{s}_*)(n-1)-n(n-1)N}\Big)\nonumber
  \\=\;&\nu_{\mb{s},\mb{k}}(\X')\Big((1-\eta)n(n-1)N-\frac{\ln 2}{\ln q_v}(n-1)-
  (n-\ell)(n-1)+ \sum_{i=2}^n\log_{q_v}\frac{\varphi_v(s_i)}{|\mb{s}_*|}\Big)\,.
  \label{eq:calcpartielHnuPM}
\end{align}
By Equation \eqref{eq:minovarphiv}, by computations similar to the
ones done in the proof of Remark \ref{rem:controlkappaprims} and by
the assumptions on $\mb{s}$ in the statement of Lemma
\ref{lem:minoentronuskN}, we have
\begin{align*}
  &\sum_{i=2}^n\log_{q_v}\frac{\varphi_v(s_i)}{|\mb{s}_*|}\geq\sum_{i=2}^n
  \log_{q_v}\frac{c_{\varphi_v}\,|s_i|}{\max\{1,\ln(-v(s_i))\}\,|\mb{s}_*|}
  \\=\;&\sum_{i=2}^n - \log_{q_v}\big(\max\{1,\ln(-v(s_i))\}\big)
  +\log_{q_v}c_{\varphi_v}-v(s_i)+v(\mb{s}_*)
  \\\geq\;& -(n-1)\Big(\log_{q_v}(\max\{1,\ln(-v(\mb{s}_*))\})
  -\log_{q_v}c_{\varphi_v}+c_0\frac{-v(\mb{s}_*)}{\max\{1,\,
    \ln(-v(\mb{s}_*))\}}\Big)
  \\\geq\;& -(n-1)\Big((c_0+1)\frac{-v(\mb{s}_*)}{\max\{1,\,
    \ln(-v(\mb{s}_*))\}}-\log_{q_v}c_{\varphi_v}\Big)
  \\\geq\;& -\frac{\eta}{4(1-\eta)}N+(n-1)\log_{q_v}c_{\varphi_v}\,.
\end{align*}
Since $N\geq N_0$, we have
\begin{align*}
&\frac{\ln 2}{\ln q_v}(n-1)+(n-\ell)(n-1)-(n-1)
\log_{q_v}c_{\varphi_v}\\\leq\;& (n-1)(n+1-\log_{q_v}c_{\varphi_v})
\leq \frac{\eta}{4(1-\eta)}\,N\,.
\end{align*}
Therefore Equation \eqref{eq:calcpartielHnuPM} becomes
\[
\frac{1}{N}\,H_{\nu_{\mb{s},\mb{k}}}(\P^N)\geq
\nu_{\mb{s},\mb{k}}(\X')\Big((1-\eta)n(n-1)-\frac{\eta}{2(1-\eta)}\Big)\,.
\]

We have $N\geq \frac{\kappa'(\mb{s})}{c_2\, q_v^{-mn}}=
\frac{\kappa'(\mb{s})}{c_2 \,\epsilon^{n}}$ by the assumptions on
$\mb{s}$ and the definition of $\epsilon$ at the beginning of this
proof. By Lemma \ref{lem:cover_X_1_dynamic_balls}
\eqref{item2:cover_X_1_dynamic_balls} applied with
$\mu=\nu_{\mb{s},\mb{k}}$, by Corollary \ref{coro:decaythin2} proving
that there is no escape of mass (the only place in this proof that
requires the principal assumption on $R_v$) whose assumptions
\eqref{eq:hyplemfuitemasse2} are satisfied, and by the definitions of
$\kappa$ and $m=m(\eta)$ at the beginning of this proof, we have
\[
\nu_{\mb{s},\mb{k}}(\X')\geq 1-\frac{1}{\kappa} \nu_{\mb{s},[\mb{k},N]}
(\X^{<q_v^{-m}})\geq 1-\frac{2\,c_2\,\epsilon^{n}}{\kappa} =
1-\frac{2\,c_2\,q_v^{-mn}n^3}{n(n-1)\,\eta} \;\geq 1-\eta\,.
\]
Hence Lemma \ref{lem:minoentronuskN} follows using Equation
\eqref{eq:ramifpullback}. \cqfd

\medskip
\noindent{\bf End of the proof of Theorem \ref{theo:main}. }
Let us fix a sequence $(\mb{s}^{(j)})_{j\in\NN}$ of elements of $(R_v\ssm
\{0\})^{n-1}$ satisfying Equations \eqref{eq:proprisbold} and
\eqref{eq:unifbehavsi} and tending to infinity such that we have
$\nu={\displaystyle\lim_{j\ra+\infty}}\;\nu_{\mb{s}^{(j)}}$.

Let us fix $\eta>0$, that will tend to $0$ at the very end of the
proof.  Let $m=m(\eta)$, $\ell=\max\{\ell_m,m+1\}$, $\P$ a
$(m,\ell)$-partition of $\X$, $M\in\NN\ssm\{0\}$ and $N_0=
N_0(\eta,\P,M)$ be as in Lemma \ref{lem:minoentronuskN}.  Since
$h_{\nu}(\aaa)$ is the upper bound of $h_{\nu}(\aaa,\P')$ where $\P'$
varies over all finite measurable partitions of $\X$, and by Equation
\eqref{eq:defientrorelpart} applied with $\mu'=\nu$ and $\phi=\aaa$,
if $M$ is large enough, we have
\[
h_{\nu}(\aaa)\geq h_{\nu}(\aaa,\P)\geq
\frac{1}{M}\,H_{\nu}(\P^M)-\eta\,.
\]

For all $\mb{s}\in(R_v\ssm \{0\})^{n-1}$ satisfying Equations
\eqref{eq:proprisbold} and \eqref{eq:unifbehavsi}, and
$\mb{k}=(k_1,\dots, k_n)\in\loz_{\mb{s}}$, let $N_{\mb{k}}=
\max\{\ell'\in\NN\ssm\{0\}: \mb{k}+ (\ell'-1) \mb{w}\in \loz_{\mb{s}}
\}$. For every $\ell''\in\NN$, the point $\mb{k}+ \ell'' \mb{w}$
belongs to $\loz_{\mb{s}}$ if and only if, for every $i\in \llbracket
2,n \rrbracket$, we have $0\leq k_i+\ell''\leq k_1+ (1-n) \ell''-
v(\mb{s}_*)$. Hence
\[
N_{\mb{k}}-1=\Big\lfloor\frac{1}{n}\big(-
v(\mb{s}_*)+k_1-\max_{i\in\llbracket2,n
\rrbracket}k_i\big)\Big\rfloor\,.
\]
Let
\[
\partial^+\loz_{\mb{s}}=\{\mb{k}=(k_1,\dots,k_n)\in\loz_{\mb{s}}:
\min_{i\in \llbracket2, n\rrbracket} k_i=0\}
\]
be the upper part of the boundary of $\loz_{\mb{s}}$ (see the picture
above Remark \ref{rem:controlkappaprims}). Since $\loz_{\mb{s}}$ is
the disjoint union of the maximal vertical (directed by $\mb{w}$)
segments contained in it, we have $\loz_{\mb{s}}=
\bigsqcup_{\mb{k}\in\partial^+\loz_{\mb{s}}}\{\mb{k}+ \ell'\mb{w}:
\ell'\in\llbracket0, N_{\mb{k}}-1\rrbracket\}$. Let
\[
N'_{\mb{k}}=\Big\lfloor\frac{1}{n}\big(- v(\mb{s}_*)-\ell+ k_1-
\max_{i\in\llbracket2,n \rrbracket} k_i\big)\Big\rfloor+1\,.
\]
Note that $0\leq N_{\mb{k}}-N'_{\mb{k}}\leq\lceil\frac{\ell}{n}\rceil$
is uniformly bounded for $\eta$ fixed and that $N'_{\mb{k}}$ satisfies
the upper bound assumption on $N$ in Lemma \ref{lem:minoentronuskN}.
With $c_{\varpi_v}$ the constant introduced in Equation
\eqref{eq:majoomegav}, let $c_\eta= \max\big\{\frac{4(1-\eta)}{\eta}
(n+1)(c_0+1),\; \frac{c_{\varpi_v} +c_0+1}{n\, c_2\,q_v^{-m(\eta)\,n}}
\big\}$ and
\[
\Omega_{\mb{s}}=\Big\{\mb{k}\in\partial^+\loz_{\mb{s}}:
N'_{\mb{k}}\geq\max\big\{N_0, c_\eta \;
\frac{-v(\mb{s}_*)}{\max\{1,\,\ln(-v(\mb{s}_*))\}}\big\}\Big\}\,.
\]
Note that for every $\mb{k}\in\Omega_{\mb{s}}$, by Remark
\ref{rem:controlkappaprims} whose assumption holds true since $\mb{s}$
verifies Equation \eqref{eq:unifbehavsi}, the number $N'_{\mb{k}}$
satisfies the lower bound assumption on $N$ in Lemma
\ref{lem:minoentronuskN}.
\begin{center}
  \!\!\!\!\!\!\!\!\!\!\!\!\!\!\!\!\!\!\!\!\!\!\!\!\!\!\!\!
  \!\!\!\!\!\!\!\!\!\!\!\!\!\!\!\!\!\!
  \begin{picture}(0,0)%
\includegraphics{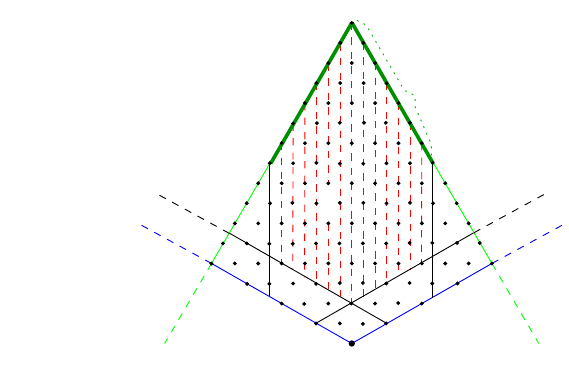}%
\end{picture}%
\setlength{\unitlength}{3812sp}%
\begingroup\makeatletter\ifx\SetFigFont\undefined%
\gdef\SetFigFont#1#2#3#4#5{%
  \reset@font\fontsize{#1}{#2pt}%
  \fontfamily{#3}\fontseries{#4}\fontshape{#5}%
  \selectfont}%
\fi\endgroup%
\begin{picture}(4695,3054)(-1546,-298)
\put(2960,-143){\makebox(0,0)[lb]{\smash{{\SetFigFont{11}{13.2}{\rmdefault}{\mddefault}{\updefault}{\color[rgb]{0,.69,0}$k_3=0$}%
}}}}
\put(-142,-136){\makebox(0,0)[lb]{\smash{{\SetFigFont{11}{13.2}{\rmdefault}{\mddefault}{\updefault}{\color[rgb]{0,.69,0}$k_2=0$}%
}}}}
\put(3129,780){\makebox(0,0)[lb]{\smash{{\SetFigFont{11}{13.2}{\rmdefault}{\mddefault}{\updefault}{\color[rgb]{0,0,1}$k_2=k_1-v(\mb{s_*})$}%
}}}}
\put(-1525,715){\makebox(0,0)[lb]{\smash{{\SetFigFont{11}{13.2}{\rmdefault}{\mddefault}{\updefault}{\color[rgb]{0,0,1}$k_3=k_1-v(\mb{s_*})$}%
}}}}
\put(2987,1187){\makebox(0,0)[lb]{\smash{{\SetFigFont{11}{13.2}{\rmdefault}{\mddefault}{\updefault}{\color[rgb]{0,0,0}$k_2=k_1-v(\mb{s_*})-\ell$}%
}}}}
\put(-1531,1195){\makebox(0,0)[lb]{\smash{{\SetFigFont{11}{13.2}{\rmdefault}{\mddefault}{\updefault}{\color[rgb]{0,0,0}$k_3=k_1-v(\mb{s_*})-\ell$}%
}}}}
\put(1936,1964){\makebox(0,0)[lb]{\smash{{\SetFigFont{11}{13.2}{\rmdefault}{\mddefault}{\updefault}{\color[rgb]{0,.69,0}$\Omega_{\mb{s}}$}%
}}}}
\put(1068,989){\makebox(0,0)[lb]{\smash{{\SetFigFont{11}{13.2}{\rmdefault}{\mddefault}{\updefault}{\color[rgb]{1,0,0}$\loz'_{\mb{s}}$}%
}}}}
\put(1398,-234){\makebox(0,0)[lb]{\smash{{\SetFigFont{11}{13.2}{\rmdefault}{\mddefault}{\updefault}{\color[rgb]{0,0,0}$k_{\mb{s}}$}%
}}}}
\put(1374,2609){\makebox(0,0)[lb]{\smash{{\SetFigFont{11}{13.2}{\rmdefault}{\mddefault}{\updefault}{\color[rgb]{0,0,0}$(0,0,0)$}%
}}}}
\put(1505,2391){\makebox(0,0)[lb]{\smash{{\SetFigFont{11}{13.2}{\rmdefault}{\mddefault}{\updefault}{\color[rgb]{0,0,0}$(-1,1,0)$}%
}}}}
\put(456,2384){\makebox(0,0)[lb]{\smash{{\SetFigFont{11}{13.2}{\rmdefault}{\mddefault}{\updefault}{\color[rgb]{0,0,0}$(-1,0,1)$}%
}}}}
\end{picture}%

\end{center}
Let $\loz'_{\mb{s}}= \bigsqcup_{\mb{k}\in\Omega_{\mb{s}}} \big\{\mb{k}
+\ell' \mb{w}: \ell'\in\llbracket0,N'_{\mb{k}}-1 \rrbracket$, which is
obtained from $\loz_{\mb{s}}$ by removing a bounded size neighborhood
of the lower part of the boundary of $\loz_{\mb{s}}$ and a
comparatively small part of the vertical side of $\loz_{\mb{s}}$ (see
the above picture). More precisely, $\loz'_{\mb{s}}=\loz_{\mb{s}}\ssm
\big(\loz''_{\mb{s}}\cup \loz'''_{\mb{s}}\big)$ where $\loz''_{\mb{s}}
=\bigsqcup_{\mb{k}\in\Omega_{\mb{s}}} \big\{\mb{k} +\ell' \mb{w}:\ell'
\in \llbracket N'_{\mb{k}},N_{\mb{k}}-1\rrbracket$, whose cardinality
is $\bigO((-v(\mb{s}_*))^{n-2})$ as seen in the proof of Lemma
\ref{lem:nuaaainvar}, and $\loz'''_{\mb{s}} =\bigsqcup_{\mb{k}\in
  \partial^+\loz_{\mb{s}}\,\ssm\,\Omega_{\mb{s}}}\big\{\mb{k}+\ell'\mb{w}:
\ell'\in\llbracket 0,N_{\mb{k}}-1\rrbracket$. Since we have
$N_{\mb{k}}= \bigO \big(\frac{-v(\mb{s}_*)}{\max\{1,\, \ln(
  -v(\mb{s}_*))\}}\big)$ when $\mb{k}\in \partial^+ \loz_{\mb{s}}
\ssm\,\Omega_{\mb{s}}$, and since $\card\;(\partial^+
\loz_{\mb{s}})=\bigO((-v(\mb{s}_*))^{n-2})$ as seen in the proof of
Lemma \ref{lem:nuaaainvar}, the cardinality of $\loz'''_{\mb{s}}$ is
$\bigO\big(\frac{(-v(\mb{s}_*) )^{n-1}}{\max\{1,\,\ln(-v(\mb{s}_*))
  \}} \big)$. Modifying Equation \eqref{eq:definuloz}, let
\[%\begin{equation}\label{eq:definuloprim}
\nu_{\mb{s}}^{\loz'}=
\frac{1}{\card\;\Lambda_{\mb{s}}\;\card\;\loz'_{\mb{s}}}\;
\sum_{\mb{t}\in\Lambda_{\mb{s}},\;\mb{k}\in\loz'_{\mb{s}}}\;\int_{a\in A(\OOO_v)}
\;\delta_{a\exp(\mb{k})\,x_{\mb{t}}}\;da\;.
\]%\end{equation}
Since the cardinalities of $\loz''_{\mb{s}}$ and $\loz'''_{\mb{s}}$ are
negligible compared to the one of $\loz_{\mb{s}}$ (given by
Proposition \ref{prop:descripxt} \eqref{item1_5:descripxt}) we have
$\displaystyle{\lim_{j\ra+\infty}}\; \nu_{\mb{s}^{(j)}}^{\loz'} =
\displaystyle{\lim_{j\ra+\infty}}\; \nu_{\mb{s}^{(j)}}^{\loz}=\nu$.
In particular, for every $j\in\NN$ large enough, we have
\[
h_{\nu}(\aaa)\geq 
\frac{1}{M}\,H_{\nu_{\mb{s}^{(j)}}^{\loz'}}(\P^M)-2\eta\,.
\]
Let $\omega_{\mb{s}}$ be the probability measure on the finite
(discrete) set $\Omega_{\mb{s}}$ defined by $\omega_{\mb{s}}(\mb{k})=
\frac{N'_{\mb{k}}} {\card\; \loz'_{\mb{s}}}$ for every $\mb{k}\in
\Omega_{\mb{s}}$.  Then by Equation \eqref{eq:definuskN}, we have
\[
\nu_{\mb{s}}^{\loz'}=\int_{a\in A(\OOO_v)}\int_{\mb{k}\in\Omega_{\mb{s}}}
a_*\nu_{\mb{s},[\mb{k},N'_{\mb{k}}]}\; d\omega_{\mb{s}} (\mb{k})\,da\,.
\]
By Lemma \ref{lem:entromino} \eqref{item2:entromino} applied with
$(\Omega,\omega)=(A(\OOO_v)\times\Omega_{\mb{s}},da\otimes\omega_{\mb{s}})$,
since $a^{-1}\P$ is also an \mbox{$(m,\ell)$-partition} for every
$a\in A(\OOO_v)$, and by Lemma \ref{lem:minoentronuskN} applied with
$N=N'_{\mb{k}}$ and integrated over $(a,\mb{k})\in
A(\OOO_v)\times\Omega_{\mb{s}}$ for the probability measure
$da\otimes\omega_{\mb{s}}$, we have
\[
\frac{1}{M}\,H_{\nu_{\mb{s}}^{\loz'}}(\P^M)\geq
\int_{a\in A(\OOO_v)}\int_{\mb{k}\in\Omega_{\mb{s}}}
\frac{1}{M}\,H_{a_*\nu_{\mb{s},[\mb{k},N'_{\mb{k}}]}}(\P^M)\;
d\omega_{\mb{s}}(\mb{k})\,da \geq
%\int_{\Omega_{\mb{s}}}\big((1-\eta)^2\;n(n-1)-\eta\big)\;
%d\omega_{\mb{s}}=
(1-\eta)^2\;n(n-1)-\eta\,.
\]
Thus $h_{\nu}(\aaa)\geq (1-\eta)^2\;n(n-1)-3\,\eta$. By letting
$\eta\ra0$, we have $h_\nu(\aaa) \geq n(n-1)$. By the
Einsiedler-Lindenstrauss Theorem \ref{theo:EL}, we hence have
$h_\nu(\aaa) = n(n-1)$ and then $\nu=\frac{m_{\X}}{\|m_{\X}\|}$, as
wanted at the beginning of the proof of Theorem \ref{theo:main}.
\cqfd

\medskip
The following result follows by averaging Theorem \ref{theo:main} over
the permutations of $\llbracket 2,n\rrbracket$ and over the compact
probability space $(\OOO_v^{\,\times}/R_v^{\,\times},
\frac{q_v(q+1)}{q_v-1} \vol'_v)$ as in the proof of Lemma
\ref{lem:deducmainuintro2}.

\bcoro\label{coro:mainsgras}
%%
%\todo{\tiny à voir}
%%
Assume that $R_v$ is principal. For every $\mb{s}$ in the set $S_n$
(endowed with the Fréchet filter) of elements $(s_2,\ldots s_n)\in
(R_v\ssm\{0\})^{n-1}$ with $s_2\mid s_3\mid\ldots\mid s_n$, $v(s_n)\in
n\ZZ$ and $v(s_2)-v(s_n)\leq \frac{-v(s_n)}{\max\{1,\, \ln(-v(s_n))}$,
let us define
\[
\Lambda'_{\mb{s}}=
\left\{\Big(\frac{r_2}{s'_2},\ldots,\frac{r_n}{s'_n}\Big)\mod R_v^{\;n-1}:
\begin{array}{l}r_2,\ldots, r_n, s'_2,\ldots,s'_n\in R_v,\\
  \forall\;j\in\llbracket2,n\rrbracket,\;r_jR_v+s'_jR_v=R_v,\\
  \{s'_2,\ldots,s'_n\}=\{s_2,\ldots,s_n\}\,.\end{array}\right\}
\]
For the weak-star convergence of Radon measures on the locally compact
space $\X_1$, we have
\[
\lim_{\mb{s}\in S_n,\;|s_n|\ra+\infty}\;\;
\frac{1}{\card\;\Lambda'_{\mb{s}}}\;\;\sum_{\mb{t}\,\in\,\Lambda'_{\mb{s}}}\;
\ov\mu_{\uuu_{\mb{t}} R_v^{\;n}}\;=\; \frac{{\tt
m}_{\X_1}}{\|{\tt m}_{\X_1}\|}\,.\;\;\;\Box
\]
\ecoro

{\small \bibliography{../biblio} }

\begin{thebibliography}{BAPP}


\bibitem[AK]{AraKim25}
N.~S.~Aranov and T.~Kim.
\newblock {\it Hausdorff dimension of singular vectors in function fields}.
\newblock {Adv. in Math. {\bf 461} (2025) 110084}.

\bibitem[BKL]{BanKimLim25}
G.~Bang, T.~Kim, S.~Lim.
\newblock {\it Singular linear forms over global function fields}.
\newblock {Preprint {\tt [arXiv:2404.07752]}}.

\bibitem[BoT]{BorTit65}
A.~Borel and J.~Tits.
\newblock {\it Groupes réductifs}.
\newblock {Publ. Math. I.~H.~É.~S. {\bf 27} (1965)  55--150}.

\bibitem[Bow]{Bowen72b}
R.~Bowen.
\newblock {\it The equidistribution of closed geodesics}.
\newblock {Amer. J. Math. {\bf 94}  (1972) 413--423}.
  
\bibitem[BPP]{BroParPau19}
A.~Broise-Alamichel, J.~Parkkonen and F.~Paulin.
\newblock {\it Equidistribution and counting under equilibrium states in
  negative curvature and trees. Applications to non-Archimedean Diophantine
  approximation}.
\newblock {With an Appendix by J.~Buzzi. Prog. Math. {\bf 329}, Birkhäuser,
  2019}.

\bibitem[CG]{CheChe16}
Y.~Cheung and N.~Chevallier.
\newblock {\it Hausdorff dimension of singular vectors}.
\newblock {Duke Math. J. {\bf 165} (2016), 2273--2329}.

\bibitem[DaL]{DanLi22}
N.-T.~Dang and J.~Li.
\newblock {\it Equidistribution and counting of periodic tori
  in the space of Weyl chambers}.
\newblock {Preprint {\tt [arXiv:2202.08323]}, to appear in
  Commentarii Math. Helv.}.

\bibitem[DFSU]{DasFisSimUrb24}
T.~Das, L.~Fishman, D.~Simmons and M.~Urbański.
\newblock {\it A variational principle in the parametric geometry of numbers}.
\newblock {Adv. Math. {\bf 437} (2024) 109435}.

\bibitem[DKMS]{DavKimMorSha25}
O.~David, T.~Kim, R.~Mor and U.~Shapira.
\newblock {\it On the rate of convergence of continued fraction statistics
  of random rationals}.
\newblock {Selecta Math. {\bf 31} (2025) n$^0$ 33}.
%Preprint {\tt [arXiv:2401.15586]}
      
\bibitem[DS1]{DavSha18}
O.~David and U.~Shapira.
\newblock {\it Equidistribution of divergent orbits and continued
fraction expansion of rationals}.
\newblock {J. London Math. Soc. {\bf 98} (2018) 149--176}.

\bibitem[DS2]{DavSha20}
O.~David and U.~Shapira.
\newblock {\it Equidistribution of divergent orbits of the diagonal
group in the space of lattices}.
\newblock {Erg. Theo. Dyn. Syst. {\bf 40} (2020) 1217--1237}.

\bibitem[EL]{EinLin10}
M.~Einsiedler and E.~Lindenstrauss.
\newblock {\it Diagonal actions on locally homogeneous spaces}.
\newblock {In "Homogeneous flows, moduli spaces and arithmetic", 
pp 155--241, Clay Math. Proc. {\bf 10}, Amer. Math. Soc. 2010}.

\bibitem[ELMV1]{EinLinMicVen11}
M.~Einsiedler, E.~Lindenstrauss, P.~Michel and A.~Venkatesh.
\newblock {\it Distribution of periodic torus
orbits and Duke's theorem for cubic fields}.
\newblock {Ann. of Math. {\bf 173} (2011) 815--885}.

\bibitem[ELMV2]{EinLinMicVen12}
M.~Einsiedler, E.~Lindenstrauss, P.~Michel and A.~Venkatesh.
\newblock {\it The distribution of closed
geodesics on the modular surface, and Duke's theorem}.
\newblock {L'Enseign. Math. {\bf 58} (2012) 249--313}.

\bibitem[ELW]{EinLinWar22}
M.~Einsiedler, E.~Lindenstrauss and T.~Ward.
\newblock {\it Entropy in ergodic theory and topological dynamics}.
\newblock {Book in preparation, https:/\!/tbward0.wixsite.com/books/entropy}.

\bibitem[Gos]{Goss98}
D.~Goss.
\newblock {\it Basic structures of function field arithmetic}.
\newblock {Erg. Math. Grenz. {\bf 35}, Springer Verlag, 1998}.

\bibitem[HaW]{HarWri08}
G.~H.~Hardy and E.~M.~Wright.
\newblock {\it An introduction to the theory of numbers}.
\newblock {Oxford Univ. Press, sixth ed., 2008}.
  
\bibitem[HoP]{HorPau24}
T.~Horesh and F.~Paulin.
\newblock {\it Joint effective equidistribution of partial lattices
in positive characteristic}.
\newblock {Preprint {\tt [arXiv:2404.04368]}}.

\bibitem[KaKLM]{KadKleLinMar17}
S.~Kadyrov, D.~Kleinbock, E.~Lindenstrauss and G.~Margulis.
\newblock {\it Singular systems of linear forms and non-escape
  of mass in the space of lattices}.
\newblock {J. Anal. Math. {\bf 133} (2017) 253--277}.

\bibitem[KePS]{KemPauSha17}
A.~Kemarsky, F.~Paulin and U.~Shapira.
\newblock {\it Escape of mass in homogeneous dynamics in 
positive characteristic}.
\newblock {J. Modern Dyn. {\bf 11} (2017) 369--407}.

\bibitem[KiKL]{KimKimLim21}
T.~Kim, W.~Kim and S.~Lim.
\newblock {\it Dimension estimates for badly approximable affine forms}.
\newblock {Erg. Theo. Dyn. Syst (2024), online}.
%Preprint {\tt [arXiv:2111.15410]}}.

\bibitem[KiLP]{KimLimPau23}
T.~Kim, S.~Lim and F.~Paulin.
\newblock {\it On Hausdorff dimension in inhomogeneous Diophantine
  approximation over global function fields}.
\newblock {J. Numb. Theo. {\bf 251} (2023) 102--146}.
  
\bibitem[KlST]{KleShiTom17}
D.~Kleinbock, R.~Shi and G.~Tomanov.
\newblock {\it S-adic version of Minkovski's geometry of numbers
  and Mahler's compactness criterion}.
\newblock {J. Numb. Theo. {\bf 174} (2017) 150--163}.

\bibitem[KW]{KleWei05}
D.~Kleinbock and B. Weiss.
\newblock {\it Friendly measures, homogeneous flows, and
  singular vectors in algebraic and topological dynamics}.
\newblock {Contemp. Math. {\bf 385} (2005) 281--292}.

\bibitem[Las]{Lasjaunias00}
A.~Lasjaunias.
\newblock {\it A survey of Diophantine approximation in fields of power
  series}.
\newblock {Monat. Math. {\bf 130} (2000) 211--229}.

\bibitem[LSS]{LimSaxSha18}
S.~Lim, N.~de Saxcé and U.~Shapira.
\newblock {\it Dimension bound for badly approximable grids}.
\newblock {Int. Math. Res. Not.  {\bf 20} (2019) 6317--6346}.

\bibitem[MaQ]{MadQue72}
M.~Madan and C.~Queen.
\newblock {\it Algebraic function ﬁelds of class number one}.
\newblock {Acta Arith. {\bf 20} (1972) 423--432}.

\bibitem[Mar]{Margulis69}
G.~Margulis.
\newblock {\it Applications of ergodic theory for the 
investigation of manifolds of negative curvature}.
\newblock {Funct. Anal. Applic. {\bf 3} (1969) 335--336}.

\bibitem[MS]{MerSti15}
P.~Mercuri and C.~Stirpe.
\newblock {\it Classiﬁcation of algebraic function ﬁelds with class
number one}.
\newblock {Journal of Number Theory  {\bf 154} (2015) 365--374}.

\bibitem[Nag]{Nagao59}
H.~Nagao.
\newblock {\it On $\GL(2,K[x])$}.
\newblock {J. Inst. Polytech. Osaka City Univ. Ser. A {\bf 10} (1959)
  117--121}.

\bibitem[Nar]{Narkiewicz04}
W.~Narkiewicz.
\newblock {\it Elementary and analytic theory of algebraic numbers}.
\newblock {3rd Ed., Springer Verlag, 2004}.

\bibitem[PaPS]{ParPauSay25}
J.~Parkkonen, F.~Paulin and R.~Sayous.
\newblock {\it Equidistribution of divergent geodesics
  in negative curvature}.
\newblock {Preprint {\tt [arXiv:2501.03925]}}.
  
\bibitem[Pau]{Paulin02}
F.~Paulin.
\newblock {\it Groupe modulaire, fractions continues et approximation
  diophantienne en caract\'eristique $p$}.
\newblock {Geom. Dedi. {\bf 95} (2002) 65--85}.

\bibitem[PauPS]{PauPolSha15}
F.~Paulin, M.~Pollicott and B.~Schapira.
\newblock {\it Equilibrium states in negative curvature}.
\newblock {Astérisque {\bf 373}, Soc. Math. France, 2015}.

\bibitem[Poo]{Poonen97}
B.~Poonen.
\newblock {\it Torsion in rank $1$ Drinfeld modules and the uniform
boundedness conjecture}.
\newblock {Math. Ann. {\bf 308} (1997) 571–586}.

\bibitem[Ros1]{Rosen99}
M.~Rosen.
\newblock {\it A generalization of Mertens' theorem}.
\newblock {J. Ramanujan Math. Soc. {\bf 14} (1999) 1–19}.

\bibitem[Ros2]{Rosen02}
M.~Rosen.
\newblock {\it Number theory in function fields}.
\newblock {Grad. Texts Math. {\bf 210}, Springer Verlag, 2002}.

\bibitem[Sch]{Schmidt00}
W.~M. Schmidt.
\newblock {\it On continued fractions and diophantine approximation in power
  series fields}.
\newblock {Acta Arith.~{\bf XCV} (2000) 139--166}.

\bibitem[Ser1]{Serre71}
J.-P.~Serre.
\newblock {\it Cohomologie des groupes discrets}.
\newblock {In Ann. Math. Stud. {\bf 70}, pp 77-169.
  Princeton Univ Press, 1971}.

\bibitem[Ser2]{Serre83}
J.-P.~Serre.
\newblock {\it Arbres, amalgames, SL$_2$}.
\newblock {3\`eme \'ed. corr., Ast\'erisque {\bf 46},
  Soc. Math. France, 1983}.

\bibitem[Sha]{Shapira17}
U.~Shapira.
\newblock {\it Full escape of mass for the diagonal group}.
\newblock {Inter. Math. Res. Not. {\bf 15} (2017) 4704--4731}.

\bibitem[ST]{SolTam23}
O.~N.~Solan and N.~Tamam.
\newblock {\it On topologically big divergent trajectories}.
\newblock {Duke Math. J. {\bf 172} (2023) 3429--3474}.

\bibitem[SY]{SolYif23}
O.~N.~Solan and Y.~Yifrach.
\newblock {\it Tori approximation of families of diagonally invariant
  measures}.
\newblock {Geom. Funct. Anal. {\bf 33} (2023) 1354--1378}.

\bibitem[Tam]{Taman22}
N.~Tamam.
\newblock {\it Existence of non-obvious divergent trajectories in
  homogeneous spaces}.
\newblock {Israel J. Math. {\bf 247} (2022) 459--478}.

\bibitem[Tom1]{Tomanov07}
G.~Tomanov.
\newblock {\it Divergent orbits on S-adic homogeneous spaces}.
\newblock {Pure Appl. Math. Q. {\bf 3} (2007) 969--985}.

\bibitem[Tom2]{Tomanov21}
G.~Tomanov.
\newblock {\it Closures of locally divergent orbits of maximal tori
  and values of homogeneous forms}.
\newblock {Erg. Theo. Dyn. Syst. {\bf 41} (2021) 3142--3177}.

\bibitem[TW]{TomWei03}
G.~Tomanov and B.~Weiss.
\newblock {\it Closed orbits for actions of maximal tori on homogeneous
  spaces}.
\newblock {Duke Math. J. {\bf 119} (2003) 367--392}.

\bibitem[Weil]{Weil70}
A.~Weil.
\newblock {\it On the analogue of the modular group in characteristic $p$}.
\newblock {In "Functional Analysis and Related Fields" (Chicago, 1968),
  pp.~211--223, Springer Verlag, 1970}.

\bibitem[Wei1]{Weiss04}
B.~Weiss.
\newblock {\it Divergent trajectories on noncompact parameter spaces}.
\newblock {Geom. Funct. Anal. {\bf 14} (2004) 94--149}.

\bibitem[Wei2]{Weiss06}
B.~Weiss.
\newblock {\it Divergent trajectories and $\QQ$-rank}.
\newblock {Israel J. Math. {\bf 152} (2006), 221--227}.

\end{thebibliography}

\bigskip
{\small
\noindent \begin{tabular}{l}
  N-T.~D, F.~P. \& R.~S. :
  Laboratoire de math\'ematique d'Orsay, UMR 8628 CNRS,\\
Universit\'e Paris-Saclay, 91405 ORSAY Cedex, FRANCE\\
{\it nguyen-thi.dang@universite-paris-saclay.fr,}\\
{\it frederic.paulin@universite-paris-saclay.fr,}\\
{\it rafael.sayous@universite-paris-saclay.fr}
\end{tabular}

\medskip
\noindent \begin{tabular}{l}
  R.~S. : Department of Mathematics and Statistics, P.O. Box 35\\
  40014 University of Jyv\"askyl\"a, FINLAND.\\
  {\it sayousr@jyu.fi}
\end{tabular}
}
\end{document}